\documentclass[trsc,nonblindrev]{informs3noheader}

\OneAndAHalfSpacedXI

 \usepackage[dvipsnames]{xcolor}
\usepackage{booktabs}
\usepackage{eurosym}
\usepackage{mathrsfs}
\usepackage{natbib}
\bibpunct[, ]{(}{)}{,}{a}{}{,}%
\def\BIBand{and}%

 \usepackage[short, nocomma]{optidef}
 \usepackage{graphicx}
\usepackage{subcaption}
 \usepackage{bm}
 \usepackage{multirow}
 \usepackage{hyperref}
\hypersetup{colorlinks,allcolors=blue}
\usepackage{makecell}
\usepackage[para,online,flushleft]{threeparttable}
\usepackage{csquotes}

\usepackage{filecontents}
\usepackage{tikz}
\usepackage{pgfplots}
\pgfplotsset{compat=1.12}

\TheoremsNumberedThrough    
\EquationsNumberedThrough    
\usepackage[linesnumbered,ruled,vlined]{algorithm2e}

\begin{document}
\RUNAUTHOR{Xia et al.} %

\RUNTITLE{Robust charging station location and routing-scheduling for electric modular autonomous units}

\TITLE{Robust charging station location and routing-scheduling for electric modular autonomous units}

\ARTICLEAUTHORS{%
\AUTHOR{Dongyang Xia$^{a}$, Lixing Yang$^{b,*}$, Yahan Lu$^{a}$, Shadi Sharif Azadeh$^{a}$ }
\AFF{$^{a}$ Department of Transport \& Planning, Delft University of Technology, Netherlands; \\
$^{b}$ School of Systems Science, Beijing Jiaotong University, Beijing, 100044, China  \\
*Corresponding author. \EMAIL{E-mail address: lxyang@bjtu.edu.cn (L. Yang)} }

} 

\ABSTRACT{%
\textit{\textbf{Problem definition:}} Motivated by global electrification targets and the advent of electric modular autonomous units (E-MAUs), this paper addresses a robust charging station location and routing-scheduling problem (E-RCRSP) in an inter-modal transit system, presenting a novel solution to traditional electric bus scheduling. The system integrates regular bus services, offering full-line or sectional coverage, and short-turning services. Considering the fast-charging technology with quick top-ups, we jointly optimize charging station locations and capacities, fleet sizing, as well as routing-scheduling for E-MAUs under demand uncertainty. E-MAUs can couple flexibly at different locations, and their routing-scheduling decisions include sequences of services, as well as charging times and locations. \textit{\textbf{Methodology:}} The E-RCRSP is formulated as a path-based robust optimization model, incorporating the polyhedral uncertainty set. We develop a double-decomposition algorithm that combines column-and-constraint generation and column generation armed with a tailored label-correcting approach. To improve computational efficiency and scalability, we propose a novel method that introduces \enquote{super travel arcs} and \enquote{network downsizing} methodologies. \textit{\textbf{Results:}} Computational results from real-life instances, based on operational data of advanced NExT E-MAUs with cutting-edge batteries provided by our industry partner, indicate that charging at both depots and en-route fast-charging stations is necessary during operations. Moreover, our algorithm effectively scales to large-scale operational cases involving entire-day operations, significantly outperforming state-of-the-art methods. Comparisons with fixed-composition buses under the same fleet investment suggest that our methods are able to achieve substantial reductions in passengers' costs by flexibly scheduling units. \textit{\textbf{Implications:}} Our approach is appealing to operators, providing entire-day solutions with enhanced service quality. Its robustness improves the practicality and sustainability of both investments and operational planning. Additionally, it is attractive to developers of advanced vehicle technologies, as we demonstrate the benefits of NExT E-MAUs. }


\maketitle

\section{Introduction}

The climate change mitigation goals established by the \cite{ipcc2023ar6} highlight the need for widespread electrification across societies. From a public transportation perspective, replacing traditional internal combustion engine buses with electric vehicles can greatly reduce pollution and greenhouse gas emissions, aligning with these objectives. This shift has driven a global surge in electric bus adoption, with a notable 53\% increase in European e-bus registrations in 2023 \citep{sustainableBus}. A further way to reduce pollution is to adopt \textit{Electric Modular Autonomous Vehicles} (E-MAVs). An E-MAV comprises multiple \textit{Electric Modular Autonomous Units} (E-MAUs, hereafter also referred to as \textit{units} for simplicity), which can be flexibly coupled and decoupled on roads to better match the transportation capacity with time-varying passenger demand \citep{NextBus}. The number of units forming an E-MAV is referred to as its \textit{compositions}. However, revolutionizing the electrification of public transportation systems requires large investments from the construction of charging stations, to the purchase of E-MAUs, to their operations. Additionally, the passenger demand in public transportation systems is widely recognized as uncertain, presenting a further challenge to large-scale electrification efforts. To address these challenges, advanced optimization tools are important for efficiently and reliably integrating electrified technologies into transit planning processes.

This paper studies the robust charging station location and routing-scheduling problem (E-RCRSP) in an inter-modal transit system that integrates \textit{regular bus services} and \textit{short-turning services}. Regular bus services can either cover the entire line or serve part of it, while short-turning services operate on partial segments. Both services on partial segments are designed to provide faster commuting options, reduce the travel times of passengers, and improve the circulation efficiency of units. Each regular bus service starts and ends in the same operational direction, whereas short-turning services have starting and ending stations situated in opposite operational directions. We jointly determine the optimal locations of charging stations and the number of charging posts to build, purchase a fleet of E-MAUs, and manage these units to complete services by planning the sequence of passing stations and charging schedules for each E-MAU. 

\begin{figure}[h]
    \centering
\includegraphics[width=0.9\linewidth]{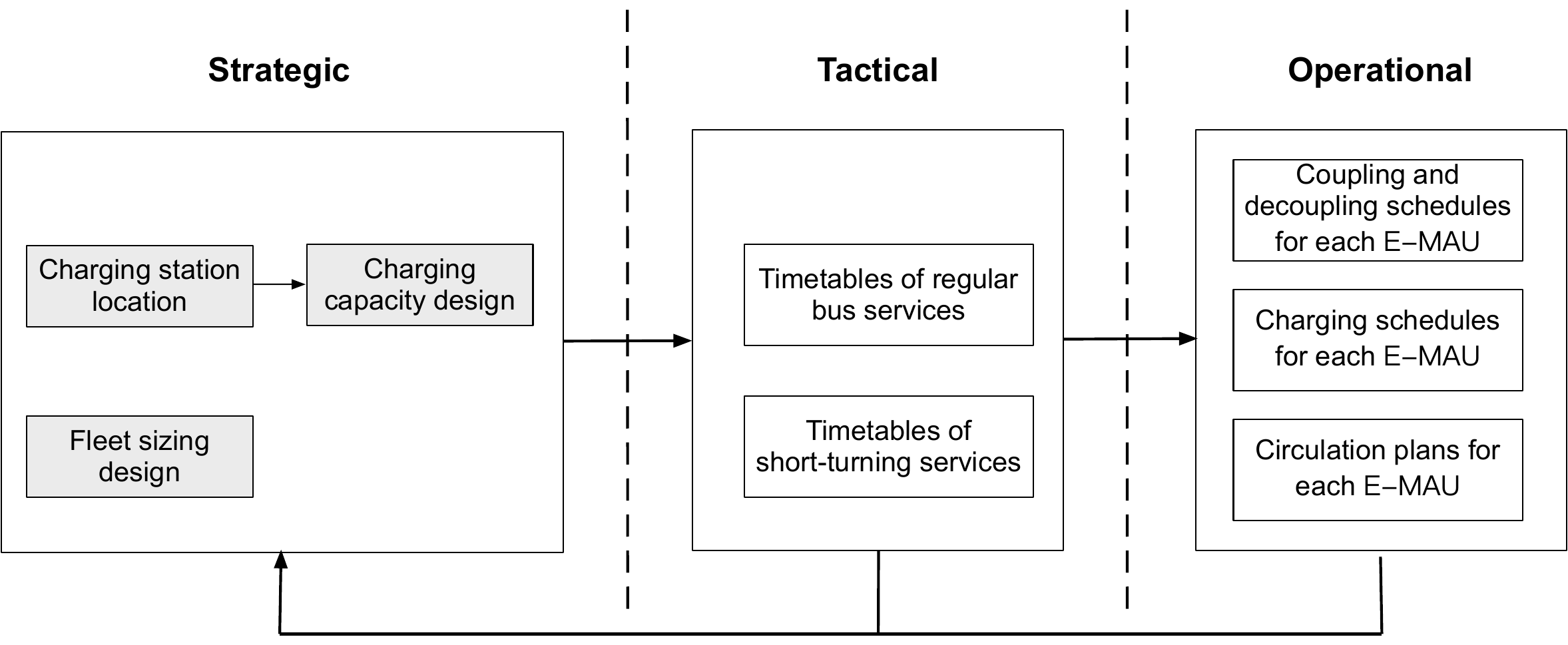}
    \caption{The proposed decision-making framework.}
    \label{fig:process}
\end{figure}

The E-RCRSP addresses the flexible scheduling and routing for units between two types of services and considers when, where, and for how long each E-MAU should recharge between assigned services. As is a variant of the \textit{electric vehicle scheduling problem}, which aims to construct feasible duties to cover a set of timetabled regular bus services from the first to the final stations \citep[e.g.,][]{JANOVEC2019, Wu2022, Rolf2023}, the E-RCRSP determines the strategic-level planning decisions and evaluates them through optimizing tactical-level and operational-level plans, as illustrated in Figure \ref{fig:process}. The strategic-level decisions aim to minimize investment costs, including the construction of capacitated charging stations and the purchase of E-MAUs. Tactical-level and operational-level planning layers aim to minimize operational costs by optimizing charging and scheduling plans. The tactical layer determines when and where to operate which type of services. The operational level optimizes charging schedules, coupling schedules, and circulation plans for each E-MAU to support the timetables generated by the tactical layer. The shared objective across these three planning stages is to minimize passengers' costs and penalties for unserved passengers. The E-RCRSP presents a complex optimization challenge, integrating location decisions with discrete routing-scheduling dynamics, alongside dynamics of time-varying and uncertain passenger flows.

This integrated decision-making framework has been widely studied in railway systems \citep[e.g.,][]{Huisman2005, Bach2016, Rolf2019, Rolf2021}. However, important differences distinguish E-MAU operations from railways. First, unlike railway rolling stock, where multiple carriages typically form a single operational unit with shared power systems \citep{Arianna2006}, each E-MAU possesses an independent power unit, enabling individual scheduling and recharging. Second, the operational necessity for stationary recharging periods is unique to E-MAUs and absent from railway practices. Third, E-MAUs operate within shorter decision-making time windows compared to railway vehicles, reflecting their flexibility and responsiveness to operational changes. Additionally, the state-of-charge (SOC) level of batteries significantly impacts coupling and decoupling decisions in E-MAUs, a consideration that does not exist for rolling stock in railways. The proposed three-stage planning framework addresses these unique challenges and provides substantial benefits. First, by integrating tactical and operational decisions, it improves strategic planning evaluations, reducing risks of over- or under-investment in charging infrastructure and fleet resources. Second, joint optimization of charging station locations, capacities, fleet size, and operational schedules under demand uncertainty ensures reliability and resilience, particularly in \textit{worst-case} scenarios. Finally, generating robust operational plans minimizes the necessity for costly re-planning efforts, thereby enhancing long-term efficiency and operational stability.

We formulate the E-RCRSP as a unit-based and path-based mixed-integer linear programming (MILP) model based on the time-space-SoC network, incorporating a polyhedral uncertainty set to address demand uncertainty. The time-space-SoC network is built through the discretization of both SoC levels and time periods, comprising nodes related to depots, non-charging stations, and charging stations. Arcs within this network capture unit flows between these nodes. The model assigns each E-MAU to a path, which represents a sequence of services and charging decisions.The system adopts two service patterns, regular services and short-turning services, to better align passenger demand with the available capacity of vehicles and charging infrastructure. Our approach allows each E-MAV to flexibly couple or decouple at every charging station and depot, providing theoretical
support to truly exploit the advantages of emerging technologies in practice. E-MAUs are scheduled for two types of services and are charged both at depots with standard facilities and at stations equipped with fast-charging facilities, considering capacity limitations and partial charging strategies. 

The main algorithmic contribution of this paper is the development of a double-decomposition solution method that combines the column and constraint generation (C\&CG) and column generation (CG) algorithms to provide high-quality solutions within a manageable computational time. To accelerate the solution of the pricing problem (PP) in the CG algorithm, we propose three innovative and generalizable techniques: (i) A tailored label-correcting algorithm explicitly designed for handling scheduling and SoC dimensions jointly. (ii) Two model properties are proposed, which systematically reduce the problem size in both the time and SoC dimensions by leveraging the \textit{subpath} concept to ensure scalability. (iii) \textit{Super travel arcs} are introduced to substantially decrease the number of arcs without compromising the optimality, making the algorithm highly scalable and adaptable to other routing-scheduling contexts. These three solution methods provide a generalized methodological framework applicable beyond the E-RCRSP context, offering a broadly relevant approach for a wide range of transportation optimization problems involving roust scheduling and routing.

Extensive computational experiments on real-world instances, incorporating Beijing bus lines and the technical specifications of state-of-art E-MAVs developed by NExT Future Transportation Inc. with high-performance batteries, demonstrate the effectiveness of the proposed approaches. We find that the three methodologies for downsizing the space-time-SoC network achieve speedups of 80.20\% to 94.97\% in the average computational time for solving the PP, with these improvements remaining robust across various problem scales. Furthermore, the proposed path-based model and solution method consistently deliver higher-quality solutions within shorter computational times compared to both the arc-based model solved using the C\&CG approach and the path-based model solved with the algorithm combining  Benders decomposition and CG. Additionally, the algorithm scales effectively to instances with entire-day operations, generating a near-optimal solution with an optimality gap of 0.96\%. From a practical perspective, the methodology offers significant benefits by jointly optimizing charging station locations and capacity, fleet sizing, routing, scheduling, and charging decisions for flexible-composition E-MAVs. In comparison to the proposed approach, the employment of  electric fixed-composition buses with two and four units that cannot be flexibly (de)coupled results in a 16.22\% and 28.22\% rise in passengers' costs, respectively, with the same optimized fleet investment generated by using the flexible-composition E-MAVs. These findings show that despite recent advancements in battery range and fast-charging technology, en-route fast-charging remains essential for efficient operations of E-MAUs in inter-modal transit systems.

This paper provides five key contributions: (1) introducing a novel charging station location and routing-scheduling problem for E-MAUs, integrating strategic, tactical, and operational decisions; (2) formulating a path-based MILP model built on a space-time-SoC network; (3) developing an efficient double-decomposition algorithm that combines C\&CG with CG; (4) proposing three generalizable acceleration techniques for solving the pricing problem, leveraging model structure and mathematical properties; (5) demonstrating the practical effectiveness of the proposed approach through real-world case studies with full-day operations in the Beijing bus system. The experiments are supported by industry data from NExT Future Transportation Inc. on E-MAVs with cutting-edge batteries.

The remainder of this paper is organized as follows. Section~\ref{sec:literature} discusses literature related to this research. Section~\ref{sec:problem} gives a detailed description of the problem and the construction of the space-time-SoC network. In Section~\ref{sec:model}, we formulate mathematical models and propose mathematical properties. In Section~\ref{sec:algorithm}, we propose solution methods. In Section~\ref{sec:results}, we present the numerical results. Lastly, we conclude in Section~\ref{sec:conclusion}.

\section{Literature review}
\label{sec:literature}

This paper contributes to the literature on the routing-scheduling of flexible-composition vehicles in electrified public transport systems under uncertainty. In this section, we provide an overview of research on charging and scheduling for electric buses, as well as scheduling for vehicles with flexible compositions.

\subsection{Charging and scheduling of electric buses}

Over the last five years, the scheduling problem of electric buses has attracted research attention \citep[e.g.,][]{Ma2021, Zhou2022, Gkiotsalitis2023}. For a comprehensive overview of related problems and methodologies, we refer to \cite{Perumal2022}. One research stream addresses the strategic-level charging station location problem but does not integrate other decision-making levels \citep[e.g.,][]{Okan2019}.  The other branch incorporates charging station capacity into the scheduling problem but has limitations with respect to time-varying passenger demand, flexible vehicle compositions, or charging station locations. For instance, \cite{Tang2019} proposed two models for scheduling electric fixed-composition buses with stochastic travel times, aiming for robust schedules under variable traffic conditions. This research considers charging station capacity but only fast charging. \cite{Wu2022} explored the multi-depot electric vehicle scheduling problem with capacitated charging stations, assuming each charging event is a full charge with a fixed duration. Another line of work considers scheduling and partial charging but excludes charging station location decisions \citep[e.g.,][]{Wen2016, Parmentier2023}.

Research on the charging and scheduling of electric buses with both capacitated charging stations and partial charging remains limited. \cite{JANOVEC2019} was among the first to address this topic, considering a single depot and fixed-composition vehicles. \cite{Zhang2021} studied electric bus scheduling with battery degradation and a non-linear charging profile, assuming a single terminal as the charging depot. Recently, \cite{Rolf2023} studied electric vehicle scheduling with capacitated charging stations, partial charging, and multi-type buses. They proposed a path-based model solved by a column generation-based heuristic but assumed fixed vehicle compositions during operations and pre-given timetables. \cite{Jacquillat2024} proposed a sub-path column generation algorithm for the electric routing-scheduling problem with uncapacitated charging stations, which jointly optimizes charging and routing-scheduling decisions of fixed-composition vehicles.

Our problem differs from previous studies in three main ways. First, we integrate routing, scheduling, and charging station location decisions in an electrified public transport system. Second, leveraging advanced vehicle technologies, we incorporate E-MAVs with flexible compositions across different times and locations. We further consider two operational patterns and allow cross-pattern unit circulation to enhance resource utilization. Third, we account for time-varying and uncertain passenger demand to better evaluate solution quality in terms of passengers' and operators' costs. This paper proposes a novel modeling and solution framework for charging station location and routing-scheduling in public transit with E-MAVs, extending beyond the traditional electric vehicle scheduling problem.

\subsection{Scheduling of vehicles with flexible compositions}

In rail transit, bus, and on-demand mobility systems, vehicles with flexible compositions offer operators the opportunity to save on operating costs due to their decoupling and coupling capabilities at different times and locations. Initially explored in the field of railways, this concept has recently been extended to urban rail transit, bus, and on-demand mobility systems.

In railway and urban rail systems, as noted in \cite{Arianna2006}, trains are typically composed of four carriages powered by a single power device. In this context, these four-carriage configurations are commonly treated as a single operational unit, consistent with the definition of a unit in this study. Consequently, each train has two composition options: either one unit or two units. Relevant studies on this topic fall into two categories. In one category, trains are permitted to couple and decouple at depots \citep[e.g.,][]{Yin2023, Chai2024, Wang2024}. In the other one, trains can couple and decouple at both transfer stations and depots \citep[e.g.,][]{Ralf2016, Entai2024}. For example, \cite{Arianna2006} developed a multicommodity flow model for efficient railway rolling stock circulation. In addition, relaxing the strong assumption on the two options in compositions, \cite{Cacchiani2010} addressed the train unit assignment problem, allowing multiple types of units and proving this problem to be strongly NP-hard. Later, \cite{Cacchiani2019} introduced a heuristic algorithm to solve this assignment problem. 

In bus systems, various studies have focused on integrating the timetabling and composition design of E-MAVs on a one-way line with deterministic passenger demand. Some permit coupling and decoupling only at depots, see, \cite{Chen2019, Shi2021, Chen2022}, while others adopt a more flexible approach by allowing units to be coupled and decoupled at stations along the line \citep[e.g.,][]{Chen20212, Liu2023}. Solution methods in this field include commercial solvers, heuristics, dynamic programming, and simulation. Recently, \cite{Xia2023} and \cite{Xia2024} addressed this problem under uncertain passenger demand. \cite{Xia2023} proposed distributionally robust optimization models for timetabling and dynamic capacity allocation on one-way bus lines. Furthermore, \cite{Xia2024} formulated models for integrated timetabling and vehicle scheduling in inter-modal transit systems combining fixed-line and on-demand services. Building on these, \cite{Xia2024integrated} extends the timetabling and vehicle scheduling problem from the line level to the network level, incorporating flexible compositions at each station, cross-line circulation of units, and in-vehicle passenger transfers. Additionally, \cite{Cats2023} demonstrated the advantages of E-MAVs in on-demand mobility.

Our approach leverages recent advancements in modular autonomous vehicle technologies to extend train unit assignment problems in railways, enabling vehicle compositions with greater flexibility than the conventional restriction to one or two units. Moreover, by jointly considering charging station locations, timetabling, and dynamic composition adjustments across different times and locations, our approach better supports decision-making in the context of continuously evolving vehicle technologies. We also combine two service patterns and integrate electrification requirements, resulting in routing-scheduling decisions. Additionally, by explicitly modeling time-dependent segment running times, our approach provides a more accurate representation of real-world operational dynamics. Finally, charging station locations and duration requirements are directly embedded into our scheduling process, capturing essential electrification considerations.

\section{Problem description}
\label{sec:problem}
In this section, we formally define the E-RCRSP for electric modular
autonomous units. Thereafter, we construct a space–time-SoC network to model the investigated problem.

\subsection{Definitions of the E-RCRSP for E-MAUs }

We consider a bi-directional bus line with a set of stations denoted by $\mathcal{S} = \overline{\mathcal{S}} \cup \underline{\mathcal{S}}$, where $\overline{\mathcal{S}}=\{1,2,...,\left|\overline{\mathcal{S}}\right|\}$ and $\underline{\mathcal{S}}=\{\left|\overline{\mathcal{S}}\right|+1,\left|\overline{\mathcal{S}}\right|+2,...,\left|\overline{\mathcal{S}}\right|+\left|\underline{\mathcal{S}}\right|\}$ represent stations in the up and down directions, respectively. The physical road between two adjacent stations $s$ and $s+1$ for all $s \in \mathcal{S} \setminus \{\left|\overline{\mathcal{S}}\right|, \left|\overline{\mathcal{S}}\right|+\left|\underline{\mathcal{S}}\right| \}$ is defined as a \textit{segment}, with the set of segments represented by $\mathcal{Q}$. To capture real-life operations, we incorporate time-dependent running times of E-MAUs on each segment and account for time-varying and uncertain segmental demand, defined as the number of passengers arriving at each segment over time. Stations are categorized as either \textit{non-charging stations} without charging capability or \textit{charging stations}. Depots, located at both ends of the line and denoted by $\mathcal{D}$, are equipped with a predefined number of standard charging facilities and serve as storage locations for E-MAUs. The set of other stations with fast-charging capabilities, excluding the depots, is represented by $\mathcal{U}_{fast}$. To ease the notation, we define $\mathcal{U} = \mathcal{D} \bigcup \mathcal{U}_{fast}$ as the set of all charging stations on the line, each of which has strict charging capacity constraints. The number of fast-charging posts at each charging station $u \in \mathcal{U}_{fast}$ needs to be determined, subject to infrastructure restrictions $r_u$, which denote the maximum number of charging posts that can be installed at station 
$u$ due to physical space or electrical capacity constraints. 

The bus manager is responsible for allocating a fleet of E-MAUs to execute a set of regular bus and short-turning services throughout the studied time horizon $\mathcal{T}$, with the fleet size being a decision variable. Each regular bus and short-turning service is operated by an E-MAV composed of one or more E-MAUs, with the number of E-MAUs defining the vehicle's \textit{composition}. Each E-MAV has a \textit{flexible composition} at different times and locations, allowing it to consist of multiple E-MAUs that can be dynamically coupled or decoupled at every charging station. The maximum composition of an E-MAV is denoted as $N^{com}$. The E-MAUs within an E-MAV may exhibit various states of charge (SoCs), denoted as $e \in \mathcal{E}$. They can recharge at charging stations, where the minimum and maximum SoC values are represented by $e_{min}$ and $e_{max}$, respectively. We incorporate the partial charging strategy, where each \textit{charging action} does not require a full recharge of the unit's battery. Furthermore, when an E-MAU is decoupled from an E-MAV at a fast-charging station $u \in \mathcal{U}_{fast}$, it can be flexibly reallocated to one of the following three options:

(i) \textit{Recharge}. The E-MAU decouples at a charging station $u$ and recharges at this station.

(ii) \textit{Regular bus services}.  The E-MAU decouples at charging station $u$ and is assigned to execute a subsequent service that continues in the same travel direction as its previously completed services, ending at either another charging station or the depot.

(iii) \textit{Short-turning services}. The E-MAU decouples at charging station $u$ and is assigned to execute a service in the opposite direction of its finished services, ending at either a charging station or the depot in that operational direction.

\begin{example}
\label{ex:Services}
An illustration of the decoupled and reallocation options of E-MAUs is illustrated in Figure~\ref{fig:reallocation}, where Figure~\ref{fig:reallocation}(a) and (b) represent the aforementioned reallocation options, respectively. For instance, in Figure~\ref{fig:reallocation}(a), the operated E-MAV is decoupled into a new E-MAV consisting of two units, while a single E-MAU with the lowest SoC is separated for recharging at the station. Meanwhile, the newly formed two-unit E-MAV continues to execute the regular bus service. Figure~\ref{fig:reallocation}(b) shows that the E-MAV is decoupled into three individual E-MAUs at charging station 5. The orange E-MAU, with the lowest SoC, is recharged at this station. The blue E-MAU, with the highest SoC, executes a regular bus service to charging station 8, followed by a short-turning service back to charging station 5 and recharge here. The brown E-MAU continues its task for the regular bus service.
\begin{figure}[h]
     \centering
     \begin{subfigure}{\textwidth}
         \centering
\includegraphics[width=\textwidth]{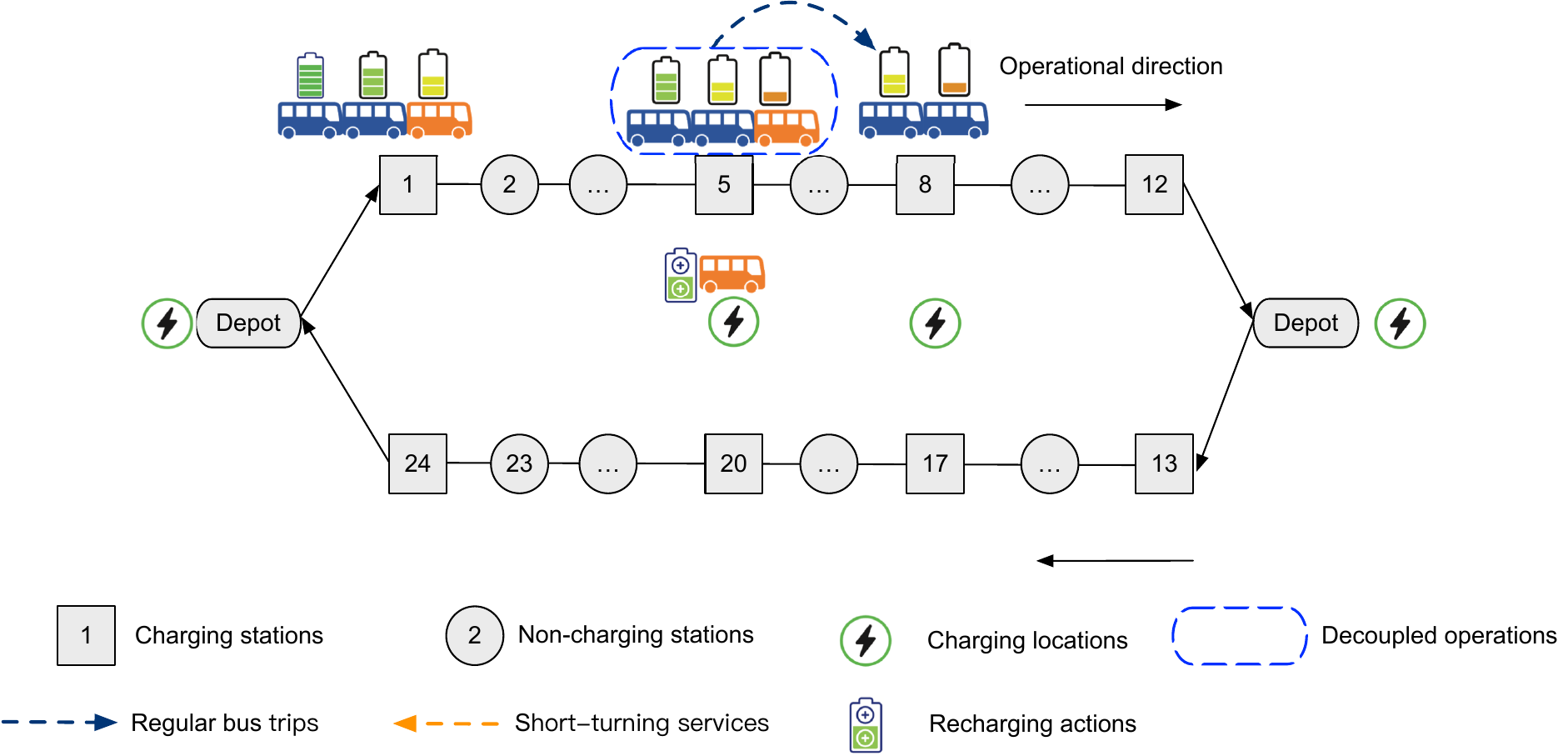}
         \caption{The first decoupled and reallocation option}
     \end{subfigure}
     
\vspace{1em}    

\begin{subfigure}{\textwidth}
         \centering
\includegraphics[width=0.85\textwidth]{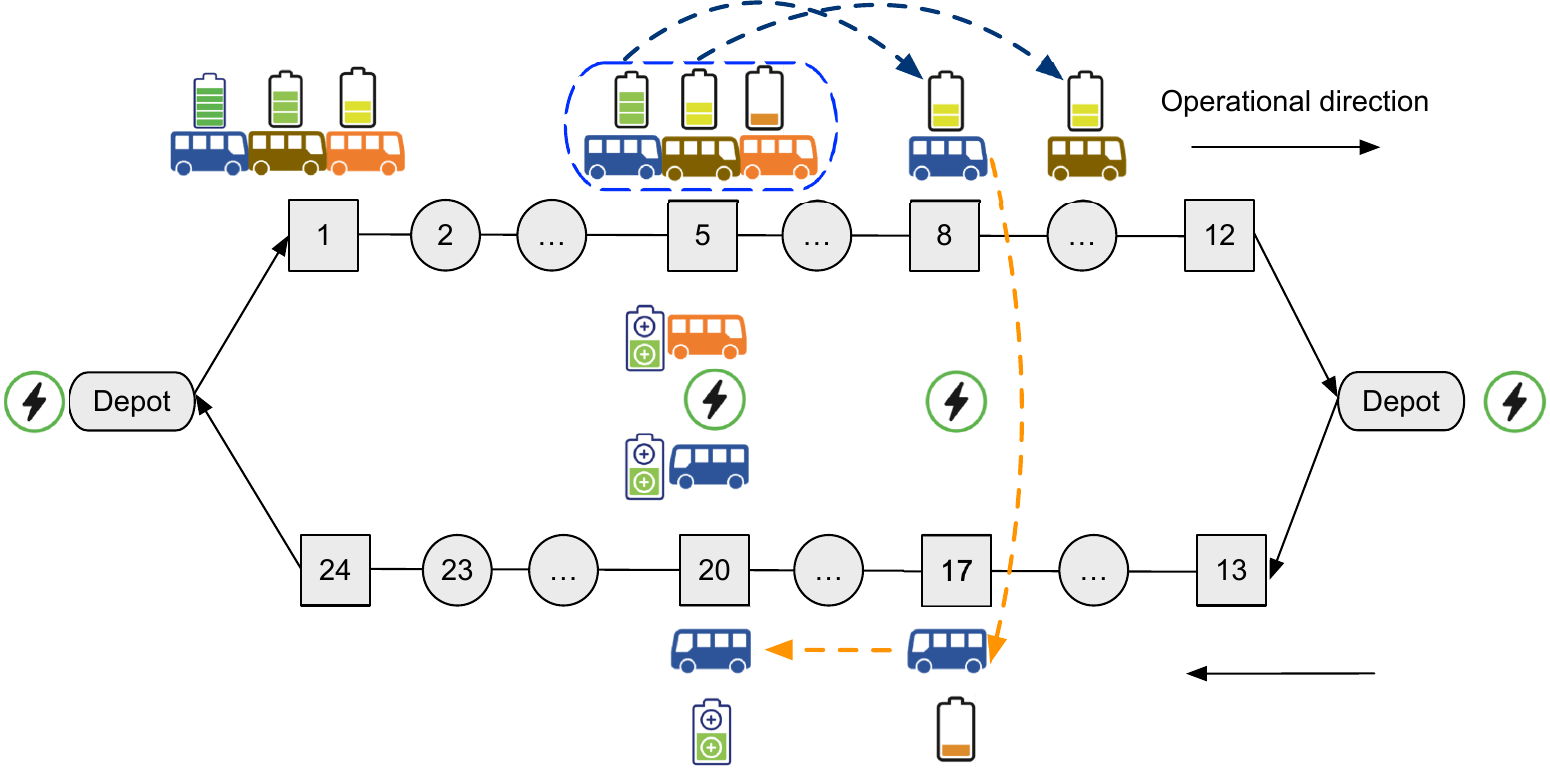}
         \caption{The second decoupled and reallocation option}
     \end{subfigure}
        \caption{Illustration of the decoupled and reallocation options of E-MAUs.}
    \label{fig:reallocation}
\end{figure}
\end{example}

To summarize, this paper addresses the robust charging station location and scheduling problem for E-MAUs, incorporating two types of services, flexible compositions, partial charging, and capacitated charging stations. Unlike common assumptions in the literature, we relax the requirements that a vehicle’s SoC must be sufficient to complete an entire service from the first to the final station, that charging actions can only occur at depots, and that all units within a vehicle must have uniform SoCs. By adopting flexible compositions at different stations and times, we allow for varying SoCs among the E-MAUs within an E-MAV. Furthermore, we ensure flexible circulation of E-MAUs across different operational directions and two types of services to improve the efficiency of unit utilization. The goal of this paper is to determine robust charging station locations and fleet sizing over the planning period with integrated circulation schedules for each E-MAU and capacity design for each built charging station. Our objective is to minimize the weighted sum of passengers' and operator's costs, including fleet investment costs, charging infrastructure costs, and variable operating costs.

To formulate the described problem, we make the following assumptions without loss of generality: 

\textbf{Assumption 1.} 
The depots of the investigated bus line are located at the two terminal stations, which are equipped with regular charging facilities.

\textbf{Assumption 2.} 
All E-MAUs are homogeneous, with the same purchasing costs, battery capacities, electricity consumption rates during travel, charging dynamics, and fixed charging costs.

\textbf{Assumption 3.} 
Each E-MAV can be decoupled at any charging station, and the decoupled E-MAUs can be coupled with other vehicles. 

Assumption 1 is consistent with the reality that bus depots are generally located at the first and last stations of a line. Assumptions 2 and 3 are commonly adopted in the scheduling field of E-MAUs, as seen in studies such as \cite{Xia2023} and \cite{Cats2023}.

\subsection{Construction of the space-time-SoC network}\label{section:network_constraction}

To facilitate modeling the movements of trains, time-dependent passenger demand, and dynamics of SoCs, we first discretize both the time horizon and SoC into discrete sets with finite numbers of elements, denoted as $\mathcal{T}$ and $\mathcal{E}$, respectively. Based on this discretization, we construct a time-space-SoC network, represented as $\mathcal{G} = \{\mathcal{N}, \mathcal{A}\}$, where $\mathcal{N}$ and $\mathcal{A}$ denote the sets of nodes and arcs. This network includes three types of nodes and various types of arcs. An example of the time-space-SoC network created for the timetables and vehicle schedules of an E-MAV made up of two E-MAUs is illustrated in Figure~\ref{time_space_soc_network} in Appendix \ref{sec:example}. 

\subsubsection{Nodes.} 

In the built network $\mathcal{G}$, we define depot nodes, station nodes, and charging nodes. These sets are denoted as $\mathcal{N}^{depot}$, $\mathcal{N}^{station}$, and $\mathcal{N}^{charge}$, respectively. The complete set of nodes is given by $\mathcal{N}=\mathcal{N}^{depot}\bigcup\mathcal{N}^{station} \bigcup \mathcal{N}^{charge}$.

\textbf{Depot nodes.} For each depot $d \in \mathcal{D}$, we include a source node $n^{source}_d$ and a sink node $n^{sink}_d$ to keep track of the numbers of incoming and outgoing E-MAUs in this depot. The set consisting of all $n^{source}_d$ and $n^{sink}_d$ for each depot 
$d \in \mathcal{D}$
is denoted as $\mathcal{N}^{depot}_d$: 
\begin{align*}
    \mathcal{N}^{depot}_d=\left\{n^{source}_d, n^{sink}_d  \right\}.
\end{align*}
Hence, we have $\mathcal{N}^{depot} = \bigcup\limits_{d \in \mathcal{D}}\mathcal{N}^{depot}_d$.

\textbf{Station node.} Each station node is characterized by three-dimensional attributes: the station $s \in \mathcal{S}$, the discretized time interval $t \in \mathcal{T}$, and the discretized SoC $e \in \mathcal{E}$. It is worthy noticing that the station node at the depot is different from the depot nodes within the space-time-SoC network. This results in the set of station nodes $\mathcal{N}^{station}$:
\begin{align*}
    \mathcal{N}^{station}=\big\{(s, t, e) \mid s\in\mathcal{S}, t\in\mathcal{T}, e\in\mathcal{E} \big\}.
\end{align*}

\textbf{Charging node}. To characterize charging operations of each E-MAU at every charging station $u \in \mathcal{U}$ and time interval $t$, we introduce the set of charging node $ \mathcal{N}^{charge}$, which can be expressed as
\begin{align*}
    \mathcal{N}^{charge}=\big\{(u,t,e) \mid u\in\mathcal{U}, t\in\mathcal{T}, e\in\mathcal{E} \big\}.
\end{align*}
Here, the value $e$ indicates the SoC of the E-MAU that stays at the charging station $u$ at time interval $t$.

\subsubsection{Arcs.}

The set of arcs $\mathcal{A}$ represents connections between nodes within the time-space-SoC network. To capture the dynamics of E-MAUs, we construct arcs corresponding to depot activities, E-MAU operations (including travel, dwell, and short-turning), and charging actions.

\textbf{Arcs related to the depot.} We first introduce the entering-depot arcs and leaving-depot arcs to track the time-varying numbers of E-MAUs that enter and leave each depot during operations. These arcs are denoted by the sets $\mathcal{A}^{in}_{d}$ and $\mathcal{A}^{out}_{d}$ at depot $ d\in \mathcal{D}$, respectively. Each entering-depot arc represents an E-MAU departing at time \( t \) with SoC \( e \) from station \( s \), and destined for \( d \). For each \( d \in \mathcal{D} \), the set of such arcs from the closest station (denoted as $f_d$) is referred to as:
\begin{align*}
    &\mathcal{A}^{in}_{d}=\big\{(n, v) \mid n=(f_d,t, e)\in\mathcal{N}^{station},v=n^{source}_d, t\in\mathcal{T}, e\in\mathcal{E} \big\}.
\end{align*}

For each depot \( d \in \mathcal{D} \), the set of leaving-depot arcs connecting this depot to the nearest station is:
\begin{align*}
    &\mathcal{A}^{out}_{d}=\big\{(n, v) \mid n=n^{source}_d, v=(f_d,t,e)\in\mathcal{N}^{station},  t\in\mathcal{T}, e=e^{full}-\theta^{soc}_{(n,v)}\big\},
\end{align*}
where $e^{full}$ represents the fully charged state. Each leaving-depot arc represents an E-MAU departing from the depot $d$ with a full SoC, arriving at station $f_d$ at time $t$ with SoC $e$. Here, $\theta^{soc}_{(n,v)}$ indicates the consumption of electricity on this arc.

\textbf{Arcs related to operations.} We develop three types of arcs that are related to operations, i.e., travel, dwell, and short-turning arcs. The sets are defined as $\mathcal{A}_{travel}$, $\mathcal{A}_{dwell}$, and $\mathcal{A}_{turn}$, respectively. 

(i) \textit{Travel arcs.} Travel arcs represent the movements of E-MAUs between stations. The cost of these arcs is positively correlated with the time-varying running time $\chi_{n,v,t}$ from node $n$ to node $v$ at time $t$. The dwell time at stations to ensure sufficient time for passengers to board and alight from the E-MAU is incorporated into $\chi_{n,v,t}$. The set of travel arcs is defined as:
\begin{align*}
    \mathcal{A}_{travel}=\big\{(n, v) \mid n=(s,t,e)\in\mathcal{N}^{station}, v=(s+1,t',e')\in\mathcal{N}^{station},(s,s+1)\in\mathcal{Q},
    e'=e-\theta^{soc}_{(n,v)},  \\ t'=t+\chi_{n,v,t}\big\},
\end{align*}
where $\theta^{soc}_{(n,v)}$ indicates the consumption of electricity from node $n$ to $v$. This type of arc describes an E-MAU arriving at station $s$ at time $t$ with SoC $e$ and subsequently arriving at station $s+1$ at time $t'$ with SoC $e'$. 

(ii) \textit{Dwell arcs.} Note that the dwell time for passengers boarding and alighting at all stations in each operational direction is aggregated into the travel time of the corresponding segments. To represent the extra dwell time of each E-MAU between the completion of one service and the execution of the next service at terminal stations (i.e., $s\in \{\left|\overline{\mathcal{S}}\right|,\left|\overline{\mathcal{S}}\right|+\left|\underline{\mathcal{S}}\right| \}$), we define the set of dwell arcs as:
\begin{align*}
    \mathcal{A}_{dwell}=\Big\{(n, v) \mid n=(s,t,e)\in\mathcal{N}^{station}, v=(s,t+\Delta, e)\in\mathcal{N}^{station}, s\in \{\left|\overline{\mathcal{S}}\right|,\left|\overline{\mathcal{S}}\right|+\left|\underline{\mathcal{S}}\right| \} \Big\},
\end{align*}
where $\Delta$ represents the duration of a discretized time interval. The cost of dwell arcs is 0 because the units are static and not being operated.

(iii) \textit{Short-turning arcs.} The short-turning arcs are developed to represent an E-MAU that changes its operating direction by turning around at a station $s$ with fast-charging capabilities, connecting to its counterpart in the opposite operational direction (denoted as station $\dot{f}_s\in\mathcal{S}$). The cost for these arcs is positively correlated with the running time $\chi_{n,v,t}$. This results in the set
\begin{align*}
    \mathcal{A}_{turn}=\big\{(n, v) \mid n=(s,t,e)\in\mathcal{N}^{station}, v=(\dot{f}_s,t',e')\in\mathcal{N}^{station}, s\in\mathcal{U}_{fast}, e'=e-\theta^{soc}_{(n,v)}, t'=t+\chi_{n,v,t}\big\}.
\end{align*}

\textbf{Arcs related to charging actions.} To model processes of charging actions, we build start-charging arcs, charging arcs, and end-charging arcs separately. The corresponding sets are denoted as $\mathcal{A}_{charge}^{start}$, $\mathcal{A}_{charge}$, and $\mathcal{A}_{charge}^{end}$. 

(i) \textit{Start-charging arcs.} The start-charging arcs are used to depict the operations of the E-MAU preparing for charging at the charging station $u$ at time $t$ and beginning the charging process at time $t'$. The SoC $e$ of E-MAUs on this arc cannot be at full charge (denoted as $e^{full}$), as charging would not be necessary in such cases. Following the approach in \cite{Rolf2023}, the costs of start-charging arcs are defined as $c^{fix}_{charge}$. This set of arcs is denoted as
\begin{align*}
\mathcal{A}_{charge}^{start} = \left\{ (n, v) \,\middle|\, n = (s, t, e) \in \mathcal{N}^{station}, v = (u, t', e) \in \mathcal{N}^{charge}, t' = t + \chi_{n,v,t}, e \in \mathcal{E} \setminus \{e^{full}\} \right\}.
\end{align*}

(ii) \textit{Charging arcs.}
Charging arcs represent the charging actions performed at charging stations, which are defined as
\begin{align*}
    \mathcal{A}_{charge}=\Big\{(n, v) \mid n=(u,t,e)\in\mathcal{N}^{charge}, v=(u,t+\Delta,e')\in\mathcal{N}^{charge}, e \in \mathcal{E} \setminus \{e^{full}\}, e'=e+E_{e,u}\Big\},
\end{align*}
where $E_{e,u}$ represents the charging quantity of the E-MAU with an initial SoC of $e$ at charging station $u$. The cost for these arcs is positively correlated with the increase in SoC, that is, $E_{e,u}$.

(iii) \textit{End-charging arcs.}
The end-charging arcs are used to describe the operation where an E-MAU completes charging at the charging station $u$ at time $t$ and becomes ready for dispatch at time $t'$. The time cost on this type of arcs is denoted as $\chi_{n,v,t}$, representing the time required to complete the preparation for being dispatched after charging at station $u$ at time $t$. This results in the set
\begin{align*}
    \mathcal{A}_{charge}^{end}=\big\{(n, v) \mid n=(u,t,e)\in\mathcal{N}^{charge}, v=(u,t',e)\in\mathcal{N}^{station}, t'=t+\chi_{n,v,t}\big\}.
\end{align*}

\begin{definition}[Path]
A path \( p \in \mathcal{P} \) is defined as an ordered sequence of nodes $
p = \{n_p^1, n_p^2, \dots, n_p^m\},
$ subject to the following conditions:
(\expandafter{\romannumeral1}) The first and last nodes of the path are depot nodes, i.e., $n_p^1 \in \left\{ n_d^{source} \mid d \in \mathcal{D} \right\}$ and $n_p^m \in\left\{ n_d^{sink} \mid d \in \mathcal{D} \right\}$.
(\expandafter{\romannumeral2}) For any $1 \leq j \leq m$, the arc $(n_p^{j-1}, n_p^{j})$ belongs to the set of arcs $\mathcal{A}$.
\end{definition}

Based on the aforementioned problem definition, the detailed mathematical formulations and solution methods for the E-RCRSP will be introduced in the subsequent sections.

\section{Mathematical formulations}
\label{sec:model}
We now present models for the E-RCRSP with E-MAUs in inter-modal transit systems. Section~\ref{sec:Notations} introduces the notations used throughout the models. In Section~\ref{sec:determinsticModel}, we formulate a MILP model that accounts for deterministic passenger demand, incorporating the variables and constraints foundational to the following robust optimization model. Finally, a robust optimization formulation under demand uncertainty for the studied problem is developed in Section~\ref{sec:robustModel}.

\subsection{Notations}
\label{sec:Notations}

To model the E-RCRSP, we summarize all sets and parameters in Table~\ref{table-para} in Appendix \ref{sec:notations}. Additionally, we define four sets of decision variables as follows.

(i) $x_{p}$ for all $p \in \mathcal{P}$. This binary variable is set to 1 if an E-MAU is assigned to path $p$; otherwise, $x_{p}=0$. It is used to determine the vehicle scheduling plan in the solution. 

(ii) $r_{u}$ for all $u \in \mathcal{U}_{fast}$. This integer variable represents the capacity of fast-charging station $u$, i.e., the number of fast-charging posts at station $u$. This variable is used to determine the charging station locations and capacities in the solution.

(iii) $y_q^t$ for all $q \in \mathcal{Q}$ and $t \in \mathcal{T}$. This continuous variable indicates the number of waiting passengers on section $q$ at time $t$.

(iv) $b_q^t$ for all $q \in \mathcal{Q}$ and $t \in \mathcal{T}$. This continuous variable denotes the number of on-board passengers on section $q$ at time $t$.

\subsection{Formulation under deterministic demand}
\label{sec:determinsticModel}

We now formulate the E-RCRSP for E-MAUs with deterministic passenger demand. The parameters $c_p$ and $m_u$ represent the operational cost of path $p$ and the costs of installing a fast-charging post at a fast-charging station $u$. Here, the operational cost $c_p$ includes both the purchase and operating costs of an E-MAU along path $p$. Each coefficient $\beta_{d,p}^{-}$ is 1 if path $p$ begins at depot $d$, and 0 otherwise.
Similarly, the coefficient $\beta_{d,p}^{+}$ is 1 if path $p$ ends at depot $d$. The parameter $l_p^a$ indicates whether path $p$ passes through arc $a$. We define the parameter $g_{u,t,p}$ is 1 if path $p$ is charged at charging station $u$ at time $t$. The parameters $\hat{\alpha}_{q}^{t}$ and $C$ represent the number of newly arrival passenger demand on section $q$ at time $t$ under the deterministic condition and the capacity of an E-MAU. We use the following formulation for the E-RCRSP.
\begin{mini!}[2]
{\mathbf{x}, \mathbf{r}, \mathbf{y}, \mathbf{b}}
{\theta_{1}(\sum_{p\in\mathcal{P}}c_p x_p + \sum_{u\in\mathcal{U}_{fast}}m_u r_u) + \theta_{2}\sum_{q\in\mathcal{Q}}\sum_{t\in\mathcal{T}} y^t_{q} + \theta_{3}\sum_{q\in\mathcal{Q}}y^{\left|\mathcal{T}\right|}_{q}\label{eq:objective}}
{\label{formulation:D}} 
{} 
\addConstraint
{\sum_{p\in\mathcal{P}}\beta_{d,p}^{-} x_p}
{=\sum_{p\in\mathcal{P}}\beta_{d,p}^{+} x_p \label{cons_balance}}
{~ \quad\forall d\in\mathcal{D},}
\addConstraint
{\sum_{a\in\mathcal{A}_{q,t}}\sum_{p\in\mathcal{P}}l^a_p x_p}
{\leq N^{com} \label{eq:max_coposition}}
{~ \quad \forall q\in\mathcal{Q},t\in\mathcal{T},}
\addConstraint
{}{\sum_{p\in\mathcal{P}}g_{u,t,p} x_p \leq
    \begin{cases}
        r_u, & \text{if } u \in \mathcal{U}_{fast} \\
        N_d^{charge} , & \text{if } u \in \mathcal{U}  \setminus \mathcal{U}_{fast}
    \end{cases} \label{eq:capacity_charge_station}
    }{~ \quad t\in\mathcal{T},}{}
\addConstraint
{y_{q}^t}
{= y_{q}^{t-1} + \hat{\alpha}_{q}^t - b_{q}^t \label{eq:passenger_flow}} 
{~ \quad \forall q\in\mathcal{Q}, t\in\mathcal{T},}
\addConstraint
{ b_{q}^t}
{\leq  C \sum_{a\in\mathcal{A}_{q,t}}\sum_{p\in\mathcal{P}}l_{p}^{a} x_p \label{eq:capacity_mv}} 
{~ \quad \forall q\in\mathcal{Q},t\in\mathcal{T},}
\addConstraint
{x_p }
{ \in \{0, 1\} \label{eq:domain_x}}
{~ \quad \forall p \in \mathcal{P},}
\addConstraint
{r_u}
{\in \{0,1,...,N^{charge}_u\}}
{~ \quad \forall u\in\mathcal{U}_{fast}, \label{eq:domain_r}}
\addConstraint
{ b_{q}^t}
{\geq 0 \label{eq:domain_b}}
{~ \quad \forall q\in\mathcal{Q}, t\in\mathcal{T},} 
\addConstraint
{y_{q}^t}
{\geq 0  \label{eq:domain_y}}
{~ \quad \forall q\in\mathcal{Q}, t\in\mathcal{T}.}
\end{mini!}

Objective function~\eqref{eq:objective} minimizes the weighted sum of the operator's costs and passengers' costs. The operator's costs include investment costs for fast-charging posts and E-MAUs, as well as operating costs for dispatching E-MAUs. Passengers' costs include waiting time and penalties associated with unserved passengers. Here, $\theta_1$, $\theta_2$, $\theta_3$ represent the weighting monetary coefficients for operator’s costs, passengers' waiting costs, and penalties associated with unserved passengers, respectively. Constraints~\eqref{cons_balance} ensure that the number of units in each depot at the beginning and end of operations remains the same. Constraints~\eqref{eq:max_coposition} is formulated to ensure the maximum number of units that can compose an E-MAV does not exceed the upper limitation $N^{com}$. Constraints~\eqref{eq:capacity_charge_station} guarantee that the total number of E-MAUs charging at the charging station $u$ and the depot $d$ in any time interval $t \in \mathcal{T}$ cannot exceed its capacity of fast-charging and standard charging posts, respectively. Constraints~\eqref{eq:passenger_flow} are formulated to model the dynamics of passengers. Constraints~\eqref{eq:capacity_mv} ensure that the number of on-board passengers cannot exceed the capacity of the E-MAV. Lastly, constraints~\eqref{eq:domain_x} - (\ref{eq:domain_y}) define the domains of the decision variables.

\subsection{Robust optimization model under demand uncertainty}
\label{sec:robustModel}
Passenger demand in public transit systems is highly uncertain \citep{Lu2023}. To enhance the practical adaptability of strategic decision-making, we incorporate uncertain passenger demand into the charging station location and scheduling problem for E-MAUs. Specifically, the uncertain time-varying demand is denoted as $\alpha_{q}^t$, where $\alpha_{q}^t=\overline{\alpha}_{q}^t+ \zeta_{q}^t\cdot \Tilde{\alpha}_{q}^t$. Here, $\overline{\alpha}_{q}^t$ represents the nominal value of passenger demand, $\Tilde{\alpha}_{q}^t$ represents its maximum deviation, and the variable $\zeta_{q}^t \in \Xi$ is a scaling factor that determines the amplitude of the fluctuation for section $q$ and time $t$, with $\Xi$ representing the domain. 

In robust optimization, an uncertainty set must be employed to adequately capture the uncertain nature of the parameters while ensuring the problem remains tractable.  Here, we introduce a polyhedral uncertainty set inspired from \cite{bertsimas2004price} to characterize the uncertainty of passenger demand. This polyhedral uncertainty set has two advantages: (1) its corresponding robust counterpart provides an adjustable level of robustness, enabling a trade-off between cost and robustness, and (2) its robust counterpart can be formulated as a linear program. In our problem, we define the polyhedral uncertainty set (denote as $\mathscr{U}$) for the passenger demand as $
    \mathscr{U} = \{\bm{\alpha}: \alpha_{q}^t = \overline{\alpha}_{q}^t+ \zeta_{q}^t\cdot \Tilde{\alpha}_{q}^t, \forall q\in\mathcal{Q},t\in\mathcal{T}, \bm{\zeta} \in \Xi \}$, where the domain of the continuous variable $\zeta_{q}^t$ is formulated as
\begin{align*}
    \Xi=\left\{ \bm{\zeta}\in\mathbf{R}^{\left|\mathcal{Q}\right|\times\left|\mathcal{T}\right|}: \sum_{t\in\mathcal{T}}\zeta_{q}^t\leq\Pi_q,  \sum_{q\in\mathcal{Q}}\zeta_{q}^t\leq\Lambda_t,0\leq\zeta_{q}^t\leq 1,\forall q\in\mathcal{Q}, t\in\mathcal{T} \right\}.
\end{align*} 
This polyhedral uncertainty set adjusts the robustness of the method against the conservatism of the solution using the parameters $\Pi_q$ for all $q \in \mathcal{Q}$ and $\Lambda_t$ for all $t \in \mathcal{T}$. The value of $\Pi_q$ represents the maximum number of time intervals in which the uncertain demand can deviate from the nominal values for each section $q \in \mathcal{Q}$ at the same time. The value of $\Lambda_t$ represents the maximum number of segements where the uncertain demand can simultaneously deviate from the nominal values at each time $t \in \mathcal{T}$.

The two-stage robust optimization model for the E-RCRSP with E-MAUs is formulated as:
\begin{mini!}[2]
{}
{\theta_{1}(\sum_{p\in\mathcal{P}}c_p x_p + \sum_{u\in\mathcal{U}_{fast}}m_u r_u)  + opt[F(\mathbf{x})]}
{\label{formulation:U}} 
{} 
\addConstraint
{}
{\eqref{cons_balance} - \eqref{eq:capacity_charge_station},\eqref{eq:domain_x} - \eqref{eq:domain_r}.}
{}
\end{mini!}
where $opt[F(\mathbf{x})]$ is the optimal solution of the second-stage recourse problem. The problem $F(\mathbf{x})$ is formulated as 
\begin{align}
\max_{\bm{\zeta}\in\Xi}\ \min_{\mathbf{y},\mathbf{b}} \quad&\sum_{q\in\mathcal{Q}}\sum_{t\in\mathcal{T}} \theta_{2}y^t_{q} + \sum_{q\in\mathcal{Q}}\theta_{3}y^{\left|\mathcal{T}\right|}_{q}  \tag{2c}\\
\text{s.t.}\qquad &y_{q}^t \geq y_{q}^{t-1} + \overline{\alpha}_{q}^t+ \zeta_{q}^t\cdot \Tilde{\alpha}_{q}^t - b_{q}^t   &\forall q\in\mathcal{Q}, t\in\mathcal{T},\label{eq:passenger_flow_uncetain} \tag{2d}\\
& (\ref{eq:capacity_mv}),  (\ref{eq:domain_b})- (\ref{eq:domain_y})\label{eq:domain_uncetain}. \tag{2e}
\end{align}
To ensure the non-negativity of the dual variable, constraints (\ref{eq:passenger_flow_uncetain}) are formulated as inequalities. Besides, the inner-minimization, i.e., $F(\mathbf{x},\bm{\zeta})=\min\limits\limits_{\mathbf{y},\mathbf{b}}\left\{\sum\limits_{q\in\mathcal{Q}}\sum\limits_{t\in\mathcal{T}} \theta_{2}y^t_{q} + \theta_{3}\sum\limits_{q\in\mathcal{Q}}y^{\left|\mathcal{T}\right|}_{q}: (\ref{eq:passenger_flow_uncetain}) - (\ref{eq:domain_uncetain})\right\}$, is evidently a linear program.

Let $\kappa_{q}^t$ represent the dual variable associated with constraints (\ref{eq:capacity_mv}), and let $\delta_{q}^t$ denote the dual variable corresponding to constraints (\ref{eq:passenger_flow_uncetain}). The dual problem of $F(\mathbf{x},\bm{\zeta})$, denoted as $\text{IDF}(\mathbf{x},\bm{\zeta})$, is formulated as follows
\begin{align}
    \text{IDF}(\mathbf{x},\bm{\zeta})=
    \max_{\bm{\delta},\bm{\kappa}} &\sum_{q\in\mathcal{Q}}\sum_{t\in\mathcal{T}}\left[(\overline{\alpha}_{q}^t+ \zeta_{q}^t\cdot \Tilde{\alpha}_{q}^t)\delta_{q}^t +   C(\sum_{a\in\mathcal{A}_{q,t}}\sum_{p\in\mathcal{P}}l_{p}^a x_p)\kappa_{q}^{t}\right]\label{bilinear_obj} \tag{3a}\\
    & \delta_{q}^t - \delta_{q}^{t+1} \leq \theta_2 &\forall q\in\mathcal{Q}, t\in\mathcal{T}\backslash\left|\mathcal{T}\right|,\label{eq:dual1} \tag{3b}\\
     & \delta_{q}^{\left|\mathcal{T}\right|} \leq \theta_2 + \theta_3 &\forall q\in\mathcal{Q}, t\in\mathcal{T},\label{eq:dual1_1} \tag{3c}\\
      & \delta_{q}^t\geq 0 \label{eq:dual3} &\forall q\in\mathcal{Q}, t\in\mathcal{T}, \tag{3d}\\
    &   \delta_{q}^t+\kappa_{q}^t  \leq 0 &\forall q\in\mathcal{Q}, t\in\mathcal{T}, \label{eq:dual2} \tag{3e}\\
    &\kappa_{q}^t\leq 0 &\forall q\in\mathcal{Q}, t\in\mathcal{T}. \label{eq:dual4} \tag{3f}
\end{align}
According to inequality (\ref{eq:dual2}), it always holds that $\kappa_{q}^t\leq -\delta_{q}^t$. By merging the two maximum terms in the objective function, $F(\mathbf{x})$ can be reformulated as a single maximization problem $\text{DF}(\mathbf{x})$:
\begin{align}
    \text{DF}(\mathbf{x})=
       \max_{\bm{\zeta},\bm{\delta}}
    &\sum_{q\in\mathcal{Q}}\sum_{t\in\mathcal{T}}\left[(\overline{\alpha}_{q}^t+ \zeta_{q}^t\cdot \Tilde{\alpha}_{q}^t)\delta_{q}^t -  C\cdot \delta_{q}^t\sum_{a\in\mathcal{A}_{q,t}}\sum_{p\in\mathcal{P}}l_{p}^a  x_p \right]\label{bilinear_obj_df} \tag{4a}\\
    &\sum_{t\in\mathcal{T}}\zeta_{q}^t\leq\Pi_q  &\forall q\in\mathcal{Q},\label{df1} \tag{4b}\\
    &\sum_{q\in\mathcal{Q}}\zeta_{q}^t\leq\Lambda_t &\forall t\in\mathcal{T},\label{df_uncertain}\tag{4c}\\
    &0\leq\zeta_{q}^t\leq 1 &\forall q\in\mathcal{Q}, t\in\mathcal{T}, \label{cons:bound}\tag{4d}\\
    &  \eqref{eq:dual1} - \eqref{eq:dual3}.\label{df2} \tag{4e}
\end{align}
It is worth mentioning that the objective function contains the bilinear term $\zeta_{q}^t \delta_{q}^t$. As a result, the second-stage recourse problem becomes a bilinear programming problem, which is known to be NP-hard. When both $\zeta_{q}^t$ and $\delta_{q}^t$ for all $q \in \mathcal{Q}$ and $t \in \mathcal{T}$ are bounded and $\zeta_{q}^t$ is a binary variable in the optimal solution, $\text{DF}(\mathbf{x})$ \eqref{bilinear_obj_df} - \eqref{df2} can be equivalently transformed into a MILP form. Notice that $\zeta_{q}^t$ itself is bounded according to constraints \eqref{cons:bound}. Next, we first present Lemma \ref{lemma:bound} to prove the boundedness of $\bm{\delta}$, and then develop Proposition \ref{pp1} to prove that the optimality condition for Problem $\text{DF}(\mathbf{x})$ ensures $\bm{\zeta}$ is either zero or one in the optimal solution.

\begin{lemma}[Boundedness]
\label{lemma:bound}
    Let $\Psi=\{\bm{\delta}\in\mathbb{R}^{\left|\mathcal{T}\right|\times\left|\mathcal{Q}\right|}: 0\leq\delta^t_{q}\leq M_t, (\ref{eq:dual1}) - (\ref{eq:dual3})\}$. It holds that $\delta_q^t\leq M_t=(\left|\mathcal{T}\right|-t+1)\theta_2 + \theta_3,\forall t\in\mathcal{T}, q\in\mathcal{Q}$ and $\Psi$ is a bounded polyhedron.
\end{lemma}
\proof{Proof.} See Appendix \ref{sec:proof}.
 \endproof

\begin{proposition}[Optimality conditions]\label{pp1}
    Following \cite{gabrel2014robust}, if $\Pi_q \in \mathbb{Z}^+, \forall q\in\mathcal{Q}$ and $\Lambda_t \in \mathbb{Z}^+, \forall t\in\mathcal{T}$, then an optimal solution $(\bm{\zeta}^{*},\bm{\delta}^{*})$ of $DF(\mathbf{x})$ exists such that $\bm{\zeta}^{*}\in\{0, 1\}^{\left|\mathcal{T}\right|\times\left|\mathcal{Q}\right|}$.
\end{proposition}
\proof{Proof.} See Appendix \ref{sec:proof}.

According to Proposition \ref{pp1}, we can now reformulate the $\text{DF}(\mathbf{x})$ as follows
\begin{align}
    \text{DF}(\mathbf{x})=\max_{\bm{\zeta},\bm{\delta}} &\sum_{t\in\mathcal{T}}\sum_{q\in\mathcal{Q}}\left[\overline{\alpha}_{q}^t  \delta_{q}^t+ \Tilde{\alpha}_{q}^t \phi_{q}^t -
C\cdot\delta_{q}^t\sum_{a\in\mathcal{A}_{q,t}}\sum_{p\in\mathcal{P}}l_{p}^a x_p \right]\label{bilinear_obj_1} \tag{5a}\\
    & \phi_{q}^t\leq \delta_{q}^t,&\forall q\in\mathcal{Q}, t\in\mathcal{T}, \label{df_1_1}\tag{5b} \\
    & \phi_{q}^t\leq M_t\zeta_{q}^t &\forall q\in\mathcal{Q}, t\in\mathcal{T}, \tag{5c}\\
    & \phi_{q}^t\geq \delta_{q}^t - M_t(1-\zeta_{q}^t)&\forall q\in\mathcal{Q}, t\in\mathcal{T}, \tag{5d}\\
     & \phi_{q}^t\geq 0&\forall q\in\mathcal{Q}, t\in\mathcal{T}, \tag{5e}\\
     &\zeta_{q}^t\in\{0, 1\}&\forall q\in\mathcal{Q}, t\in\mathcal{T}, \label{df_2_1}\tag{5f}\\
    &  (\ref{eq:dual1}) - (\ref{eq:dual3}), (\ref{df1}) - (\ref{df_uncertain}).\label{df_3_1} \tag{5g}
\end{align}

Lastly, the MILP formulation for the E-RCRSP with E-MAUs is expressed as follows:
\begin{align}
\min \quad &\theta_{1}(\sum_{p\in\mathcal{P}}c_p x_p + \sum_{u\in\mathcal{U}}m_u r_u)  + \max_{\bm{\zeta},\bm{\delta}} \sum_{t\in\mathcal{T}}\sum_{q\in\mathcal{Q}}\big[\overline{\alpha}_{q}^t  \delta_{q}^t+ \Tilde{\alpha}_{q}^t \phi_{q}^t -
C\cdot\delta_{q}^t\sum_{a\in\mathcal{A}_{q,t}}\sum_{p\in\mathcal{P}}l_{p}^a x_p \big] \nonumber \\
\text{s.t.}\quad &\eqref{cons_balance} - \eqref{eq:capacity_charge_station},\eqref{eq:domain_x} - \eqref{eq:domain_r}, \eqref{df_1_1} - \eqref{df_3_1}. \label{model:MILP}\tag{6}
\end{align}
 
\section{Solution methods}
\label{sec:algorithm}

In this section, we develop a double-decomposition algorithm that combines the column-and-constraint generation (C\&CG) method with the column-generation approach. Given the strong performance of the C\&CG method developed by \cite{zeng2013solving} in addressing two-stage robust optimization problems, we adopt this approach to tackle the E-RCRSP studied in this paper. Specifically, the problem is decomposed into a master problem (MP) and a subproblem (SP). The MP determines decisions related to charging station locations and capacities, fleet sizing, and routing-scheduling plans for E-MAUs. The SP, in turn, identifies the worst-case scenario and evaluates the performance of the decisions generated by MP by solving the model with passenger-related constraints. Based on the worst-case scenario identified by the SP, recourse variables and their associated constraints are generated and added to the MP. This iterative process continues until convergence. In this study, the SP is solved using GUROBI. To accelerate the solution process, we design an outer approximation algorithm that generates a high-quality initial solution, which is fed into GUROBI for further optimization.

In our problem, the MP is formulated as a path-based model with an exponential number of variables, making direct solution impractical. To address this issue, we develop a column generation (CG) algorithm to dynamically include only the most effective paths into the problem, thereby significantly enhancing scalability and solvability. Within this framework, the solution to the linear programming (LP) relaxation of the MP is obtained iteratively by solving a \textit{restricted master problem} (RMP) and a \textit{pricing problem} (PP). The RMP is a simplified version of the LP relaxation of the MP, containing only a subset of paths. Initially, it is initialized with a set of paths that guarantee a feasible solution to the RMP. After solving the RMP in each iteration, the values of the dual variables are passed as input for the PP. The PP identifies the paths with most negative reduced costs and feedback to the RMP. The RMP is then resolved with this expanded set of paths to improve the solution quality. If the PP cannot identify a path with a negative reduced cost, the RMP is solved to optimality, and the column generation process terminates.

The PP in this study is a shortest path problem with negative weights, which is computationally challenging to solve on large-scale networks. To solve it efficiently, we adopt a two-pronged approach. First, we design a customized label-correcting algorithm tailored to our problem. Second, acknowledging that the size of the space-time-SoC network significantly impacts the efficiency of solving the PP, we propose three methodologies to downsize the network and aggregate arcs, thereby reducing its scale. These methodologies, which are based on the concept of \textit{subpaths} and \textit{super travel arcs}, have been proven to maintain the optimality of solutions. 

We develop an arc-based model and its corresponding C\&CG algorithm to serve as a benchmark for the path-based model, as presented in Appendix~\ref{sec:arc-based}. Besides, for the path-based model, we use the Benders decomposition (BD) method as a benchmark for the C\&CG algorithm. Consequently, the benchmark for the proposed algorithm, which combines C\&CG and CG, is a solution method that integrates BD and CG. Details of this benchmark solution method are provided in Appendix~\ref{sec:BD}.

In the remainder of this section, the C\&CG approach for the E-RCRSP is developed in Setion~\ref{sec:CCG}. In Section~\ref{section_cg}, we develop the CG algorithm to solve the MP in the C\&CG framework. Thereafter, the methodologies for downsizing the space-time-SoC network and aggregating arcs are introduced in Section~\ref{sec:aggregating}. The outer approximation method, designed to accelerate the solution of the SP within the C\&CG framework, is presented in Section~\ref{sec:Outer}. Lastly, the overall framework and the two procedures of the solution algorithm are presented in Section~\ref{sec:framework}.

\subsection{Column and constraints generation algorithm}
\label{sec:CCG}

Here, we introduce the C\&CG approach in detail. The C\&CG algorithm is implemented in a master-subproblem framework. At each iteration \( k \), the MP is formulated based on a subset of passenger demand scenarios identified in the previous iterations. We first solve the MP to obtain the solutions for decision variables \( \hat{\mathbf{x}} \) and \( \hat{\mathbf{r}} \), which provide a lower bound for the original two-stage robust optimization model~\eqref{model:MILP}. Thereafter, we solve \( \text{DF}(\hat{\mathbf{x}}) \) \eqref{bilinear_obj_1}--\eqref{df_3_1} to identify the worst-case scenario \( \hat{\bm{\zeta}}_k \). The upper bound is then updated accordingly. If the termination criteria are not met, the recourse variables \( \mathbf{y} \) and \( \mathbf{b} \) are generated on the fly, and the corresponding constraints associated with this scenario are added to the MP. A stronger lower bound is obtained by solving the updated MP. This process is repeated until the termination criteria are satisfied. The detailed procedure of the proposed C\&CG algorithm is presented in Algorithm~\ref{alg:ccg}. We define the set \( \mathcal{W} = \{w \mid w = 1, 2, \dots, k\} \) to represent all iterations up to the \( k \)-th, where each \( w \in \mathcal{W} \) corresponds to a previously identified worst-case scenario \( \hat{\bm{\zeta}}_w \). The MP is formulated as follows:

\begin{align}
\label{mp_model}
{\rm{\big[MP\big]}}\quad
\min\limits_{\mathbf{x}, \mathbf{r}, \mathbf{y}, \mathbf{b}}  & \theta_{1}(\sum\limits_{p\in\mathcal{P}}c_p x_p + \sum\limits_{u\in\mathcal{U}_{fast}}m_u r_u) + \eta  \tag{7a}\\
\mbox{s.t.}\quad\quad
&\sum_{p\in\mathcal{P}}\beta_{d,p}^{-} x_p=\sum_{p\in\mathcal{P}}\beta_{d,p}^{+} x_p & \forall d\in\mathcal{D},\label{mp_balance} \tag{7b}\\
&\sum_{a\in\mathcal{A}_{q,t}}\sum_{p\in\mathcal{P}}l^a_p x_p \leq N^{com} & \forall q\in\mathcal{Q},t\in\mathcal{T}\label{mp_max_coposition}, \tag{7c}\\
&\sum_{p\in\mathcal{P}}g_{u,t,p} x_p \leq
    \begin{cases}
        r_u, & \text{if } u \in \mathcal{U}_{fast} \\
        N_d^{charge} , & \text{if } u \in \mathcal{U}  \setminus \mathcal{U}_{fast}
    \end{cases} & t\in\mathcal{T},\label{mp_capacity_charge_station} \tag{7d}\\
&\eta \geq \sum\limits_{q\in\mathcal{Q}}\sum\limits_{t\in\mathcal{T}} \theta_{2}y^t_{q}(w) + \sum\limits_{q\in\mathcal{Q}}\theta_{3}y^{\left|\mathcal{T}\right|}_{q}(w) & \forall w \in \mathcal{W}, \label{mp_eta}\tag{7e}\\
&b_{q}^t(w)  \leq  C \sum_{a\in\mathcal{A}_{q,t}}\sum_{p\in\mathcal{P}}l_{p}^{a} x_p & \forall q\in\mathcal{Q}, t\in\mathcal{T},w \in \mathcal{W}, \label{mp_capacity_mv} \tag{7f}\\
&y_{q}^t(w) \geq y_{q}^{t-1}(w) + \overline{\alpha}_{q}^t+ \zeta_{q}^t(w)\cdot \Tilde{\alpha}_{q}^t - b_{q}^t(w)   & \forall q\in\mathcal{Q}, t\in\mathcal{T}, w\in \mathcal{W}, \label{mp_y}\tag{7g}\\
&y_{q}^t(w),b_{q}^t(w) \geq 0 & \forall q\in\mathcal{Q}, t\in\mathcal{T},w\in \mathcal{W}, \label{mp_y_b_domain}\tag{7h}\\
& x_p \in \{0, 1\} & \forall p\in\mathcal{P}, \label{mp_domain_x} \tag{7i}\\
& r_u \in \{0,1,...,N^{charge}_u\} & \forall u\in\mathcal{U}_{fast}. \label{mp_domain_r} \tag{7j}
\end{align}

Note that the MP is in general difficult to be solved, as the set $\mathcal{P}$ of all possible paths, is extremely large. To address this challenge, we develop a column generation algorithm that dynamically incorporates paths into the MP, as detailed in the following section.
\begin{algorithm}[]
\caption{C$\&$CG algorithm}\label{alg:ccg}
\SetKwInOut{Input}{Input}
\SetKwInOut{Output}{Output}
\SetKw{KwGoTo}{go to}
\SetKw{KwSet}{set}
\SetKw{KwBreak}{break}
\SetKw{Return}{return}
\SetAlgoLined

\Input{Space-time-SoC network $\mathcal{G} = \{\mathcal{N}, \mathcal{A}\}$; time-dependent and uncertain passenger demand.}
\Output{Optimal solutions $(\mathbf{x}^*, \mathbf{r}^*)$.}
\BlankLine

Initialize $LB = -\infty$, $UB = +\infty$; initialize iteration counter $k = 0$; initialize tolerance $\epsilon$.\\
\While{$(UB - LB)/LB > \epsilon$}{
    Solve the MP~\eqref{mp_model} - \eqref{mp_domain_r} to obtain $(\mathbf{x}_{k+1}^*, \mathbf{r}_{k+1}^*, \eta, \mathbf{y}_{1}^*, \dots, \mathbf{y}_{k}^*, \mathbf{b}_{1}^*, \dots, \mathbf{b}_{k}^*)$;\\
    Update $LB = \theta_1 (\mathbf{c}^{\mathbf{T}} \mathbf{x}_{k+1}^* + \mathbf{m}^{\mathbf{T}} \mathbf{r}_{k+1}^*) + \eta$;\\
    Call GUROBI to solve $DF(\mathbf{x}_{k+1}^*)$ \eqref{bilinear_obj_1} - \eqref{df_3_1}; obtain the optimal solution $(\zeta_{q}^t(k+1))^*$; update $UB = \min \{UB, \theta_1 (\mathbf{c}^{\mathbf{T}} \mathbf{x}_{k+1}^* + \mathbf{m}^{\mathbf{T}} \mathbf{r}_{k+1}^*) + DF(\mathbf{x}_{k+1}^*)\}$;\\
    Add the recourse variables $\mathbf{y}$ and $\mathbf{b}$ and the following constraints to the MP:
    \Indp
\[
\left\{
\begin{aligned}
    &\eta \geq \sum\limits_{q \in \mathcal{Q}} \sum\limits_{t \in \mathcal{T}} \theta_2 y^t_q(k+1) + \sum\limits_{q \in \mathcal{Q}} \theta_3 y^{\left|\mathcal{T}\right|}_q(k+1);\\
    &y_q^t(k+1) \geq y_q^{t-1}(k+1) + \overline{\alpha}_q^t + \zeta_q^t(k+1)^* \cdot \Tilde{\alpha}_q^t - b_{q}^t(k+1),  \quad \forall q \in \mathcal{Q}, t \in \mathcal{T};\\
    &b_q^t(k+1) \leq C \sum\limits_{a\in\mathcal{A}_{q,t}}\sum\limits_{p\in\mathcal{P}}l_{p}^{a} x_p, \quad \forall q \in \mathcal{Q}, t \in \mathcal{T};\\
    &y_q^t(k+1), b_{q}^t(k+1) \geq 0, \quad \forall q \in \mathcal{Q}, t \in \mathcal{T}.
\end{aligned}
\right.
\]
\Indm
    Update $k = k + 1$;\\
}
\Return{Optimal solutions $(\mathbf{x}_{k+1}^*, \mathbf{r}_{k+1}^*)$ and terminate.}
\end{algorithm}

\subsection{Column generation algorithm}
\label{section_cg}

The core idea of the CG algorithm lies in iteratively solving a RMP, which initially considers only a restricted subset of feasible paths for E-MAUs, denoted as $\hat{\mathcal{P}} \subseteq \mathcal{P}$, within the MP~\eqref{mp_model}--\eqref{mp_domain_r}. Our column generation algorithm proceeds iteratively
as follows. First, we begin with a restricted subset $\hat{\mathcal{P}}$, ensuring that the RMP has a feasible solution. We solve the RMP on this restricted set to obtain a primal solution and the corresponding dual solution. Then, the dual information obtained from solving the RMP is used as input to a PP to see if this primal solution is optimal for
the full problem by identifying whether paths with the negative reduced costs exist. The paths with the negative reduced costs can improve the quality of the solution. If existing, these newly identified paths with the lowest reduced costs are then added to the subset $\hat{\mathcal{P}}$, gradually expanding the solution space represented in the RMP. This iterative process continues until no paths with negative reduced costs remain. In this section, we introduce the RMP, the PP, and the customized label-correcting algorithm developed to efficiently solve the PP in detail.

\textbf{Restricted master problem.} At the $k$-th iteration of the C\&CG, the RMP is the LP relaxation of the MP, which is based on a subset of paths $\hat{\mathcal{P}}$ and continuous decision variables $\mathbf{x}$ and $\mathbf{r}$. We now can formulate the RMP as follows:
\begin{align}\label{Rmp_model}
{\rm{\big[RMP-}}k\big]\quad
\min\limits_{\mathbf{x}, \mathbf{r}, \mathbf{y}, \mathbf{b}} \quad & \theta_{1}(\sum\limits_{p\in\hat{\mathcal{P}}}c_p x_p + \sum\limits_{u\in\mathcal{U}_{fast}}m_u r_u) + \eta \tag{8a}\\
\mbox{s.t.}\quad\quad
& \eqref{mp_eta} - \eqref{mp_y}, \eqref{mp_y_b_domain}, \tag{8b}\\
&b_{q}^t(w) \leq  C \sum_{a\in\mathcal{A}_{q,t}}\sum_{p\in \hat{\mathcal{P}}}l_{p}^{a} x_p & \forall q\in\mathcal{Q}, t\in\mathcal{T},w\in \mathcal{W}, \label{rmp_capacity_mv} \tag{8c}\\
&\sum_{p\in\hat{\mathcal{P}}}\beta_{d,p}^{-} x_p=\sum_{p\in\hat{\mathcal{P}}}\beta_{d,p}^{+} x_p & \forall d\in\mathcal{D},\label{rmp_balance} \tag{8d}\\
&\sum_{a\in\mathcal{A}_{q,t}}\sum_{p\in \hat{\mathcal{P}}}l^a_p x_p\leq N^{com} & \forall q\in\mathcal{Q},t\in\mathcal{T},\label{eq:max_coposition_rmp} \tag{8e}\\
&\sum_{p\in\hat{\mathcal{P}}}g_{u,t,p} x_p \leq
    \begin{cases}
        r_u, & \text{if } u \in \mathcal{U}_{fast} \\
        N_d^{charge} , & \text{if } u \in \mathcal{U}  \setminus \mathcal{U}_{fast}
    \end{cases} & \forall t\in\mathcal{T},\label{eq:capacity_charge_station_rmp} \tag{8f}\\
& 0\leq x_p \leq 1 & \forall p\in \hat{\mathcal{P}}, \label{eq:domain_x_rmp} \tag{8g}\\
& r_u \in [0, N^{charge}_u] &\forall u\in\mathcal{U}_{fast}. \label{eq:domain_r_rmp} \tag{8h}
\end{align}

\textbf{Pricing problem.} 
In each iteration of the CG algorithm, the PP identifies new variables corresponding to paths with negative reduced costs. This process is based on the current values of the dual variables obtained from the RMP. Specifically, we define the dual variable for constraints (\ref{rmp_capacity_mv}) as $\varsigma_{qt}(w)$, for all $q \in \mathcal{Q}$, $t \in \mathcal{T}$, and $w \in \mathcal{W}$; the dual variable for constraints (\ref{rmp_balance}) as $\iota_{d}$ for all $d \in \mathcal{D}$; the dual variable for constraints (\ref{eq:max_coposition_rmp}) as $\pi_{qt}$, where $\forall q \in \mathcal{Q}, t \in \mathcal{T}$; and the dual variable for constraints (\ref{eq:capacity_charge_station_rmp}) as $\rho_{ut}$, where $\forall u \in \mathcal{U}$ and $t \in \mathcal{T}$. Using this dual information from the RMP, the reduced cost corresponding to path $p$ (denoted as $R_p$) is calculated as follows:
\begin{align}\label{reduced_cost_pp}
  R_p = \theta_{oper} c_p &- \sum_{d\in\mathcal{D}}(\beta_{d,p}^{-}-\beta_{d,p}^{+})\iota_d- \sum_{q\in\mathcal{Q}}\sum_{t\in\mathcal{T}}\sum_{a\in\mathcal{A}_{q,t}}l_p^a\pi_{qt}- \sum_{u\in\mathcal{U}}\sum_{t\in\mathcal{T}}g_{u,t,p}\rho_{ut}  +C\sum_{w\in\mathcal{W}}\sum_{q\in\mathcal{Q}}\sum_{t\in\mathcal{T}}\varsigma_{qt}(w)\sum_{a\in\mathcal{A}_{q,t}}l^a_p. \tag{9}
\end{align}

\textbf{Tailored label-correcting algorithm.} To find paths $p\in\mathcal{P}$ with the lowest reduced cost, we solve a time-dependent shortest path problem based on the space-time-SoC network with negative weights. This problem is efficiently addressed using a tailored label-setting algorithm, the details of which are provided in Algorithm~\ref{alg:spfa-specific} in Appendix \ref{sec:labelcorrecting}.

We denote the label of each node as $\Psi_n$. Given two nodes, $n$ and $v$, if there exists an arc $a = (n, v) \in \mathcal{A}$, the label update mechanism can be expressed as follows:
\begin{align}
\Psi_v&= \Psi_n + \theta_{oper} * c_{(n,v)} - \beta_{(n,v)}\iota_{d(n,v)}- \pi_{l(n)l(v)t(n)} - \rho_{l(v)t(v)} + C\sum_{w\in\mathcal{W}}\varsigma_{l(n)l(v)t(n)}(w). \tag{10}
\end{align}
The term $\beta_{(n,v)}$ is defined as follows: If $(n,v)\in\cup_{d\in\mathcal{D}}\mathcal{A}_{d}^{in}$, then $\beta_{(n,v)}=1$; if $(n,v)\in\cup_{d\in\mathcal{D}}\mathcal{A}_{d}^{out}$, then $\beta_{(n,v)}=-1$; otherwise, $\beta_{(n,v)}=0$. For clarity, dual variables are assumed to be zero when undefined within their respective domains. For example, the dual variable $\pi_{(n,v,t)}$ is defined only on arcs in $\mathcal{A}_{travel}$, implying $\pi_{(n,v,t)} = 0$ for all $a \in \mathcal{A} \setminus \mathcal{A}_{travel}$. 

\subsection{Network downsizing and arc aggregation}
\label{sec:aggregating}

The key to efficiently solving the PP in large-scale instances lies in reducing the scale of the space-time-SoC network. Rather than explicitly modeling all possible connections between nodes across the three dimensions of SoC, time and space, redundant nodes and arcs are both removed. Specifically, to eliminate redundant nodes and arcs, we first determine the theoretically possible maximum and minimum values of SoC at each station when every unit arrives. Similarly, we identify the theoretically possible earliest and latest arrival times of units at each station. Both downsizing methods are derived by defining the concept of a \textit{subpath} and are proven to preserve the feasible region of solutions. In addition, by aggregating travel arcs that enable a unit to execute consecutively without involving decoupling and coupling possibilities to \textit{super travel arcs}, the number of arcs within the space-time-SoC network is greatly reduced. The downsizing and aggregating process is carried out through the following three-step procedure.

\textbf{First-stage downsizing from the SoC dimension.} At each station $s$, for the unit that executes a service with path $p$, the minimum and maximum SoC values must be within a certain range. To rigorously derive this range, we define the subpaths, which serve as the foundation for proving Proposition~\ref{max_min_soc}. The goal is to eliminate redundant nodes and arcs in the space-time-SoC network along the SoC dimension for each discrete SoC, as detailed in Proposition~\ref{max_min_soc}.

 \begin{definition} [Subpath]
(i) A charging-station subpath $\Tilde{p}$ is defined by a node sequence $\Tilde{p}=\{n_{\Tilde{p}}^1,n_{\Tilde{p}}^2,...,n_{\Tilde{p}}^m\}$ with starting node $n_{\Tilde{p}}^1\in\mathcal{N}^{charge}\cup\mathcal{N}^{depot}$ and ending node $n_{\Tilde{p}}^m\in\mathcal{N}^{station}$. Any two adjacent nodes in the subpath satisfy $(n_{\Tilde{p}}^{j-1},n_{\Tilde{p}}^{j})\in\mathcal{A}$. We use $\textit{CS}_s$ to denote the set of all charging-station subpaths ending at the station $s$. 

(ii) A station-charging subpath $\Tilde{p}$ is defined by a node sequence $\Tilde{p}=\{n_{\Tilde{p}}^1,n_{\Tilde{p}}^2,...,n_{\Tilde{p}}^m\}$ with starting node $n_{\Tilde{p}}^1\in\mathcal{N}^{station}$ and ending node $n_{\Tilde{p}}^m\in\mathcal{N}^{charge}\cup\mathcal{N}^{depot}$. Any two adjacent nodes in the subpath satisfy $(n_{\Tilde{p}}^{j-1},n_{\Tilde{p}}^{j})\in\mathcal{A}$. We use $\textit{SC}_s$ to denote the set of all station-charging subpaths starting at station $s$. For each $s\in\mathcal{S}$, the corresponding subpath sets $\textit{CS}_s$ and $\textit{SC}_s$ can be enumerated.
\end{definition}

 \begin{proposition}[Bounds on SoCs at station nodes]
 \label{max_min_soc}
Let $(\mathbf{x},\bm{\delta},\bm{\zeta})$ denote a feasible solution to model (\ref{model:MILP}). For each station node $n$ located at station $s$ along path $p$ where $x_p=1$, the following condition holds: $\min\limits_{\Tilde{p} \in \textit{SC}_s} \left\{  e_{\text{min}}+ \sum_{j=1}^{m-1} \theta^{\text{soc}}_{(n_{\Tilde{p}}^{j}, n_{\Tilde{p}}^{j+1})} \right\}\leq e(n)\leq \max\limits_{\Tilde{p} \in \textit{CS}_s} \left\{  e_{\text{max}}- \sum_{j=1}^{m-1} \theta^{\text{soc}}_{(n_{\Tilde{p}}^{j}, n_{\Tilde{p}}^{j+1})} \right\}$, where $e(n)$ represent the SoC value at node $n$, and $\Tilde{p}$ is the subpath either ending or starting at station $s$.
 \end{proposition}

\proof{Proof.} See Appendix \ref{sec:proof}.

\textbf{Second-stage downsizing from the time dimension.} At each station $s$, for a unit executing a servcie along path $p$, the earliest and latest arrival times must fall within a specific range. To rigorously derive this range and thus removing redundant station nodes and the corresponding arcs, we develop Proposition~\ref{time} based on the subpaths.

 \begin{proposition}[Bounds on arrival times at station nodes]
 \label{time}
Let $(\mathbf{x},\bm{\delta},\bm{\zeta})$ denote a feasible solution to model (\ref{model:MILP}). For each station node $n$ located at station $s$ along path \( p \) where $x_p=1$, the following condition holds: $\min\limits_{\Tilde{p} \in \textit{CS}_s} \left\{  \delta+ \sum_{j=1}^{m-1} \chi_{(n_{\Tilde{p}}^{j}, n_{\Tilde{p}}^{j+1})} \right\}\leq t(n)\leq \max\limits_{\Tilde{p} \in \textit{SC}_s} \left\{  \left|\mathcal{T}\right|- \sum_{j=1}^{m-1} \chi_{(n_{\Tilde{p}}^{j}, n_{\Tilde{p}}^{j+1})} \right\}$, where $t(n)$ represent the time at node $n$, and $\Tilde{p}$ is the subpath either ending or starting at station $s$.
 \end{proposition}
\proof{Proof.} See Appendix \ref{sec:proof}.

\textbf{Third-stage downsizing from the space dimension through aggregation.} For each path $p$, we aggregate the travel arcs between two adjacent facilities into a single super travel arc. These facilities can either be a depot and a charging station, or two successive charging stations. An illustration of constructing super travel arcs is given in Example \ref{example:superarcs} in Appendix \ref{sec:example}. This super travel arc represents a direct connection between a depot and the nearest charging station or between two adjacent charging stations, defined as follows:

\begin{definition} [Super travel arc]
We denote the super travel arc as \( \check{a} = \{(n_0, n_k) \mid n_0, n_k \in \mathcal{N}^{station}\} \) and the set of all super travel arcs as \( \check{\mathcal{A}}_{travel} \). The corresponding sequence of normal travel arcs is defined as 
\[
\mathcal{A}_{\check{a}} = \{a_0, a_1, \dots, a_k \mid a_j \in \mathcal{A}_{travel}, j \leq k\} = \{(n_0, n_1), (n_1, n_2), \dots, (n_{k-1}, n_k)\}.
\]
For any two adjacent arcs in \( \mathcal{A}_{\check{a}} \), the following hold: (i) $t_j = t_{j-1} + \chi_{n_{j-1}, n_j, t_{j-1}}$; (ii) $e_j = e_{j-1} - \theta^{soc}_{(n_{j-1}, n_j)}$. Here, \( t_j \) and \( e_j \) represent the updated time and SoC, respectively. We replace the travel arcs in the space-time-Soc network $\mathcal{G}$ (Section \ref{section:network_constraction}) with super travel arcs, keeping all other elements unchanged, to obtain the new network $\check{\mathcal{G}}$. We denote $\check{q}$ in $\check{\mathcal{G}}$ as the super section between a depot and the nearest charging station or between two adjacent charging stations. 
\end{definition}

\begin{proposition}\label{pp4}
    Let $opt[\text{PP}(\mathcal{G})]$ be the optimal objective value of pricing problem defined in Section \ref{section_cg} with space-time-Soc network $\mathcal{G}$. It always holds that 
    $opt[\text{PP}(\mathcal{G})]=opt[\text{PP}(\check{\mathcal{G}})]$.
\end{proposition}
\proof{Proof.} See Appendix \ref{sec:proof}.

\subsection{Outer approximation algorithm}
\label{sec:Outer}

Recall that the SP in the C\&CG framework, i.e., $\text{DF}(\hat{\mathbf{x}})$ \eqref{bilinear_obj_1} - \eqref{df_3_1}, is formulated to identify the worst-case scenario $\hat{\bm{\zeta}}_k$. While the SP is solved using GUROBI, preliminary experiments suggest that there is certain room for improving the efficiency of the solution process. We propose an outer approximation algorithm to provide high-quality initial solutions, which enhance the solution efficiency of the GUROBI. The details of this approach are outlined in Appendix \ref{procedure_outer}.

\subsection{Overall framework of the algorithm combining C\&CG and CG}
\label{sec:framework}

Lastly, we design two procedures to solve the E-RCRSP, with their overall frameworks presented in Figure~\ref{fig:algorithms} in Appendix \ref{sec:overallFlowchart}.

$\bullet$ \textit{Always-integer procedure:} In each iteration of the C\&CG framework, an integer solution of the MP is obtained and passed to the SP for evaluation. To derive this integer solution, CG is first used to generate a fractional solution that satisfies the termination conditions and identifies the corresponding set of paths. Then, GUROBI is employed to solve the MP with this fixed set of paths, producing the integer solution. For the SP, the outer approximation algorithm is designed to obtain high-quality initial solutions, which are then passed to GUROBI to produce the optimal solution. If the termination criteria of the C\&CG algorithm are not met, the SP identifies the worst-case scenario and feeds it back to the MP. The partial enumeration is expanded, and the recourse variables and corresponding constraints are added to the MP. The algorithm terminates when the optimality gap falls below a predefined threshold.

$\bullet$ \textit{First-continuous procedure:} The LP relaxation of the MP is solved, and its solution is passed to the SP. This iterative process continues until the termination criteria of the C\&CG algorithm are met—specifically, when the gap between the best upper and lower bounds is less than a predefined threshold. Subsequently, the set of paths is fixed, and the C\&CG algorithm is restarted to find the integer solution. This process repeats until the termination criteria of the C\&CG algorithm are satisfied. The detailed framework of this procedure is presented in Figure~\ref{fig:algorithm_overall_first} in Appendix~\ref{sec:detailedFlowchart}.

\section{Numerical experiments}
\label{sec:results}

In this section, we evaluate the performance of our proposed model and solution methodologies through case studies constructed using real-world data. The data include the operational information of the bus line and passenger flows, as well as the real-life technical specifications of advanced MAVs provided by our industry partner NExT Future Transportation Inc. Section \ref{sec:instances} introduces the dataset, explains how instances with varying numbers of potential charging stations and discretized time intervals are constructed, and presents the resulting instances. In Section \ref{sec:sensitivity}, we analyze the impact of weight coefficients in the objective function on operators' and passengers' costs, thereby determining a reasonable setting for these weights in subsequent experiments. Section \ref{sec:benefitsAlgorithm} examines the benefits of the proposed algorithm and downsizing methodologies. The practical advantages of employing E-MAVs with flexible compositions are assessed in Section \ref{sec:benefitFlexible}. Finally, in Section \ref{sec:entireDay}, we solve an instance with entire-day operations to demonstrate the scalability of our algorithm and provide insights for practical implementation.

The proposed algorithm is implemented in Python 3.9.13, and the subproblem in the C\&CG framework is solved using GUROBI v9.5.1. All computational experiments are performed on a laptop featuring a 12th-generation Intel i9 processor at 2.2 GHz and with 64 GB of random access memory and running the Windows operating system.

\subsection{Instances}
\label{sec:instances}

\begin{figure}[tbp]
    \centering
    \includegraphics[width=0.9\linewidth]{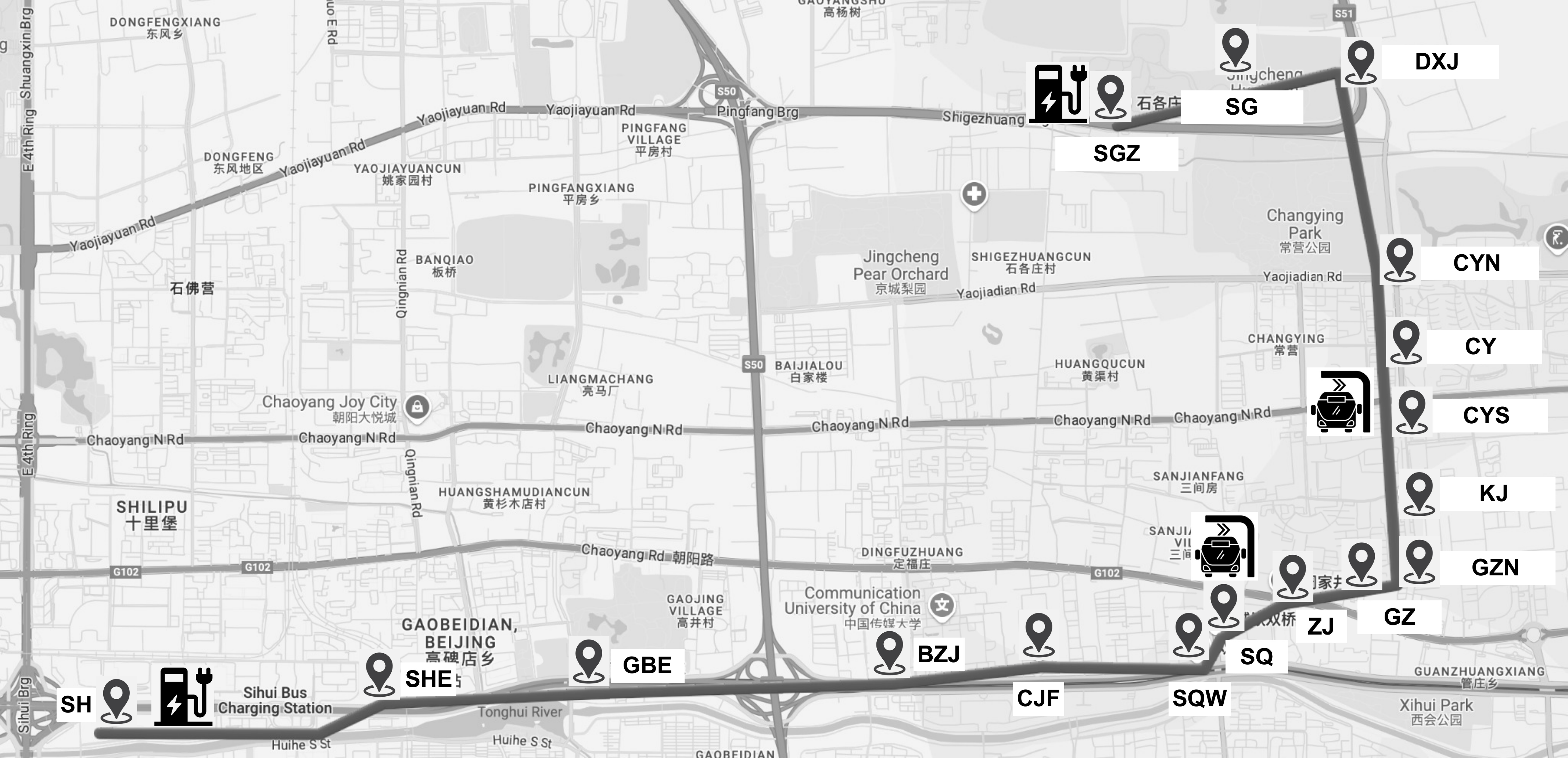}
    \caption{Real-life bus line used in the computational study.}
    \label{fig:busLine}
\end{figure}

The instances used in the experiments are based on Line 506 of the Beijing bus network. As shown in Figure~\ref{fig:busLine}, each direction of this line includes 17 stations. Two depots, located near terminal stations SH and SGZ, are equipped with 10 standard charging posts and support decoupling and coupling operations. Stations SQ and CYS enable short-turning, decoupling, and coupling operations and are potential charging stations. Table~\ref{tab:instances} summarizes all instances studied in this paper, detailing the number of potential charging stations, the start and end times, and the number of discretized time intervals considered in each instance.

\begin{table}[tbp]
\centering
\caption{Characteristics of the instances.}
\label{tab:instances}
\begin{tabular}{ccccc}
\Xhline{1pt}
Instance  & \begin{tabular}[c]{@{}c@{}}Number of potential\\ charging stations\end{tabular} & Start time & End time & \begin{tabular}[c]{@{}c@{}}Number of discretized\\  time intervals \end{tabular}\\ 
\Xhline{0.6pt}
A\_1  & 1  & 07:00 & 08:30 & 90 \\
A\_2  & 2  & 07:00 & 08:30 & 90\\
B\_1  & 1  & 07:00 & 09:00 & 120\\
B\_2  & 2  & 07:00 & 09:00 & 120\\
C\_1  & 1  & 07:00 & 09:30 & 150\\
C\_2  & 2  & 07:00 & 09:30 & 150\\
D\_1  & 1  & 07:00 & 10:00 & 180\\
D\_2  & 2  & 07:00 & 10:00 & 180\\
E\_1  & 1  & 07:00 & 11:00 & 240\\
E\_2  & 2  & 07:00 & 11:00 & 240\\
F\_1  & 1  & 07:00 & 12:00 & 300\\
F\_2  & 2  & 07:00 & 12:00 & 300\\
G\_1  & 2  & 07:00 & 14:00 & 420\\
H\_1  & 2  & 07:00 & 22:00 & 900\\
\Xhline{1pt}
\end{tabular}
\end{table}

The parameter settings for all experiments are presented in Table~\ref{tab:parameters}. Parameters related to running times on segments and passenger demand are based on real-life data provided by the operators. Parameters related to E-MAUs, such as battery capacity, idling electricity consumption, operational electricity consumption, maximum compositions of an E-MAV, and the capacity of each E-MAU, are derived from the technical specifications of NExT MAVs designed by NExT Future Transportation Inc. Some of these specifications are publicly available on this company’s official website, while others were directly provided by the company.

\begin{table}[h]
\centering
\caption{Parameter settings for all experiments.}
\label{tab:parameters}
\resizebox{\textwidth}{!}{
\begin{tabular}{@{}llll@{}}
\Xhline{1pt}
\textbf{Parameter} & \textbf{Value} & \textbf{Parameter} & \textbf{Value} \\
\Xhline{0.6pt}
Range of SoC & [20\%, 80\%] & Battery capacity$^{a}$ & 90 kWh \\
Length of unit time & 1 min & Speed of the vehicle & 30 km/h \\
Idling electricity consumption & 0.042 kW/min$^{a}$ & Operational electricity consumption & 0.25 kW/km$^{a}$ \\
Charging power of fast-charging stations & 600 kW $^{b}$ & Charging power of depots & 50 kW $^{b}$\\
Charging capacity of each depot& 10 & Maximum charging capacity of each fast-charging station  & 2  \\
Maximum compositions of E-MAVs$^{c}$ & 5 &Capacity of an E-MAU & 16 passengers $^{c}$ \\
Cost of installing a charging post & 200,000 \euro $^{b}$ &Cost of purchasing an E-MAU & 50,000 \euro  \\
Fixed charging cost & 10 \euro $^{d}$ & Electricity cost$^{d}$ & 0.1361 \euro/kWh $^{c}$ \\
Length of unit SoC & 1 \% & Number of stations  & 34  \\
\Xhline{1pt}
\end{tabular}
}
\begin{tablenotes}
        \footnotesize
        \raggedright
        \item[a] Parameters given by NExT Future Transportation Inc. \\
        \item[b] Parameters borrowed from \cite{shadi2022}. \\
        \item[c] Parameters borrowed from \cite{NextCapacity}.\\
        \item[d] Parameters borrowed from \cite{Rolf2023}.\\
    \end{tablenotes}
\end{table}

\subsection{Sensitivity analysis of objective function weights}
\label{sec:sensitivity}

Recall that the objective function of model \eqref{formulation:U} includes weighting coefficients for operators' and passengers' costs. The balance between these weights significantly influences the benefits for both parties. Consequently, for bus managers, the trade-off between passengers' costs and operators' costs in the objective function is a key factor to consider. In this subsection, we use instance G\_1 as the basis for 11 experiments with various weighting coefficients to investigate this trade-off. To facilitate comparison, we vary $\theta_1$ and $\theta_2$, while keeping $\theta_3$ set to ten times $\theta_2$ to penalize the cost of unserved passengers. The proposed solution method, combining C\&CG and CG, is applied to solve these experiments.

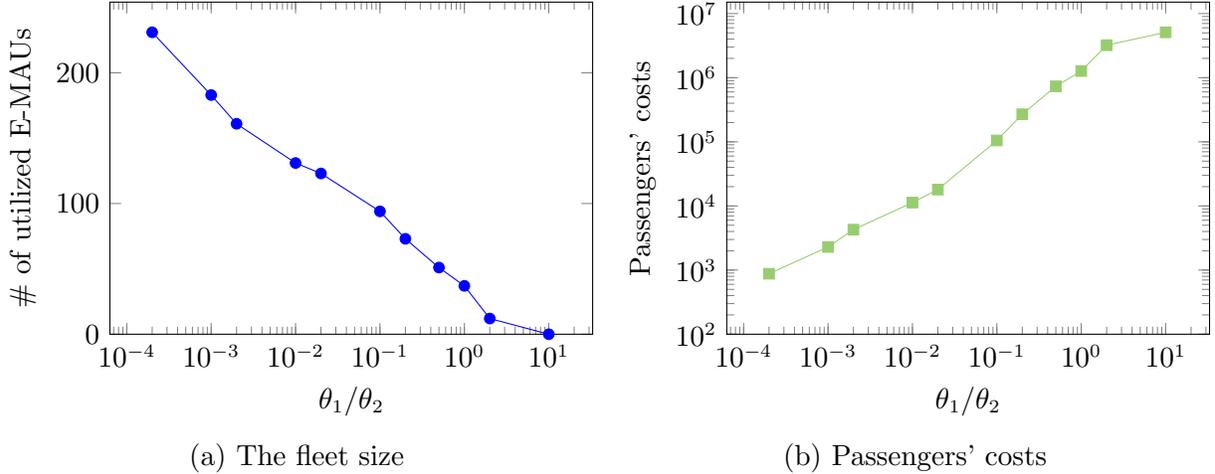
\begin{figure}[htbp]
    \centering
\begin{subfigure}[b]{0.49\textwidth}
    \centering
\begin{tikzpicture}
  \begin{axis}[
    xlabel={$\theta_1/\theta_2$},
    ylabel={\# of utilized E-MAUs},
    width=8cm,
    height=6cm,
    ymin=0,
     xmode=log, 
    log basis x={10}, 
    xmin=10^(-4.2),
    legend style={at={(0.5,-0.15)},anchor=north,legend columns=-1}
  ]
    \addplot[
      mark=*,
      color=blue
    ]
    table[col sep=comma, x=Weight1, y=FleetSize] {gainsFleet.csv};
  \end{axis}
\end{tikzpicture}
\caption{The fleet size}
\end{subfigure}
\begin{subfigure}[b]{0.49\textwidth}
    \centering
    \begin{tikzpicture}
    \begin{axis}[
      xlabel={$\theta_1/\theta_2$},
      ylabel={Passengers' costs},
      width=8cm,
      height=6cm,
      ymode=log, 
    log basis y={10}, 
    ymin=1e2, 
     xmode=log, 
    log basis x={10}, 
    xmin=10^(-4.2),
      legend style={at={(0.5,-0.15)},anchor=north,legend columns=-1}
    ]
      \addplot[
        mark=square*,
        color=YellowGreen
      ]
      table[col sep=comma, x=Weight3, y=Passenger] {gainPassenger.csv};
    \end{axis}
  \end{tikzpicture}
\caption{Passengers' costs}
\end{subfigure}
\caption{Impact of weighting coefficients on the fleet size and passengers' costs.}
\label{fig:weightLog}
\end{figure}

Figure \ref{fig:weightLog} reports the optimized fleet size and passengers' costs. The horizontal axis in Figure \ref{fig:weightLog}(a) and both axes in Figure \ref{fig:weightLog}(b) are plotted on a logarithmic scale, with values displayed in their original scale. As \( \theta_1/\theta_2 \) increases, we observe that the fleet size rapidly decreases, prioritizing lower operators' costs, while passengers' costs rise sharply, indicating reduced service quality. At \( \theta_1/\theta_2 = 1/10 \), a favorable balance is observed: the fleet size remains adequate to ensure reasonable service quality, and passengers' costs are manageable without overburdening operators. From these findings, we can derive a conclusion that as more attention is paid to operator costs (i.e., $\theta_1/\theta_2$ are getting bigger), the service quality is decreased due to a reduction in operated units. Thus, \( \theta_1/\theta_2 = 1/10 \) is selected for the following experiments as it offers a practical trade-off between operators' and passengers' interests.

\subsection{Performances of solution methodologies}
\label{sec:benefitsAlgorithm}

In this subsection, we first investigate the effectiveness of the proposed three methodologies for downsizing the space-time-SoC network, which lays the foundation for their applications in subsequent experiments. Second, we compare the computational performance of our algorithm, which combines C\&CG and CG, with the state-of-the-art solution method that integrates Benders decomposition and CG for solving the relaxed path-based model. Lastly, we evaluate the computational efficiency and solution quality of our algorithm on the path-based model compared to the arc-based model, aiming to highlight the benefits of the path-based formulation.

\subsubsection{Benefits of downsizing the network.}
\label{sec:downsizing}

We first test the effectiveness of our proposed three methodologies for downsizing the space-time-SoC network. Two benchmarks are considered: the first applies the algorithm without any downsizing methodology, and the second incorporates the first two downsizing methodologies, focusing on the time and SoC dimensions. In contrast, our algorithm employs all three downsizing methodologies, including those for the time and SoC dimensions, along with the super travel arcs. To facilitate comparison, without loss of generality, we use the model \eqref{formulation:D} presented in Section \ref{sec:determinsticModel}, setting $\alpha_{q}^t=\overline{\alpha}_{q}^t+ 0.5\cdot \Tilde{\alpha}_{q}^t$. We solve the relaxation of this model using the column generation algorithm introduced in Appendix \ref{sec:cg_appendix}. It should be noted that the space-time-SoC network used to solve the relaxation of this model is identical to the one employed in the robust optimization model. The employed model is solved to optimality in this series of experiments.

\begin{insight}
The proposed three methodologies for downsizing the space-time-SoC network can reduce the network size by up to 86.28\%, resulting in a decrease in the average computational time for solving the pricing problem by at least 80.20\% and up to 94.97\%. Moreover, these benefits are highly robust across various scales of problems.
\end{insight}

\begin{table}[h]
\centering
\caption{Performance comparison among the algorithms without any downsizing methodology (LC), with two methods (LC\_R), and with three downsizing techniques (LC\_RS).}\label{Algorithm:label_correcting}
\resizebox{\textwidth}{!}{
\begin{tabular}{llrrrrrr}
\Xhline{1pt}
Instance & Method & \makecell{\# of \\ nodes} & \makecell{\# of \\ arcs} & \makecell{Objective \\ value} & \makecell{Average \\ PP time (s)} & \makecell{\% reduction in \\ network size} & \makecell{\% reduction in \\ computational time} \\
\Xhline{0.6pt}
A\_1 & LC     & $2.09 \times 10^5$ & $2.71 \times 10^5$ & $2.18 \times 10^5$ & 0.60  & -     & -     \\
     & LC\_R  & $4.98 \times 10^4$ & $7.00 \times 10^4$ & $2.18 \times 10^5$ & 0.11  & 74.21 & 81.38 \\
     & LC\_RS & $1.70 \times 10^4$ & $3.72 \times 10^4$ & $2.18 \times 10^5$ & 0.03  & 86.28 & 94.97 \\
     \rule{0pt}{15pt}
A\_2 & LC     & $2.20 \times 10^5$ & $3.14 \times 10^5$ & $2.16 \times 10^5$ & 0.94  & -     & -     \\
     & LC\_R  & $7.27 \times 10^4$ & $1.07 \times 10^5$ & $2.16 \times 10^5$ & 0.20  & 65.94 & 79.04 \\
     & LC\_RS & $3.01 \times 10^4$ & $6.45 \times 10^4$ & $2.16 \times 10^5$ & 0.08  & 79.48 & 91.76 \\
     \rule{0pt}{15pt}
B\_1 & LC     & $2.78 \times 10^5$ & $3.65 \times 10^5$ & $7.31 \times 10^5$ & 1.95  & -     & -     \\
     & LC\_R  & $1.11 \times 10^5$ & $1.53 \times 10^5$ & $7.31 \times 10^5$ & 0.64  & 58.16 & 67.28 \\
     & LC\_RS & $3.50 \times 10^4$ & $7.66 \times 10^4$ & $7.31 \times 10^5$ & 0.18  & 79.01 & 90.60 \\
     \rule{0pt}{15pt}
B\_2 & LC     & $2.93 \times 10^5$ & $4.22 \times 10^5$ & $7.05 \times 10^5$ & 3.10  & -     & -     \\
     & LC\_R  & $1.40 \times 10^5$ & $2.07 \times 10^5$ & $7.05 \times 10^5$ & 1.12  & 51.06 & 63.93 \\
     & LC\_RS & $5.50 \times 10^4$ & $1.22 \times 10^5$ & $7.05 \times 10^5$ & 0.39  & 71.11 & 87.37 \\
     \rule{0pt}{15pt}
C\_1 & LC     & $3.48 \times 10^5$ & $4.58 \times 10^5$ & $1.31 \times 10^6$ & 4.00  & -     & -     \\
     & LC\_R  & $1.68 \times 10^5$ & $2.32 \times 10^5$ & $1.31 \times 10^6$ & 1.80  & 49.47 & 54.98 \\
     & LC\_RS & $5.19 \times 10^4$ & $1.15 \times 10^5$ & $1.31 \times 10^6$ & 0.53  & 74.85 & 86.79 \\
     \rule{0pt}{15pt}
C\_2 & LC     & $3.66 \times 10^5$ & $5.30 \times 10^5$ & $1.21 \times 10^6$ & 5.90  & -     & -     \\
     & LC\_R  & $2.04 \times 10^5$ & $3.04 \times 10^5$ & $1.21 \times 10^6$ & 2.80  & 42.57 & 52.55 \\
     & LC\_RS & $7.91 \times 10^4$ & $1.80 \times 10^5$ & $1.21 \times 10^6$ & 0.93  & 66.13 & 84.29 \\
     \rule{0pt}{15pt}
D\_1 & LC     & $4.17 \times 10^5$ & $5.52 \times 10^5$ & $1.84 \times 10^6$ & 11.89 & -     & -     \\
     & LC\_R  & $2.24 \times 10^5$ & $3.11 \times 10^5$ & $1.84 \times 10^6$ & 3.49  & 43.69 & 70.62 \\
     & LC\_RS & $6.83 \times 10^4$ & $1.54 \times 10^5$ & $1.84 \times 10^6$ & 0.99  & 72.00 & 91.70 \\
     \rule{0pt}{15pt}
D\_2 & LC     & $4.39 \times 10^5$ & $6.38 \times 10^5$ & $1.72 \times 10^6$ & 13.28 & -     & -     \\
     & LC\_R  & $2.68 \times 10^5$ & $4.01 \times 10^5$ & $1.72 \times 10^6$ & 4.81  & 37.16 & 63.79 \\
     & LC\_RS & $1.03 \times 10^5$ & $2.36 \times 10^5$ & $1.72 \times 10^6$ & 1.63  & 63.00 & 87.69 \\
     \rule{0pt}{15pt}
E\_1 & LC     & $5.56 \times 10^5$ & $7.39 \times 10^5$ & $2.89 \times 10^6$ & 14.14 & -     & -     \\
     & LC\_R  & $3.40 \times 10^5$ & $4.73 \times 10^5$ & $2.89 \times 10^6$ & 8.03  & 35.94 & 43.20 \\
     & LC\_RS & $1.03 \times 10^5$ & $2.36 \times 10^5$ & $2.89 \times 10^6$ & 2.05  & 68.04 & 85.47 \\
     \rule{0pt}{15pt}
E\_2 & LC     & $5.86 \times 10^5$ & $8.54 \times 10^5$ & $2.78 \times 10^6$ & 24.73 & -     & -     \\
     & LC\_R  & $3.95 \times 10^5$ & $5.95 \times 10^5$ & $2.78 \times 10^6$ & 10.97 & 30.30 & 55.63 \\
     & LC\_RS & $1.50 \times 10^5$ & $3.51 \times 10^5$ & $2.78 \times 10^6$ & 3.32  & 58.95 & 86.56 \\
     \rule{0pt}{15pt}
F\_1 & LC     & $6.95 \times 10^5$ & $9.27 \times 10^5$ & $3.04 \times 10^6$ & 21.38 & -     & -     \\
     & LC\_R  & $4.68 \times 10^5$ & $6.52 \times 10^5$ & $3.04 \times 10^6$ & 12.45 & 29.63 & 41.78 \\
     & LC\_RS & $1.40 \times 10^5$ & $3.24 \times 10^5$ & $3.04 \times 10^6$ & 3.21  & 65.09 & 84.97 \\
     \rule{0pt}{15pt}
F\_2 & LC     & $7.32 \times 10^5$ & $1.07 \times 10^6$ & $3.02 \times 10^6$ & 28.68 & -     & -     \\
     & LC\_R  & $5.38 \times 10^5$ & $8.10 \times 10^5$ & $3.02 \times 10^6$ & 17.88 & 24.32 & 37.65 \\
     & LC\_RS & $2.03 \times 10^5$ & $4.76 \times 10^5$ & $3.02 \times 10^6$ & 5.68  & 55.60 & 80.20 \\
\Xhline{1pt}
\end{tabular}%
}
\end{table}

Table \ref{Algorithm:label_correcting} presents the results of the algorithms without any downsizing methodology (denoted as LC), with two downsizing methodologies (denoted as LC\_R), and with three downsizing methodologies (denoted as LC\_RS). This table illustrates the objective values, differences in network size, and computational times for solving the pricing problem between LC\_R and LC, as well as between LC\_RS and LC. From the results in Table \ref{Algorithm:label_correcting}, we can observe that for smaller instances (i.e., A\_1 and A\_2), LC\_RS achieves reductions in network size of 86.28\% and 79.48\%, respectively, compared to LC, and computational time reductions of 94.97\% and 91.76\%, respectively. This highlights the significant impact of incorporating all three downsizing methodologies, even for relatively small problems. The super travel arcs introduced in LC\_RS provide additional benefits compared to LC\_R, as evidenced by the reductions in computational time (e.g., from 81.38\% for LC\_R to 94.97\% for LC\_RS in A\_1).

The other observation is that as the scale of the problem increases, the advantages of LC\_RS become even more pronounced. For example, in instance E\_2, LC\_RS reduces the network size by 58.95\% and computational time by 86.56\%, compared to LC. Similarly, for the largest instance F\_2, LC\_RS achieves reductions of 55.60\% in network size and 80.20\% in computational time. This robust performance across various problem sizes show the scalability and efficiency of LC\_RS. From these findings, we can derive a conclusion that LC\_RS is highly effective in significantly reducing the network size and computational time while maintaining the solution quality.

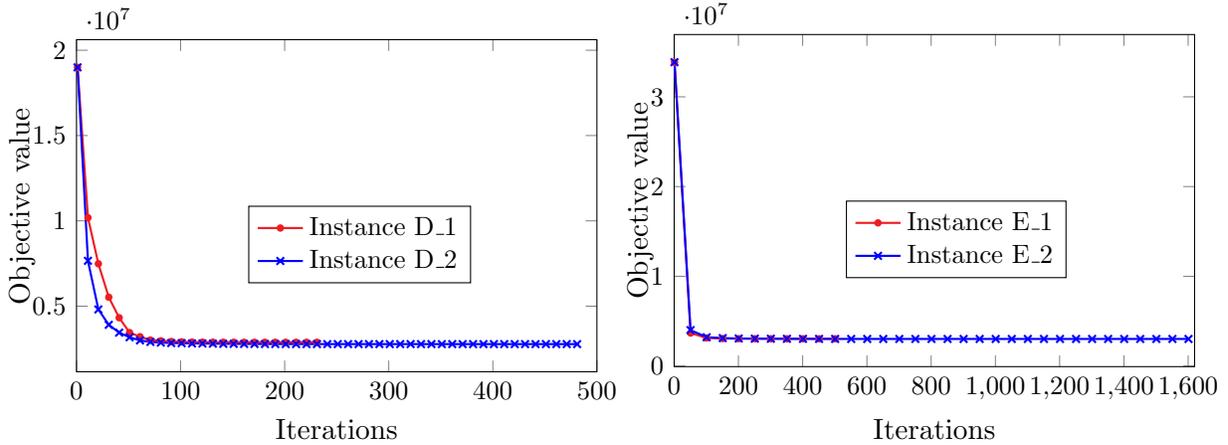
\begin{figure}[h]
    \centering
\begin{subfigure}[b]{0.49\textwidth}
    \centering
\begin{tikzpicture}
    \centering
\begin{axis}
[xlabel={Iterations},
  ylabel={Objective value},every axis plot/.append style={thick}, legend style={font=\small, at={(0.33,0.50)},anchor=north west},legend cell align={left},     width=8.5cm,
    height=6cm, xmin=0, xmax=500, ylabel near ticks,  ylabel style={yshift=-0.3cm}, tick label style={font=\small}  ]
\addplot+[Red, mark=*, mark size=1pt, mark options={solid}] table [x=Iterations, y=D1, col sep=comma] {gains.csv};
\addlegendentry{Instance D\_1}
\addplot+[blue, mark=x, mark size=2pt,mark options={fill=blue}] table [x=Iterations, y=D2, col sep=comma] {gains.csv};
\addlegendentry{Instance D\_2}
\end{axis}
\end{tikzpicture}
\caption{Instances D\_1 and D\_2}
\end{subfigure}
\begin{subfigure}[b]{0.49\textwidth}
    \centering
\begin{tikzpicture}
    \centering
\begin{axis}
[xlabel={Iterations},
  ylabel={Objective value},every axis plot/.append style={thick}, legend style={font=\small, at={(0.33,0.50)},anchor=north west},legend cell align={left},     width=8.5cm,
    height=6cm, xmin=0, xmax=1620,  ylabel style={yshift=-0.3cm}, tick label style={font=\small} ]
\addplot+[Red, mark=*, mark size=1pt, mark options={solid}] table [x=Iterations, y=E1, col sep=comma] {gains2.csv};
\addlegendentry{Instance E\_1}
\addplot+[blue, mark=x, mark size=2pt,mark options={fill=blue}] table [x=Iterations, y=E2, col sep=comma] {gains2.csv};
\addlegendentry{Instance E\_2}
\end{axis}
\end{tikzpicture}
\caption{Instances E\_1 and E\_2}
\end{subfigure}
\caption{Convergence curves of objective values among various instances.}
\label{fig:Trucked}
\end{figure}

We further analyze the evolution of the objective value throughout the solving process. Figure \ref{fig:Trucked} illustrates the convergence curves of instances D\_1, D\_2, E\_1, and E\_2. It can be observed that the objective function values decrease very fast at the beginning of iterations. However, as the iteration progresses, the rate of decrease slows. Besides, before convergence, the algorithm undergoes a prolonged iteration process during which the decrease in the objective values become minimal. This indicates a pronounced tailing effect in the column generation process in our problem. Based on these observations, introducing a truncation strategy in subsequent experiments could greatly enhance computational efficiency with minimal impact on solution quality.

\subsubsection{Performace of the C\&CG+CG algorithm.}
In this subsection, we use the classical Benders decomposition (BD) as a benchmark to evaluate the performance of the proposed C\&CG algorithm in solving the linear relaxation problem. The detailed introduction of the algorithm combining BD and CG is provided in Appendix \ref{sec:BD}. Based on the conclusions in Section \ref{sec:downsizing}, the CG algorithm is terminated if the objective function value of the RMP decreases by no more than 0.5\% for 10 consecutive iterations. For both BD+CG and C\&CG+CG, the termination conditions are set such that the algorithm stops if either the optimality gap is no greater than 3\% or the computational time exceeds 3600 seconds from the start of the iteration. The optimality gap is defined as $(\text{UB}-\text{LB})/\text{UB}\times 100\%$, where UB and LB represent the upper and lower bounds defined in Algorithm \ref{alg:ccg}.

\begin{table}[tbp]
    \centering
    \caption{Comparison of upper bound, lower bound, iterations, optimality gap, and computational time between the C\&CG+CG and the BD+CG. }
    \label{tab:BD}
    \begin{threeparttable}
    \begin{tabular}{clrrrrr}
        \Xhline{1pt}
        Instance & Method & \makecell{Upper \\ Bound} & \makecell{Lower \\ Bound} & \makecell{\# of \\ iteration} & \makecell{Optimality \\ gap (\%)}  & \makecell{Computational \\ time (s)}   \\
        \Xhline{0.6pt}
        \multirow{2}{*}{A\_1} & BD+CG & 1.10E+05 & 1.07E+05 & 26 & 2.83 & 103.62 \\
                              & C\&CG+CG & 1.10E+05 & 1.10E+05 & \textbf{2} & \textbf{0} & \textbf{4.38} \\
        \multirow{2}{*}{A\_2} & BD+CG & 1.11E+05 & 1.08E+05 & 34 & 2.92 & 307.04 \\
                              & C\&CG+CG & 1.10E+05 & 1.10E+05 & \textbf{2} & \textbf{0} & \textbf{7.26} \\
        
        \multirow{2}{*}{B\_1} & BD+CG & 9.31E+04 & 8.20E+04 & 122 & 11.96 & 3655.69 \\
                              & C\&CG+CG & 8.69E+04 & 8.69E+04 & \textbf{2} & \textbf{0} & \textbf{40.86} \\
      
        \multirow{2}{*}{B\_2} & BD+CG & 8.89E+04 & 7.61E+04 & 99 & 14.4 & 3616.77 \\
                              & C\&CG+CG & 8.31E+04 & 8.31E+04 & \textbf{2} & \textbf{0} & \textbf{67.86} \\
    
        \multirow{2}{*}{C\_1} & BD+CG & 2.13E+05 & 1.05E+05 & 50 & 50.82 & 3663.74 \\
                              & C\&CG+CG & 1.47E+05 & 1.47E+05 & \textbf{2} & \textbf{0} & \textbf{66.85} \\
      
        \multirow{2}{*}{C\_2} & BD+CG & 1.67E+05 & 1.07E+05 & 63 & 35.74 & 3610.39 \\
                              & C\&CG+CG & 1.35E+05 & 1.35E+05 & \textbf{2} & \textbf{0} & \textbf{94.58} \\
      
        \multirow{2}{*}{D\_1} & BD+CG & 4.12E+06 & 1.18E+06 & 16 & 71.32 & 3835.40 \\
                              & C\&CG+CG & 2.85E+06 & 2.80E+06 & \textbf{5} & \textbf{1.69} & \textbf{222.78} \\
       
        \multirow{2}{*}{D\_2} & BD+CG & 4.39E+06 & 8.37E+05 & 9 & 80.92 & 3885.13 \\
                              & C\&CG+CG & 2.66E+06 & 2.63E+06 & \textbf{4} & \textbf{1.24} & \textbf{267.12} \\
     
        \multirow{2}{*}{E\_1} & BD+CG & 1.16E+07 & 9.61E+05 & 8 & 91.72 & 3791.17 \\
                              & C\&CG+CG & 4.63E+06 & 4.53E+06 & \textbf{4} & \textbf{2.08} & \textbf{466.73} \\
    
        \multirow{2}{*}{E\_2} & BD+CG & 1.06E+07 & 9.68E+05 & 6 & 90.9 & 4198.42 \\
                              & C\&CG+CG & 4.41E+06 & 4.32E+06 & \textbf{3} & \textbf{2.1} & \textbf{606.33} \\
      
        \multirow{2}{*}{F\_1} & BD+CG & 1.38E+07 & 1.03E+06 & 4 & 92.51 & 4068.84 \\
                              & C\&CG+CG & 4.69E+06 & 4.59E+06 & \textbf{5} & \textbf{2.11} & \textbf{898.20} \\
  
        \multirow{2}{*}{F\_2} & BD+CG & 1.63E+07 & 8.47E+05 & 3 & 94.81 & 4629.98 \\
                              & C\&CG+CG & 4.78E+06 & 4.65E+06 & \textbf{4} & \textbf{2.71} & \textbf{1639.75} \\
        \Xhline{1pt}
    \end{tabular}
\begin{tablenotes}
    \footnotesize
    \item \textit{Notes: BD+CG refers to the algorithm combining Benders Decomposition and CG; C\&CG+CG refers to our proposed algorithm combining C\&CG and CG.}
\end{tablenotes}
\end{threeparttable}
\end{table}

Table \ref{tab:BD} presents the performance comparison between the BD+CG method and the proposed C\&CG+CG solution method across 12 instances, reporting the upper bound, lower bound, number of iterations, optimality gap, and computational time. It can be observed that C\&CG+CG consistently achieves an optimality gap of 0\% or near 0\%, ensuring higher quality of solutions compared to BD+CG, which exhibits significantly larger gaps for most instances. For example, in Instance D\_1, BD+CG results in an optimality gap of 71.32\%, whereas C\&CG+CG reduces it to 1.69\%. A second observation is C\&CG+CG outperforms BD+CG in terms of computational time and the number of iterations across all instances. Notably, in Instance E\_2, C\&CG+CG solves the problem in 606.33 seconds, while BD+CG requires 4,198.42 seconds. This demonstrates the benefits of the proposed solution method, particularly for larger-scale instances. From these findings, we can conclude that the proposed C\&CG+CG method significantly outperforms BD+CG in terms of solution quality, computational efficiency, and iterative convergence.  

\subsubsection{Algorithm performance for integer solutions.}

We compare the performance of the proposed path-based model and its two solution procedures with that of the arc-based model solved using the C\&CG algorithm. The \textit{Always-Integer Procedure} is denoted as AIP, and the \textit{First-Continuous Procedure} is denoted as FCP. The \textit{Arc-based Robust Optimization} model (denoted as ARO), formulated in Appendix \ref{sec:arc-based}, is solved using the C\&CG algorithm. Notably, it cannot be solved by the algorithm combining C\&CG and CG, as CG requires paths. The optimality gap is computed as: $\text{Optimality gap} = (\text{UB}-\text{LB})/\text{UB} \times 100\%$, where $\text{UB}$ is the upper bound and $\text{LB}$ is the lower bound. The algorithm terminates when the optimality gap is no greater than 3\% or the computation time exceeds 7,200 seconds. 

\begin{insight}
Both AIP and FIP provide significant accelerations compared to ARO and achieve higher-quality solutions. Specifically, AIP and FIP reduce computational time by 73.83\% to 94.38\% and lower the optimality gaps to 0.20\% to 2.96\%, compared to the large gaps (10.50\% to 67.66\%) observed in ARO for Instances D and F. Moreover, FIP outperforms AIP on large-scale instances, offering better computational efficiency and slightly improved solution quality.
\end{insight}

\begin{table}[h]
    \centering
    \caption{Comparison of upper bound, lower bound, optimality gap, and computational time among ARO, AIP, and FIP.}
     \label{tab:results_comparison}
\begin{threeparttable}
    \begin{tabular}{llrrrrrrr}
        \Xhline{1pt}
        & & \multicolumn{3}{c}{Solution method (Instance\_1)} & \multicolumn{3}{c}{Solution method (Instance\_2)} \\
        \cmidrule(lr){3-5} \cmidrule(lr){6-8}
        Instance & Results & \textbf{ARO} & \textbf{AIP} & \textbf{FIP} & \textbf{ARO} & \textbf{AIP} & \textbf{FIP} \\
        \Xhline{0.6pt}
        \multirow{4}{*}{A} & Upper Bound       & 1.24E+05 & 1.24E+05 & 1.24E+05 & 1.24E+05 & 1.24E+05 & 1.24E+05 \\
                           & Lower Bound       & 1.24E+05 & 1.24E+05 & 1.24E+05 & 1.24E+05 & 1.24E+05 & 1.24E+05 \\
                           & Optimality gap (\%)          & 0.00     & 0.00     & 0.00     & 0.00     & 0.00     & 0.00     \\
                           & Computational time (s) & 4.40     & 7.26     & 7.74     & 5.88     & 30.83    & 10.83    \\
        \midrule
        \multirow{4}{*}{B} & Upper Bound       & 1.39E+05 & 1.39E+05 & 1.39E+05 & 1.39E+05 & 1.39E+05 & 1.39E+05 \\
                           & Lower Bound       & 1.39E+05 & 1.39E+05 & 1.39E+05 & 1.39E+05 & 1.39E+05 & 1.39E+05 \\
                           & Optimality gap (\%)          & 0.00     & 0.00     & 0.00     & 0.00     & 0.00     & 0.00     \\
                           & Computational time (s) & 11.51    & 33.12    & 13.53    & 19.67    & 52.24    & 36.51    \\
        \midrule
        \multirow{4}{*}{C} & Upper Bound       & 2.50E+05 & 2.50E+05 & 2.50E+05 & 2.44E+05 & 2.44E+05 & 2.44E+05 \\
                           & Lower Bound       & 2.50E+05 & 2.50E+05 & 2.50E+05 & 2.44E+05 & 2.44E+05 & 2.44E+05 \\
                           & Optimality gap (\%)          & 0.00     & 0.00     & 0.00     & 0.00     & 0.00     & 0.00     \\
                           & Computational time (s) & 1148.41  & \textbf{98.28}    & \textbf{145.23}   & 3616.21  & \textbf{126.93}   & \textbf{182.56}   \\
        \midrule
        \multirow{4}{*}{D} & Upper Bound       & 3.20E+06 & 2.96E+06 & 2.95E+06 & 3.01E+06 & 2.79E+06 & 2.77E+06 \\
                           & Lower Bound       & 2.68E+06 & 2.90E+06 & 2.94E+06 & 2.57E+06 & 2.74E+06 & 2.74E+06 \\
                           & Optimality gap (\%)          & 16.31    & \textbf{1.96}     & \textbf{0.20}     & 14.66    & \textbf{1.79}     & \textbf{0.98}     \\
                           & Computational time (s) & 7234.96  & \textbf{406.75}   & \textbf{314.11}   & 7233.39  & \textbf{523.08}   & \textbf{400.81}   \\
        \midrule
        \multirow{4}{*}{E} & Upper Bound       & 4.75E+06 & 4.78E+06 & 4.73E+06 & 4.52E+06 & 4.52E+06 & 4.60E+06 \\
                           & Lower Bound       & 4.57E+06 & 4.72E+06 & 4.67E+06 & 4.34E+06 & 4.43E+06 & 4.48E+06 \\
                           & Optimality gap (\%)          & 3.88     & \textbf{1.45}     & \textbf{1.10}     & 3.99     & \textbf{1.96}     & \textbf{2.77}     \\
                           & Computational time (s) & 7244.91  & \textbf{531.89}   & \textbf{436.52}   & 7246.55  & \textbf{814.92}   & \textbf{724.04}   \\
        \midrule
        \multirow{4}{*}{F} & Upper Bound       & 5.16E+06 & 4.80E+06 & 4.88E+06 & 1.13E+07 & 4.86E+06 & 4.88E+06 \\
                           & Lower Bound       & 4.62E+06 & 4.75E+06 & 4.79E+06 & 3.66E+06 & 4.72E+06 & 4.86E+06 \\
                           & Optimality gap (\%)          & 10.50    & \textbf{1.12}     & \textbf{1.84}     & 67.66    & \textbf{2.96}     & \textbf{0.35}     \\
                           & Computational time (s) & 7260.87  & \textbf{1202.20}  & \textbf{1136.48}  & 7263.62  & \textbf{1900.69}  & \textbf{1685.73}  \\
       \Xhline{1pt}
    \end{tabular}
    \begin{tablenotes}
    \footnotesize
    \item \textit{Notes: ARO refers to the use of the C\&CG algorithm to solve the arc-based robust optimization model; AIP refers to the Always-Integer Procedure for solving the path-based model; FIP refers to the First-Continuous Procedure for solving the path-based model.}
\end{tablenotes}
\end{threeparttable}
\end{table}

Table \ref{tab:results_comparison} reports the upper bound, lower bound, optimality gap, and computational time of ARO, AIP, and FIP for 12 different problem sizes. The main observation is that the proposed path-based model, along with the AIP and FIP, provides significant accelerations compared to solving the arc-based model using the C\&CG algorithm (ARO). Notably, ARO leaves very large optimality gaps ranging from 10.50\%  to 67.66\% for Instances D and F. In contrast, the AIP and FIP not only reduce computational time by  73.83\% to 94.38\% but also bring the optimality gaps down to 0.20\% to 2.96\%.

It can also be observed that when the problem size is small, all three algorithms are capable of finding the optimal solution (e.g., in Instances A and B) within a relatively short time. However, as the problem size increases, the computational time for ARO grows significantly, and the quality of the solution cannot be guaranteed. This is mainly due to the rapid expansion of the size of the space-time-SoC network as the problem scales, which makes the problem increasingly difficult to solve efficiently. In comparison, AIP and FIP show the ability to handle larger-scale problems with both high computational efficiency and solution quality due to the integration of the column generation algorithm.

Another observation is that FIP achieves solutions of comparable quality to AIP but with shorter computational times for Instances D, E, and F. For instance, FIP solves Instance D\_1 in 320 seconds with an optimality gap of 0.2\%, whereas AIP takes over 400 seconds with a gap of 1.96\%. The difference is even more pronounced in Instance E\_2, where FIP requires 1,686 seconds to achieve a gap of 0.35\%, compared to that AIP uses 1,900 seconds and obtain a solution with a gap of 2.96\%. Thus, these observations highlight the computational benefits of our path-based optimization model and algorithm over the C\&CG method for solving the arc-based optimization model. Furthermore, based on these findings, we conclude that the FIP method is preferable for solving large-scale problems.

\subsection{Benefits of utilizing E-MAVs with flexible compositions }
\label{sec:benefitFlexible}

In this subsection, we evaluate the practical benefits of the proposed optimization methodology compared to benchmarks with traditional fixed-composition electric buses that are widely implemented in practice. To this end, we first solve instance G\_1 using flexible-composition E-MAVs, with minimum and maximum compositions set to one and five, respectively, as introduced in Section \ref{sec:instances}. Subsequently, we compare these results with two cases involving different types of traditional fixed-composition vehicles, which lack the flexibility to decouple or couple units. The first type consists of two-unit vehicles, while the second type consists of four-unit vehicles. These three experiments are solved by the proposed FIP solution method.

\begin{insight}
Under comparable fleet investment, the proposed E-MAVs with flexible compositions significantly outperform fixed-composition buses. Compared to using flexible-composition E-MAVs, fixed-composition buses with two and four units increase passengers' costs by 16.22\% and 28.22\%, respectively, and average passenger waiting time by 17.46\% and 29.14\%. This highlights the substantial benefits of flexibility in improving service quality.
\end{insight}

\begin{table}[h]
\centering
\caption{Comparison of utilizing vehicles with flexible and fixed compositions.}
\label{table:performance}
\begin{tabular}{ccccccc}
\Xhline{1pt}
 \multirow{2}{*}{Compositions}& \multirow{2}{*}{\begin{tabular}[c]{@{}c@{}}Fleet size\\ of E-MAUs\end{tabular}} & \multicolumn{2}{c}{Passengers' costs} & \multicolumn{2}{c}{Average waiting time} \\ \cline{3-6}
          &            & Value (min)    & Diff. (\%)              & Value (min)               & Diff. (\%)           \\ \Xhline{0.6pt}
Flexible         & 94         & 105,855.00    & -                       & 5.72                & -                    \\
Fixed, two-unit vehicles         & 94         & 126,343.00    & 16.22                   & 6.93                & 17.46                \\ 
Fixed, four-unit vehicles        & 96         & 176,014.00    & 28.22                   & 9.78                & 29.14                \\ 
\Xhline{1pt}
\end{tabular}%
\end{table}

Table \ref{table:performance} presents the results, reporting the optimized fleet size of E-MAUs, passengers' costs, and the average passenger waiting time. First, we solve the case using flexible-composition E-MAVs, with optimization results indicating that 94 units are required, resulting in an average passenger waiting time of 5.72 minutes. For a fair comparison, when solving the cases with two-unit and four-unit vehicles, additional constraints on fleet size are applied after obtaining the fractional solution. These constraints limit the fleet size to 94 units (equivalent to 47 buses, each comprising two units) and 96 units (equivalent to 24 buses, each comprising four units) in the integer solutions, ensuring a fair comparison of service quality under nearly the same fleet investment. From the results in Table \ref{table:performance}, we can derive the following two observations: (1) Compared to using E-MAVs with flexible compositions, using fixed-composition vehicles made up of two units results in an increase of 16.22\% in passengers' costs and 17.46\% in average waiting time. (2) Using fixed-composition vehicles made up of four units further worsens service quality, leading to an increase of 28.22\% in passengers' costs and 29.14\% in average waiting time. These findings demonstrate that, under nearly the same level of fleet investment, the flexibility of E-MAVs in adjusting compositions significantly enhances service quality and reduces passengers' costs, highlighting their practical advantages over fixed-composition vehicles.

\subsection{Managerial insights from instances with entire-day operations }
\label{sec:entireDay}

We conclude by evaluating the practical benefits of the proposed methodologies in managing operations throughout an entire operational period. To this end, we use the FIP solution method to solve Instance H\_1, where operations begin at 7:00 and end at 22:00. The instance comprises 
$2.2\times 10^{6}$ nodes and  $3.2\times 10^{6}$ arcs, a scale comparable to that of the traditional E-VSP problem analyzed in \cite{Rolf2023}.

\begin{insight}
The results demonstrate the capability of the proposed methodology to manage large-scale operations with complex scheduling requirements while maintaining computational efficiency and solution quality. Specifically, the optimized circulation schedule highlights the flexibility of E-MAUs to adapt dynamically by performing regular bus and short-turning services along with intermittent recharging and composition adjustments. These operational patterns, supported by robust investment plans for charging infrastructure, enable effective utilization of resources, improve service quality, and ensure the practicality of deploying E-MAUs in real-world transit systems.
\end{insight}

\begin{figure}[h]
    \centering
    \includegraphics[width=0.9\linewidth]{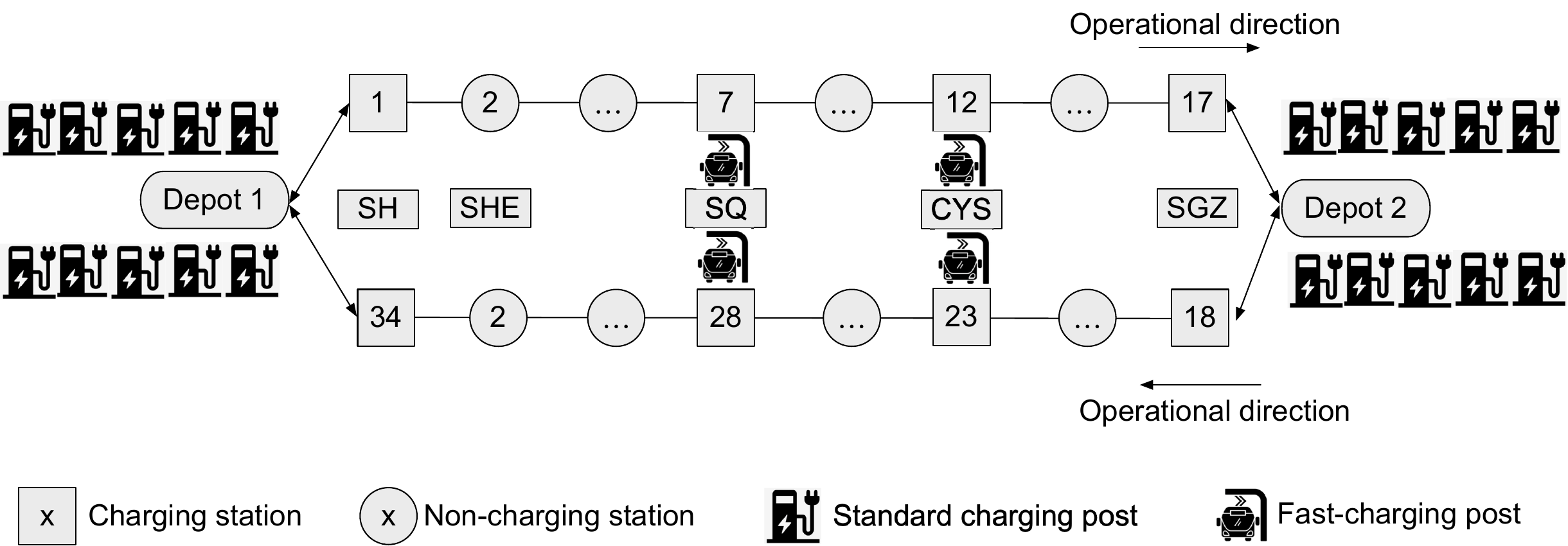}
    \caption{Optimized charging station location and capacity.}
    \label{fig:optimizedCharging}
\end{figure}

After 11 hours and 52 minutes of computation, a solution is obtained with an optimality gap of 0.96\%, featuring an optimized fleet size of 101 E-MAUs. The optimized strategic-level plans, specifically the robust investment plans for charging station locations and capacities, are presented in Figure \ref{fig:optimizedCharging}. The results indicate that both potential charging station locations are utilized, with each station equipped with two fast-charging posts, reaching the maximum allowed capacity.

\begin{figure}[h]
    \centering
    \includegraphics[width=0.9\linewidth]{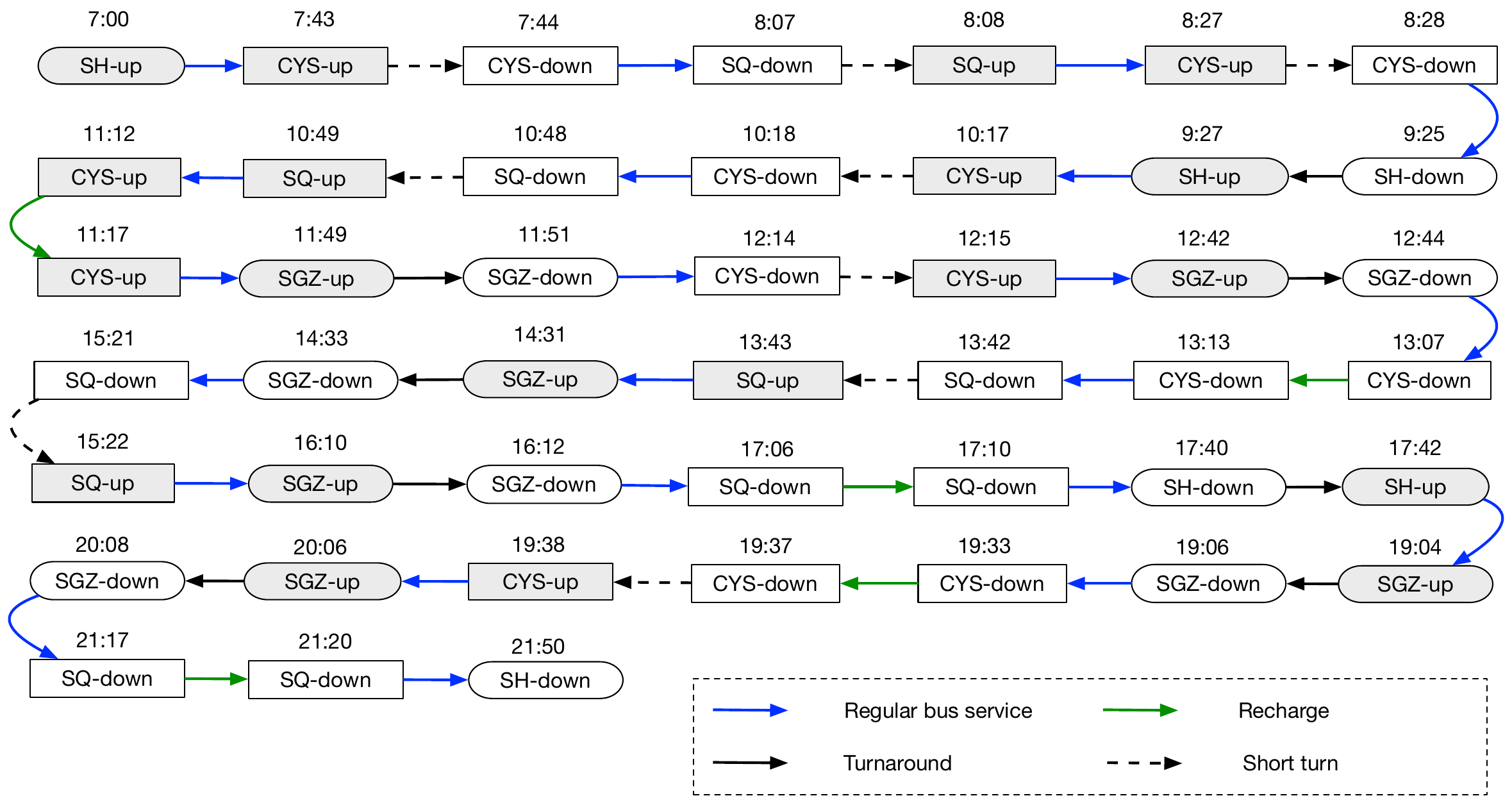}
    \caption{Optimized circulation schedule of an E-MAU between depots and charging stations from 7:00 to 22:00 with forward-only movements at non-charging stations.}
    \label{fig:optimizedCirculationAfter}
\end{figure}

Figures \ref{fig:optimizedCirculationAfter} shows the optimized circulation schedule of an E-MAU from 7:00 to 22:00, reporting the schedules between depots and charging stations. This schedule indicates that the unit undergoes multiple recharging, decoupling, and coupling operations at charging stations. Additionally, it performs a combination of regular bus and short-turning services on an intermittent basis. For instance, the unit departs from charging station SH in the up-operational direction to perform a regular bus service, arrives at station CYS, performs a short-turning service here, and returns to station SQ in the down-operational direction.

\section{Conclusion}
\label{sec:conclusion}

This paper addressed a robust charging station location and scheduling problem for electric modular autonomous units, extending traditional electric bus scheduling and vehicle routing problems by incorporating flexible circulations and integrated charging station location decisions. The problem jointly optimizes robust charging station locations and capacities, fleet sizing, as well as units' routing, scheduling, and charging decisions. It provides flexibility in deciding where, when, and for how long to charge, as well as when and where to decouple or couple units to perform specific service types. Additionally, it allows the units within an E-MAV to operate with varying SoC levels. We formulated the problem as a unit-based and path-based mixed-integer optimization model, employing a polyhedral uncertainty set to take the uncertain demand into account. The model is based on a space-time-SoC network, with both time and SoC levels discretized into finite sets. 

To address this challenging problem, we developed a double-decomposition solution methodology that integrates C\&CG with CG algorithms. In the C\&CG framework, the problem is decomposed into a master problem, which determines robust charging station locations, capacities, fleet sizing, and the routing and scheduling of units, and a subproblem, which optimizes passenger assignments. To handle the large number of paths in the master problem, a CG algorithm is employed, with a tailored label-setting algorithm designed to efficiently solve the pricing problem. The pricing problem is further accelerated by introducing three downsizing methodologies for the space-time-SoC network. These methods significantly reduce the number of nodes and arcs while preserving the optimality of the solutions. The proposed downsizing methodologies exhibit broad applicability and can be effectively employed across various models built on either space-time or space-time-SoC networks, highlighting their versatility and general effectiveness.  Additionally, an outer approximation algorithm is designed to generate high-quality initial solutions for GUROBI, accelerating the solution process of the subproblem in the C\&CG framework. Finally, we propose two procedures to effectively combine all algorithmic components.

Extensive computational experiments based on real-life instances yield four main conclusions. First, the proposed three methodologies for downsizing the space-time-SoC network achieve a considerable reduction in the network size, significantly improving the efficiency of solving the pricing problem and enabling large-scale problems with entire-day operations to be solved within manageable computation times. Second, our solution methodology outperforms state-of-the-art benchmarks by delivering higher-quality solutions in shorter computational times. Third, this study introduces a novel modeling methodology and algorithmic framework for the robust charging station location and routing-scheduling problem with electric modular autonomous units. Our proposed framework is generalizable, providing valuable insights and solution techniques applicable to a wide range of routing and scheduling optimization problems. Fourth, the proposed methodology offers substantial practical benefits, including significant improvements in service quality and the flexibility to schedule units to alternate between two service patterns multiple times during operations. Our methodology and managerial insights offer robust support for the adoption and integration of electrification solutions in logistics and transportation.

For future research, incorporating machine learning into the proposed algorithm is a promising direction. Additionally, exploring distributionally robust optimization techniques to address uncertainty in passenger demand could further enhance the practicality of the solutions.

\ACKNOWLEDGMENT{%
This work was supported by the National Natural Science Foundation of China (No. 72288101). The authors thank Alexandre Jacquillat of the Massachusetts Institute of Technology for his insightful feedback. We also thank Tommaso Gecchelin and NExT Future Transportation Inc. for providing data related to the NExT modular autonomous vehicles used in the numerical experiments of this research.
}

\newpage

\newpage

\begin{APPENDICES}
\section{Examples}
\label{sec:example}

In this section, we present examples of the reallocation options of E-MAUs, the timetable and vehicle schedule based on the space-time-SoC network, and constructing super travel arcs.

\begin{example}\label{ex:time_space_soc_network}
Figure~\ref{time_space_soc_network} illustrates the timetable and vehicle schedule for an E-MAV composed of two E-MAUs within the time-space-SoC network. At the beginning, the E-MAV, comprising two units with 80\% SoC, departs from depot 1. When it arrives at the second station, this E-MAV is decoupled into two individual units. Both units are recharged to 80\% SoC. Subsequently, as shown in Figure~\ref{time_space_soc_network}(a), one E-MAU continues to execute the regular bus service, eventually arriving at depot 2. Meanwhile, the other unit performs a short turn at the second station, returning to the first station before finally arriving back at the first depot, as depicted in Figure~\ref{time_space_soc_network}(b).
\end{example}

\begin{figure}[htbp]
     \centering
     \begin{subfigure}{\textwidth}
         \centering
\includegraphics[width=0.7\textwidth]{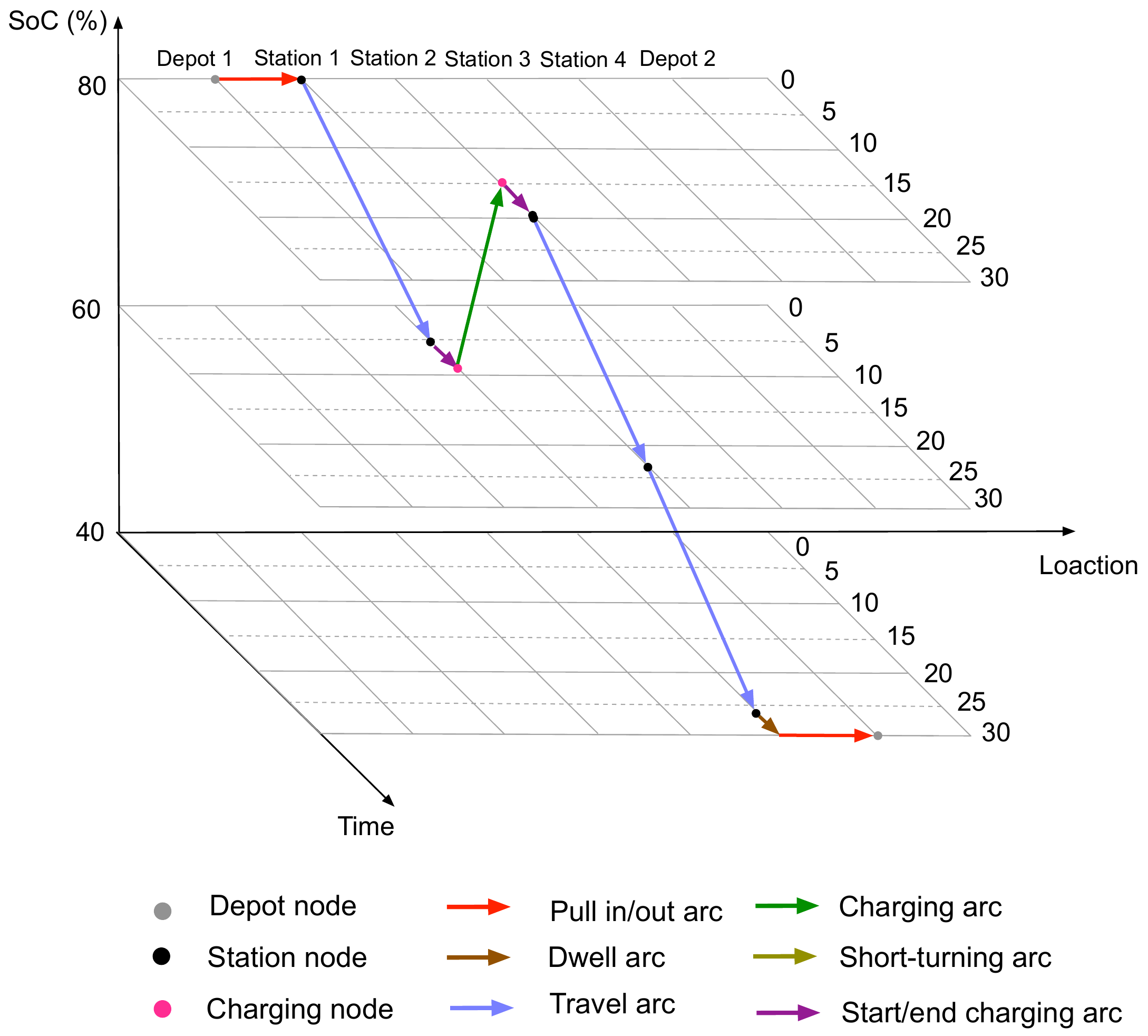}
         \caption{The E-MAU that executes a regular bus service}
     \end{subfigure}
    
\vspace{2em}

     \begin{subfigure}{\textwidth}
         \centering
\includegraphics[width=0.8\textwidth]{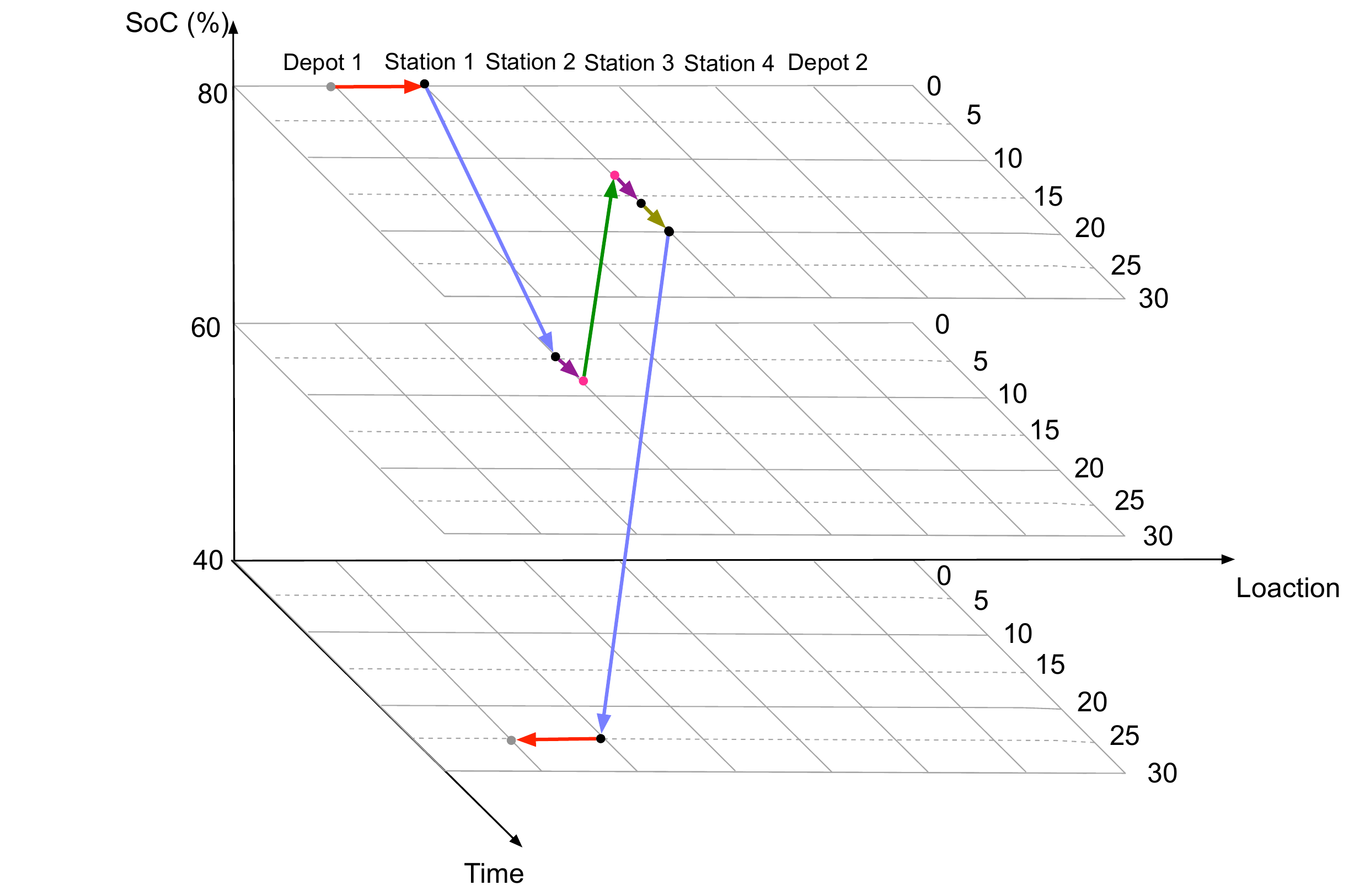}
         \caption{The E-MAU that executes a short-turning service}
     \end{subfigure}
        \caption{Illustration of the timetable and vehicle schedule for an E-MAV comprised two E-MAUs based on the space-time-SoC network.}
\label{time_space_soc_network}
\end{figure}

\begin{example}\label{example:superarcs}

Figure~\ref{super_arc} illustrates the construction of super travel arcs. 
   \begin{figure}[htbp]
     \centering     \includegraphics[width=0.7\textwidth]{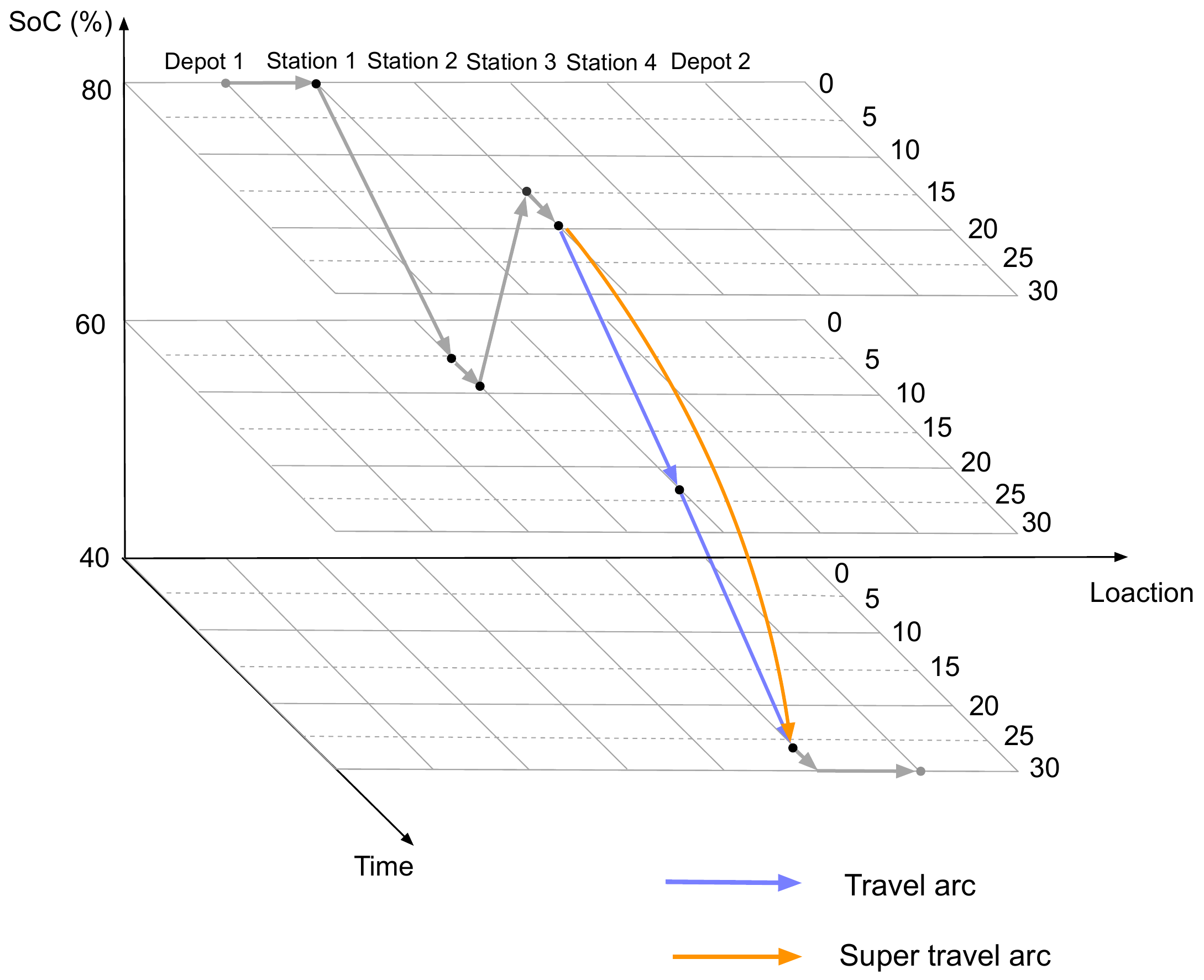}
        \caption{Illustration of constructing super travel arcs by aggregating normal travel arcs.}
\label{super_arc}
\end{figure}
 
\end{example}

\newpage

\section{Notations used in the model}
\label{sec:notations}

This section introduces the notations used in the model in Table \ref{table-para}.

\begin{table}[h]
  \centering
  \caption{Sets and parameters used in the model.} \label{table-para}
  \begin{tabular}{ll}
\Xhline{1pt} \noalign{\smallskip}
\multicolumn{2}{l}{\bf{Sets}}\\
$\mathcal{S}$& Set of stations, $\mathcal{S}=\{1, 2, \ldots, \left| \mathcal{S} \right| \},$ indexed by $s$\\
$\overline{\mathcal{S}}$& Set of upstream stations, $\overline{\mathcal{S}}=\{1, 2, \ldots, \left| \overline{\mathcal{S}} \right| \}$\\
$\underline{\mathcal{S}}$& Set of downstream stations, $\underline{\mathcal{S}}=\{\left| \overline{\mathcal{S}} \right| +1, \ldots, \left| \overline{\mathcal{S}} \right|+\left| \underline{\mathcal{S}} \right| \}$\\
$\mathcal{Q}$& Set of sections, $\mathcal{Q}=\{1, 2, \ldots, \left| \mathcal{Q} \right| \},$ indexed by $q$ or $(s, s+1)$\\
$\mathcal{D}$& Set of depots, indexed by $d$\\
$\mathcal{U}_{fast}$& Set of fast-charging stations\\
$\mathcal{U}$& Set of charging stations, indexed by $u$, , $\mathcal{U} = \mathcal{D} \cup \mathcal{U}_{fast}$\\
$\mathcal{T}$& Set of discretized time intervals, $\mathcal{T}=\{1, 2, \ldots, \left| \mathcal{T} \right|\}$, indexed by $t$\\
$\mathcal{N}$& Set of nodes, indexed by $n, v$\\
$\mathcal{N}_{st}$& Set of nodes located in station $s\in\mathcal{S}$ at time $t\in\mathcal{T}$\\
$\mathcal{A}$& Set of arcs, indexed by $a$ or $(i, j)$\\
$\mathcal{A}_{q,t}$& Set of arcs related to segment $q$ and time $t$\\
$\mathcal{E}$& Set of discretized SoC, $\mathcal{E}=\{e_{min}, e_{min}+\sigma, \ldots, e_{max}\}$, where $\sigma$ represents the step \\
& size per SoC level, indexed by $e$\\
$\mathcal{P}$& Set of paths, $\mathcal{P}=\{1, 2, \ldots, \left| \mathcal{P} \right| \},$ indexed by $p$\\
\multicolumn{2}{l}{\bf{Parameters}}\\
$\Delta$& Duration of a discretized time interval (unit: min)\\
$\chi_{n,v,t}$& Travel time from node $n$ to node $v$ at time $t$\\
$\sigma$& Step size between two discretized SoCs (unit: \%)\\
$\theta_1$, $\theta_2$, $\theta_3$ & Weighting monetary coefficients for operator's costs, passengers' waiting costs, \\ 
& and penalties associated with unserved passengers\\
$c_p$ & Operational cost of path $p$\\
$l^a_p$ & Binary indicator. $l^a_p = 1$ if path $p$ passes through arc $a$; otherwise, $l^a_p = 0$\\
$g_{u,t,p}$ & Binary indicator. If path $p\in\mathcal{P}$ is charged at charging station $u\in\mathcal{U}$ at time $t$,   \\& $g_{u,t,p}=1$;
otherwise $g_{u,t,p}=0$\\
$C$& Capacity of an E-MAU\\
$m_u$& Costs of installing a fast-charging post at a fast-charging station $u\in\mathcal{U}$\\
$N^{charge}_u$& Maximum number of fast-charging posts that can be installed at a fast-charging station $u$\\
$N^{charge}_d$& Maximum number of standard charging posts that can be used at depot $d$\\
$N^{com}$& Maximum compositions of each E-MAV\\
$\theta^{soc}_{(n,v)}$ & The SoC required for traveling from nodes $n$ to $v$, $(n,v)\in\mathcal{A}$.  $\theta^{soc}_{(n,v)}=k\cdot\sigma,k\in\mathbb{Z}$\\
$f_d$ & The nearest station to depot $d$\\
$\dot{f}_s$ & The station located opposite to station $s$\\
$E_{e,u}$ & Charging quantity of the E-MAU with an initial SoC of $e$ at charging station $u$\\
$\hat{\alpha}_q^t$ & The number of newly arrival passenger demand on section $q$ at time $t$ under the \\ & deterministic condition\\
$\overline{\alpha}_q^t$ & The nominal number of newly arrival passenger demand on section $q$ at time $t$ under the \\ & uncertain condition\\
$\Tilde{\alpha}_q^t$ & The maximum deviation between the uncertain and the nominal demand on section $q$ at \\ & time $t$ under the uncertain condition\\
$\beta_{d,p}^{-}$ & Binary indicator. $\beta_{d,p}^{-} = 1$ if path $p$ begins at depot $d$; otherwise, $\beta_{d,p}^{-} = 0$\\
$\beta_{d,p}^{+}$ & Binary indicator. $\beta_{d,p}^{+}=1$ if path $p$ ends at depot $d$; otherwise, $\beta_{d,p}^{+}=0$\\
\Xhline{1pt}
\end{tabular}
\end{table}

\newpage
\section{Proofs of Lemma 1 and Propositions 1-4}
\label{sec:proof}

\textbf{Proof of Lemma \ref{lemma:bound}.} 
From constraints (\ref{eq:dual1}), we have:
\begin{align*}
\delta_q^{t} = \delta_q^{\left|\mathcal{T}\right|} + (\delta_q^{\left|\mathcal{T}\right| - 1} - \delta_q^{\left|\mathcal{T}\right|}) + \dots + (\delta_q^{t+1} - \delta_q^{t+2}) + (\delta_q^{t} - \delta_q^{t+1}).
\end{align*}

Using the structure of the constraints and the non-negativity conditions, each difference term \((\delta_q^{t+i} - \delta_q^{t+i+1})\) is bounded by \(\theta_2\), and the last term \(\delta_q^{\left|\mathcal{T}\right|}\) is bounded by \(\theta_2 + \theta_3\). Summing these bounds yields:
\begin{align*}
\delta_q^{t} \leq (\left|\mathcal{T}\right| - t + 1)\theta_2 + \theta_3.
\end{align*}

Thus, \( 0 \leq \delta_q^t \leq M_t = (\left|\mathcal{T}\right| - t + 1)\theta_2 + \theta_3, \, \forall t \in \mathcal{T}, \, q \in \mathcal{Q}\) holds. Furthermore, since $\delta_q^t$ is bounded, $\Psi$ is a bounded polyhedron. \Halmos

\textbf{Proof of Proposition \ref{pp1}.} 
Since both the uncertainty set $\Xi$ and $\Psi$ are bounded polyhedra, an optimal solution $(\bm{\zeta}^{*},\bm{\delta}^{*})$ of $\text{DF}(\mathbf{x})$ exists such that $\bm{\delta}^{*}$ is the extreme point of $\Psi$ and $\bm{\zeta}$ is the extreme point of $\Xi$ \citep{horst2013global}. Consequently, when $\Pi_q \in \mathbb{Z}^+, \forall q\in\mathcal{Q}$ and $\Lambda_t \in \mathbb{Z}^+, \forall t\in\mathcal{T}$, it holds that $\bm{\zeta}^{*}\in\{0, 1\}^{\left|\mathcal{T}\right|\times\left|\mathcal{Q}\right|}$. \Halmos \endproof

\textbf{Proof of Proposition \ref{max_min_soc}.} 
Consider a charging-station subpath \(\Tilde{p} = \{n_{\Tilde{p}}^1, n_{\Tilde{p}}^2, \dots, n_{\Tilde{p}}^m\} \in \textit{CS}_s\). For any arc \((n_{\Tilde{p}}^j, n_{\Tilde{p}}^{j+1}) \in \mathcal{A}\) along this subpath, by the definition of the SoC transition, the SoC at node \(n_{\Tilde{p}}^{j+1}\) satisfies:
$$
e(n_{\Tilde{p}}^{j+1}) = e(n_{\Tilde{p}}^j) - \theta^{\text{soc}}_{(n_{\Tilde{p}}^j, n_{\Tilde{p}}^{j+1})}.
$$
Extending this relationship recursively over the entire subpath, the SoC at the final node \(n_{\Tilde{p}}^m\) can be expressed as:
$$
e(n_{\Tilde{p}}^m) = e(n_{\Tilde{p}}^1) - \sum_{j=1}^{m-1} \theta^{\text{soc}}_{(n_{\Tilde{p}}^j, n_{\Tilde{p}}^{j+1})}.
$$
Since the maximum SoC at \(n_{\Tilde{p}}^1\) is given by \(e_{\text{max}}\), the SoC at station node \(n\) located at station $s$ along the subpath \(\Tilde{p}\) satisfies:
$$
e(n) \leq \max_{\Tilde{p} \in \textit{SC}_s} \left\{ e_{\text{max}} - \sum_{j=1}^{m-1} \theta^{\text{soc}}_{(n_{\Tilde{p}}^j, n_{\Tilde{p}}^{j+1})} \right\}.
$$

Similarly, consider a station-charging subpath \(\Tilde{p} = \{n_{\Tilde{p}}^1, n_{\Tilde{p}}^2, \dots, n_{\Tilde{p}}^m\} \in \textit{SC}_s\). 
The SoC at the first node \(n_{\Tilde{p}}^1\) along this subpath can then be written as:
$$
e(n_{\Tilde{p}}^1) = e(n_{\Tilde{p}}^m) + \sum_{j=1}^{m-1} \theta^{\text{soc}}_{(n_{\Tilde{p}}^j, n_{\Tilde{p}}^{j+1})}.
$$
Since the minimum SoC at \(n_{\Tilde{p}}^m\) is \(e_{\text{min}}\), we have
$$
e(n) \geq \min_{\Tilde{p} \in \textit{CS}_s} \left\{ e_{\text{min}} + \sum_{j=1}^{m-1} \theta^{\text{soc}}_{(n_{\Tilde{p}}^j, n_{\Tilde{p}}^{j+1})} \right\}.
$$

This completes the proof. \Halmos \endproof

\textbf{Proof of Proposition \ref{time}.} 
In this proposition, we ignore the SoC dimension and project all points located at the same place and at the same time onto the two-dimensional plane. For each arc \((v, n) \in \mathcal{A}\), the time at node \(n\) satisfies:
$$
t(n) = t(v) + \chi_{(v, n)},
$$
where \(\chi_{(v, n)}\) represents the travel time between nodes \(v\) and \(n\).

Consider a charging-station subpath \(\Tilde{p} = \{n_{\Tilde{p}}^1, n_{\Tilde{p}}^2, \dots, n_{\Tilde{p}}^m\} \in \textit{SC}_s\). The time at the final node \(n_{\Tilde{p}}^m\) along this subpath can be expressed recursively as:
$$
t(n_{\Tilde{p}}^m) = t(n_{\Tilde{p}}^1) + \sum_{j=1}^{m-1} \chi_{(n_{\Tilde{p}}^j, n_{\Tilde{p}}^{j+1})}.
$$
Since the maximum time at the depot node \(n_{\Tilde{p}}^1\) is \(\left|\mathcal{T}\right|\), we have:
$$
t(n) \leq \max_{\Tilde{p} \in \textit{SC}_s} \left\{ \left|\mathcal{T}\right| - \sum_{j=1}^{m-1} \chi_{(n_{\Tilde{p}}^j, n_{\Tilde{p}}^{j+1})} \right\}.
$$

Now we consider a station-charging subpath \(\Tilde{p} = \{n_{\Tilde{p}}^1, n_{\Tilde{p}}^2, \dots, n_{\Tilde{p}}^m\} \in \textit{CS}_s\). The time at the first node \(n_{\Tilde{p}}^1\) along this subpath can be expressed as:
$$
t(n_{\Tilde{p}}^1) = t(n_{\Tilde{p}}^m) - \sum_{j=1}^{m-1} \chi_{(n_{\Tilde{p}}^j, n_{\Tilde{p}}^{j+1})}.
$$
Since the minimum time at the depot node \(n_{\Tilde{p}}^m\) is \(\delta\), we have:
$$
t(n) \geq \min_{\Tilde{p} \in \textit{CS}_s} \left\{ \delta + \sum_{j=1}^{m-1} \chi_{(n_{\Tilde{p}}^j, n_{\Tilde{p}}^{j+1})} \right\}.
$$

This completes the proof. \Halmos \endproof

\textbf{Proof of Proposition \ref{pp4}.}
Let \( \mathcal{G} \) denote the original space-time-SoC network and \( \check{\mathcal{G}} \) the modified network where super travel arcs replace sequences of travel arcs. For each path $p\in\mathcal{P}$, there exists and only exists one corresponding path \( p'\in\mathcal{P}^{'} \) in \( \check{\mathcal{G}} \) such that the following condition holds: if all super travel arcs \( \check{a} \in \check{\mathcal{A}}_{travel} \) in \( p' \) are replaced by their constituent arcs in \( \mathcal{A}_{\check{a}} \), then \( p' \) and \( p \) are identical. Let $\check{\mathcal{Q}}$ denote the corresponding section on \( \check{\mathcal{G}} \). Similar to equation (\ref{reduced_cost_pp}), the reduced cost of path $p'$ is calculated as follows:
\begin{align}
  R_{p'} = \theta_{oper} c_{p'} &- \sum_{d\in\mathcal{D}}(\beta_{dp}^{-}-\beta_{dp}^{+})\iota_d- \sum_{a\in\mathcal{A}_{\check{a}}}l_p^{\check{a}}\pi_{\check{a}}- \sum_{u\in\mathcal{U}}\sum_{t\in\mathcal{T}}g_{u,t,p}\rho_{ut}  +C\sum_{w\in\mathcal{W}}\sum_{\check{q}\in\check{\mathcal{Q}}}\sum_{t\in\mathcal{T}}\varsigma_{\check{q}t}\sum_{\check{a}\in\mathcal{A}_{\check{q},t}}l^{\check{a}}_p. \tag{11}
\end{align}
where $\pi_{\check{a}}=\sum_{a\in\mathcal{A}_{\check{a}}}\pi_{a}$, and $\varsigma_{\check{q}t}=\sum_{\check{q}\in\check{\mathcal{Q}}}\varsigma_{qt}$. Then for all $p\in\mathcal{P}$, there exist $p'$ that satisfies $R_{p'}=R_{p}$ and vice versa. 

The pricing problem involves finding the path with the minimum reduced cost. As paths in \( \mathcal{G} \) and \( \check{\mathcal{G}} \) are equivalent with respect to reduced cost, the optimal objective values of the pricing problem in \( \mathcal{G} \) and \( \check{\mathcal{G}} \) are identical, i.e., $opt[\text{PP}(\mathcal{G})] = opt[\text{PP}(\check{\mathcal{G}})]$.
This completes the proof. \Halmos \endproof

\newpage
\section{Procedure of the tailored label-correcting algorithm}
\label{sec:labelcorrecting}

In this section, we introduce the procedure of the tailored label-correcting algorithm in Algorithm \ref{alg:spfa-specific}.

\begin{algorithm}[h]
\caption{Tailored label-correcting algorithm}\label{alg:spfa-specific}
\SetKwInOut{Input}{Input}
\SetKwInOut{Output}{Output}
\SetKw{KwGoTo}{go to}
\SetKw{KwSet}{set}
\SetKw{KwBreak}{break}
\SetKw{Return}{return}
\SetAlgoLined

\Input{The time-space-SoC network $\mathcal{G}=\{\mathcal{N}, \mathcal{A}\}$, cost functions and dual variables.}
\Output{Paths with minimum reduced costs for each depot pair.}
\BlankLine

Initialize $Q$ to empty queue;\\
Initialize $inQueue[n] = \text{False}$ for all $n \in \mathcal{N}$;\\
\ForEach{node $n \in \mathcal{N}$}{
    \lIf{$n \in \{n^{source}_d | d \in \mathcal{D}\}$}{
        $\Psi(n) = 0$
    }
    \lElse{
        $\Psi(n) = \infty$
    }
    \If{$n \in \{n^{source}_d | d \in \mathcal{D}\}$}{
        Enqueue $n$ into $Q$;\\
        $inQueue[n] = \text{True}$;
    }
}
\While{$Q$ is not empty}{
    $n \gets$ Dequeue from $Q$;\\
    $inQueue[n] = \text{False}$;\\
    \For{each outgoing arc $a=(n, v) \in \mathcal{A}$}{
        Calculate the tentative label for $v$ as follows:\\
        \Indp
        $\Psi_{v}'= \Psi_n + \theta_{oper} * c_{(n,v)} - \beta_{(n,v)}\iota_{d(n,v)}- \pi_{l(n)l(v)t(n)} - \rho_{l(v)t(v)} + C\sum_{w\in\mathcal{W}}\varsigma_{l(n)l(v)t(n)}(w)$;\\
        \Indm
        \If{$\Psi_{v}' < \Psi_v$}{
            $\Psi_v = \Psi_{v}'$;\\
            Update the predecessor of $v$ to $n$;\\
            \If{not $inQueue[v]$}{
                Enqueue $v$ into $Q$;\\
                $inQueue[v] = \text{True}$;
            }
        }
    }
}

\For{each depot $d \in \mathcal{D}$}{
    $p_d^* \gets$ reconstruct the path from $n^{source}_d$ to $n^{sink}_d$ using the predecessor links;\\
    \Return{$p_d^*$ and $\Psi_{p_d^*}$};
}

\end{algorithm}

\newpage
\section{Detailed flowchart of the proposed algorithm combining C\&CG and CG}
\label{sec:detailedFlowchart}

In this section, we introduce the proposed algorithm combining C\&CG and CG with the first-continuous procedure in detail in Figure \ref{fig:algorithm_overall_first}.

\begin{figure}[h]
    \centering
    \includegraphics[width=\linewidth]{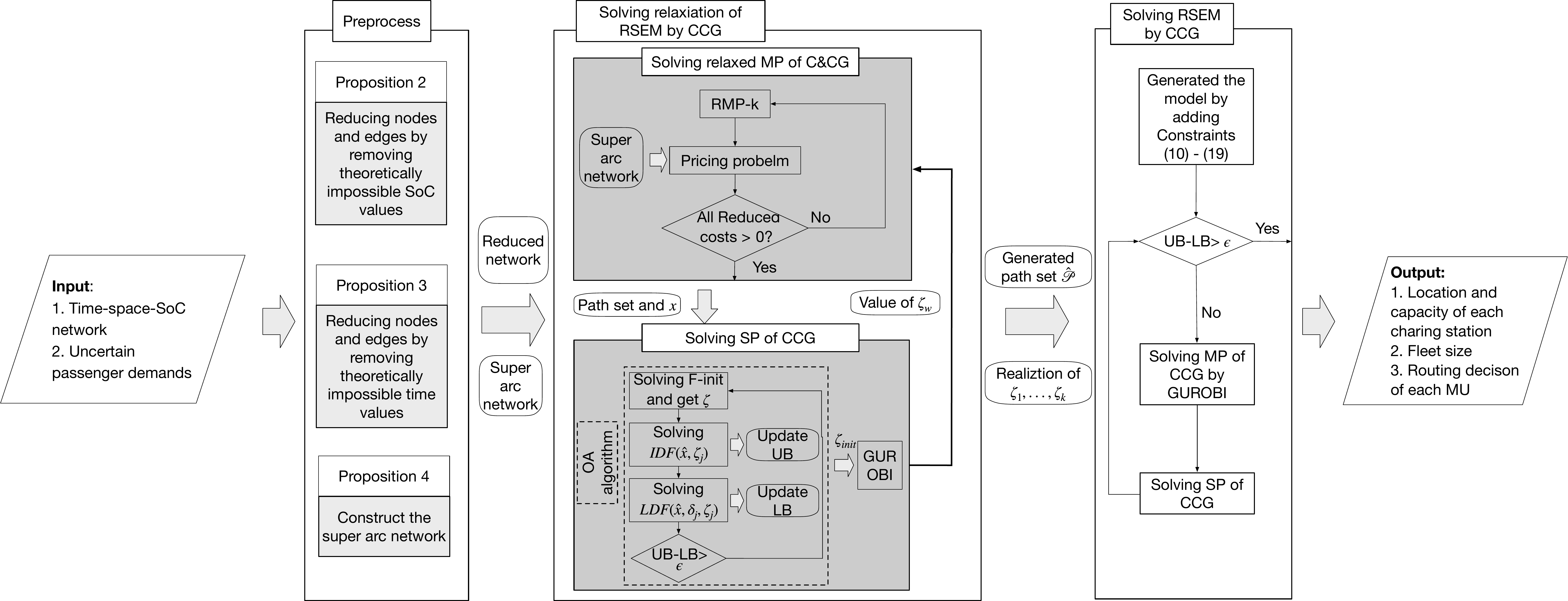}
    \caption{Detailed flowchart of the proposed algorithm with the first-continuous procedure.}
    \label{fig:algorithm_overall_first}
\end{figure}

\newpage
\section{Procedure of the outer approximation method}\label{procedure_outer}

This section presents the procedure of the outer approximation method.

\textbf{Outer approximation algorithm.}
\begin{itemize}
    \item \textbf{Step 1. Initialization.}  
    \begin{itemize}
        \item Input the solution $\hat{\mathbf{x}}$ passed from the MP.
        \item Set \( j \gets 1 \), \( UB \gets +\infty \), \( LB \gets -\infty \), and the tolerance level \( \epsilon^{OA} \).
        \item Define \(\textit{F-init} = \min\limits_{\bm{\zeta}} \{\bm{\Tilde{\alpha}}^\mathbf{T} \bm{\zeta} : \bm{\zeta} \in \Xi\}\). Solve \textit{F-init} to obtain the initial realization \(\bm{\zeta}_j\).
        \item Go to Step 2.
    \end{itemize}

    \item \textbf{Step 2. Solve the Dual Problem (\(\text{IDF}\)).}  
    \begin{itemize}
        \item Solve \(\text{IDF}(\hat{\mathbf{x}}, \bm{\zeta}_j)\) \eqref{bilinear_obj} - \eqref{eq:dual4}. Let $\bm{\delta}_j$ denote its optimal solution.
        \item Update \( LB \gets \text{IDF}(\hat{\mathbf{x}}, \bm{\zeta}_j) \).
        \item Go to Step 3.
    \end{itemize}

    \item \textbf{Step 3. Linearize and solve the SP.}  
    \begin{itemize}
        \item Linearize the bilinear term \(\bm{\delta}^\mathbf{T} \bm{\zeta}\) at \((\bm{\zeta}_j, \bm{\delta}_j)\):  
        \[
        L_j(\bm{\zeta}, \bm{\delta}) = (\bm{\delta} - \bm{\delta}_j)^\mathbf{T} \bm{\zeta}_j + (\bm{\zeta} - \bm{\zeta}_j)^\mathbf{T} \bm{\delta}_j + \bm{\delta}^\mathbf{T} \bm{\zeta}.
        \]  
        \item  Define the linearized SP (denoted as \(\text{LDF}(\hat{\mathbf{x}}, \bm{\delta}_j, \bm{\zeta}_j)\)) as:  
        \begin{align*}
    \text{LDF}(\hat{\mathbf{x}}, \bm{\delta}_j,\bm{\zeta}_j)=\max_{\bm{\zeta},\bm{\delta}}
    &\ \sum_{q\in\mathcal{Q}}\sum_{t\in\mathcal{T}}\big[\overline{\alpha}_{q}^t -  C\sum_{a\in\mathcal{A}_{q,t}}\sum_{p\in\mathcal{P}}l_{p}^a x_p )\delta_{q}^t + \Delta\big]\\
      \mbox{s.t.} \ 
    &\Delta\leq L_i(\bm{\zeta}, \bm{\delta}) \quad  i=1,2,...,j,\\
    &(\ref{df1}) - (\ref{df2}).
\end{align*}
        \item Solve \(\text{LDF}(\hat{\mathbf{x}}, \bm{\delta}_j, \bm{\zeta}_j)\).
        \item Let \((\bm{\zeta}_{j+1}, \bm{\delta}_{j+1})\) denote its optimal solution.
        \item Update \( UB \gets \text{LDF}(\hat{\mathbf{x}}, \bm{\delta}_j, \bm{\zeta}_j) \).
        \item Go to Step 4.
    \end{itemize}

    \item \textbf{Step 4. Check termination conditions.}  
    \begin{itemize}
        \item If \( UB - LB < \epsilon^{OA} \), terminate and return \(\bm{\zeta}\) as the optimal solution.
        \item Otherwise, increment \( j \gets j + 1 \) and go to Step 2.
    \end{itemize}
\end{itemize}

\newpage
\section{Overall frameworks of the proposed algorithm combining C\&CG and CG}
\label{sec:overallFlowchart}

The overall frameworks of two procedures of our algorithm combining C\&CG and CG are presented in Figure~\ref{fig:algorithms}. 

\begin{figure}[htbp]
     \centering
     \begin{subfigure}{0.49\textwidth}
         \centering
\includegraphics[width=\textwidth]{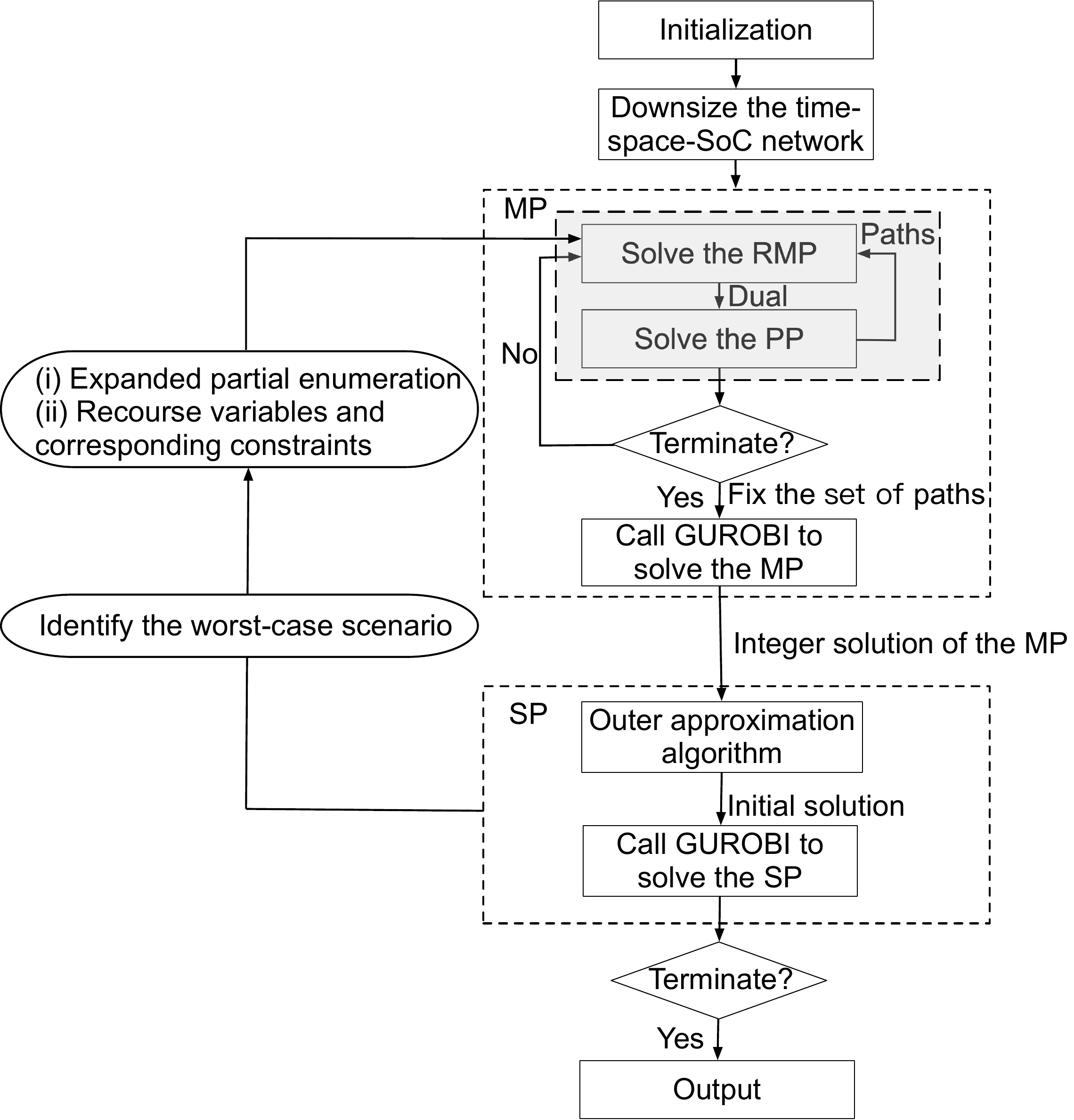}
         \caption{Always-integer procedure}
     \end{subfigure}
     \begin{subfigure}{0.49\textwidth}
         \centering
\includegraphics[width=\textwidth]{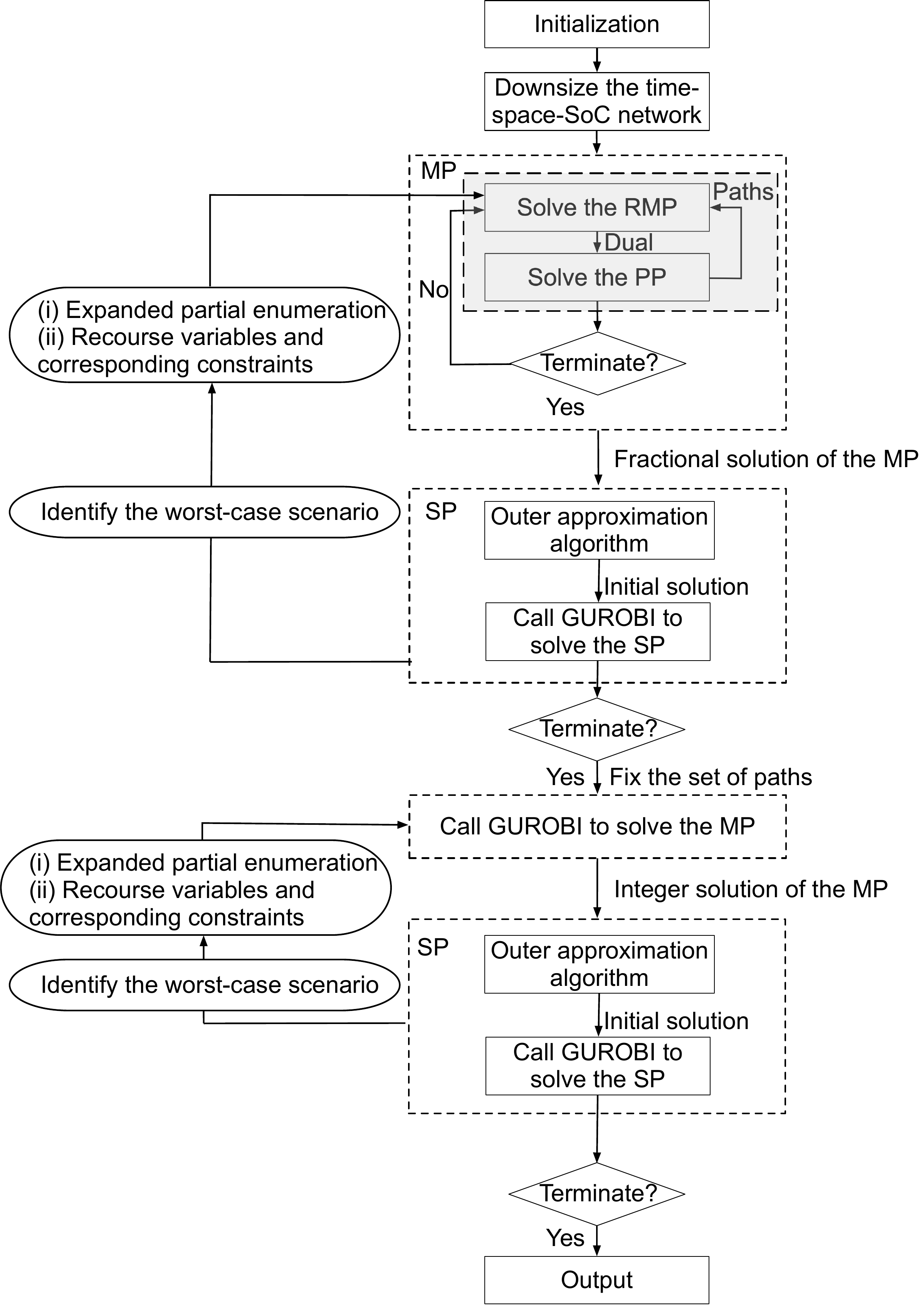}
         \caption{First-continuous procedure}
     \end{subfigure}
        \caption{Overview of two procedures of the proposed solution framework.}
\label{fig:algorithms}
\end{figure}

\newpage
\section{Arc-based robust optimization model}
\label{sec:arc-based}
To show the effectiveness of the proposed model and algorithm, we formulate the following arc-based model as the benchmark. To do so, we first introduce additional variables and sets used in this model:

\begin{itemize}
\item $z_{(n,v)}$: Number of MAUs running on arc $(n, v)$.
\item $\mathcal{A}^{+}_{n}$/$\mathcal{A}^{-}_{n}$: Set of outgoing / incoming arcs of node $n$.
\end{itemize}
Besides, we let $y^{0}_{s,s+1}=0$ for the sake of simplify modeling. We can now formulate the model as follows:
\begin{align}
\min_{\mathbf{z}, \mathbf{r}, \mathbf{y}} \quad & \theta_{1}(\sum_{(n,v)\in\mathcal{A}}c_{(n,v)}z_{(n,v)} + \sum_{u\in\mathcal{U}_{fast}}m_ur_u)  + \theta_{2}\sum_{t\in\mathcal{T}}\sum_{q\in\mathcal{Q}} y^t_{q} + \theta_{3}\sum_{q\in\mathcal{Q}} y^{\left|\mathcal{T}\right|}_{q}, \tag{12a}\\
    &\sum_{(n,v)\in\mathcal{A}^{+}_n}z_{(n,v)} - \sum_{(v,n)\in\mathcal{A}^{-}_n}z_{(v,n)} = 0
&\forall n\in\mathcal{N}\label{eq:vehicle_flow_appendix}\\
& \sum_{(n,v)\in\mathcal{A}_{q,t}}z_{(n,v)}\leq  N^{com}, &\forall q\in\mathcal{Q},t\in\mathcal{T},\\
&\sum_{(n,v)\in\mathcal{A}_{ut}} z_{(n,v)}\leq \begin{cases}
        r_u, & \text{if } u \in \mathcal{U}_{fast} \\
        N_d^{charge} , & \text{if } u \in \mathcal{U}  \setminus \mathcal{U}_{fast}
    \end{cases} & t\in\mathcal{T},\label{capacity_charge_appendix}\\
&y_{q}^t = y_{q}^{t-1} + \alpha_{q}^t - b_{q}^t, &\forall t\in\mathcal{T}, q\in\mathcal{Q} \\
&b_{q}^t \leq C  \sum_{(n,v)\in\mathcal{A}_{qt}}z_{(n,v)} , &\forall t\in\mathcal{T}, q=(s,s+1)\in\mathcal{Q} \\
&z_{(n,v)} \in \mathbb{Z}_{\geq 0}, &\forall (n, v)\in\mathcal{A}\label{z_domain_appendix}\\
& r_u \in \{0,1,...,N^{charge}_u\} , &\forall u\in\mathcal{U}_{fast}.\label{r_domain_appendix} \\
&b_{q}^t\geq 0,  &\forall t\in\mathcal{T}, q\in\mathcal{Q},\\
&y_{q}^t\geq 0, &\forall t\in\mathcal{T}, q\in\mathcal{Q} \label{y_domian_appendix}
\end{align}

The robust model can be formulated as follows
\begin{align*}
\min_{\mathbf{z}, \mathbf{r}, \mathbf{y}} \quad & \theta_{1}(\sum_{(n,v)\in\mathcal{A}}c_{(n,v)}z_{(n,v)} + \sum_{u\in\mathcal{U}_{fast}}m_ur_u)  + opt[ADF(\bm{z})]\\
& \text{Constraints} (\ref{eq:vehicle_flow_appendix}) - (\ref{capacity_charge_appendix}), (\ref{z_domain_appendix}) - (\ref{r_domain_appendix})
\end{align*}

where
\begin{align}
    [\text{ADF}(\mathbf{z})]\quad
    \max_{\bm{\zeta},\bm{\delta}} &\sum_{t\in\mathcal{T}}\sum_{q\in\mathcal{Q}}(\overline{\alpha}_{q}^t - C  \sum_{(n,v)\in\mathcal{A}_{qt}}z_{(n,v)}) \delta_{q}^t+ \Tilde{\alpha}_{q}^t \phi_{q}^t \label{bilinear_obj_appden}\\
    & \phi_{q}^t\leq \delta_{q}^t,&\forall q\in\mathcal{Q}, t\in\mathcal{T}\\
    & \phi_{q}^t\leq M_t\zeta_{q}^t,&\forall q\in\mathcal{Q}, t\in\mathcal{T}\\
    & \phi_{q}^t\geq \delta_{q}^t - M_t(1-\zeta_{q}^t),&\forall q\in\mathcal{Q}, t\in\mathcal{T}\\
     & \phi_{q}^t\geq 0,&\forall q\in\mathcal{Q}, t\in\mathcal{T}\\
    & \text{Constraints (\ref{df1}) - (\ref{df2})}
\end{align}

\newpage
\section{Algorithm combining the Benders decomposition and column generation approaches}\label{sec:BD}
By following \cite{zeng2013solving}, the optimality cuts for the proposed model can be formulated as
\begin{align}
    \eta \geq \sum_{t\in\mathcal{T}}\sum_{q\in\mathcal{Q}}(\overline{\alpha}_{q}^t - C\sum_{a\in\mathcal{A}_{q,t}}\sum_{p\in\mathcal{P}}l_{p}^a x_p) (\delta_{q}^t(w))^*+ \Tilde{\alpha}_{q}^t (\phi_{q}^t(w))^*,\forall w\leq k \label{bd_cuts_zeng}
\end{align}

We decompose the model into a master problem and a sub-problem. The master problem can be formulated as follows:
\begin{align}
\label{mp_model_BD}
{\rm{\big[BD-MP\big]}}\quad
\min\limits_{\mathbf{x}, \mathbf{r}, \mathbf{y}, \mathbf{b}} \quad & \theta_{1}(\sum\limits_{p\in\mathcal{P}}c_p x_p + \sum\limits_{u\in\mathcal{U}_{fast}}m_u r_u) + \eta \\
\mbox{s.t.}\quad\quad
&\eta \geq \sum_{t\in\mathcal{T}}\sum_{q\in\mathcal{Q}}(\overline{\alpha}_{q}^t - C\sum_{a\in\mathcal{A}_{q,t}}\sum_{p\in\mathcal{P}}l_{p}^a x_p) (\delta_{q}^t(w))^*+ \Tilde{\alpha}_{q}^t (\phi_{q}^t(w))^*,& w\leq k,\label{bd_cuts}\\
&\sum_{p\in\mathcal{P}}\beta_{d,p}^{-} x_p=\sum_{p\in\mathcal{P}}\beta_{d,p}^{+} x_p, &\forall d\in\mathcal{D}\label{mp_balance_BD},\\
&\sum_{a\in\mathcal{A}_{q,t}}\sum_{p\in\mathcal{P}}l^a_p x_p \leq N^{com}, & \forall q\in\mathcal{Q},t\in\mathcal{T}\label{mp_max_coposition_BD},\\
&\sum\limits_{p\in\mathcal{P}}g_{u,t,p} x_p\leq \begin{cases}
        r_u, & \text{if } u \in \mathcal{U}_{fast} \\
        N_d^{charge} , & \text{if } u \in \mathcal{U}  \setminus \mathcal{U}_{fast}
    \end{cases} & t\in\mathcal{T},\label{mp_capacity_charge_station_BD}\\
&  x_p\in\{0,1\}, &\forall p\in\mathcal{P}. \label{mp_domain_x_BD}\\
& r_u \in \{0,1,...,N^{charge}_u\} , &\forall u\in\mathcal{U}_{fast}. \label{mp_domain_r_BD}
\end{align}

The subproblem can be formulated as follows:
\begin{align}
    \text{BD-SP}(\hat{\mathbf{x}})=\max_{\bm{\zeta},\bm{\delta}} &\sum_{t\in\mathcal{T}}\sum_{q\in\mathcal{Q}}\big[\overline{\alpha}_{q}^t  \delta_{q}^t+ \Tilde{\alpha}_{q}^t \phi_{q}^t -
C\cdot\delta_{q}^t\sum_{a\in\mathcal{A}_{q,t}}\sum_{p\in\mathcal{P}}l_{p}^a \hat{x}_p\big]\label{bilinear_obj_sp} \\
    & \phi_{q}^t\leq \delta_{q}^t,&\forall q\in\mathcal{Q}, t\in\mathcal{T} \label{df_1} \\
    & \phi_{q}^t\leq M_t\zeta_{q}^t,&\forall q\in\mathcal{Q}, t\in\mathcal{T} \\
    & \phi_{q}^t\geq \delta_{q}^t - M_t(1-\zeta_{q}^t),&\forall q\in\mathcal{Q}, t\in\mathcal{T} \\
     & \phi_{q}^t\geq 0,&\forall q\in\mathcal{Q}, t\in\mathcal{T} \\
     &\zeta_{q}^t\in\{0, 1\},&\forall q\in\mathcal{Q}, t\in\mathcal{T} \\
    &  (\ref{eq:dual1}) - (\ref{eq:dual3}), (\ref{df1}) - (\ref{df_uncertain})\label{df_2} 
\end{align}

Similar to Section \ref{section_cg}, we use the column generation algorithm to solve the master problem. The reduced cost can be calculated as follows:
\begin{align}
    R_p = \theta_{oper} c_p- \sum_{d\in\mathcal{D}}(\beta_{dp}^{-}-\beta_{dp}^{+})\iota_d  - \sum_{q\in\mathcal{Q}}\sum_{t\in\mathcal{T}}\sum_{a\in\mathcal{A}_{q,t}}l_p^a\pi_{q,t} - \sum_{u\in\mathcal{U}}\sum_{t\in\mathcal{T}}g_{u,t,p}\rho_{ut} -C\sum_{w\in\mathcal{W}}\sum_{q\in\mathcal{Q}}\sum_{t\in\mathcal{T}}\delta_{qt}(w)\varepsilon_w\sum_{a\in\mathcal{A}_{q,t}}l_{p}^a. 
\end{align}
wherer $\varepsilon_w$ is the dual variable of constraints (\ref{bd_cuts}).

Below we give the detailed procedure of the algorithm in Algorithm \ref{alg:bd_cg}.

\begin{algorithm}[H]
\caption{Algorithm combining the Benders decomposition and column generation methods}\label{alg:bd_cg}
\SetKwInOut{Input}{Input}
\SetKwInOut{Output}{Output}
\SetKw{KwGoTo}{go to}
\SetKw{KwSet}{set}
\SetKw{KwBreak}{break}
\SetKw{Return}{return}
\SetAlgoLined

\Input{Space-time-SoC network $\mathcal{G} = \{\mathcal{N}, \mathcal{A}\}$; time-dependent and uncertain passenger demand.}
\Output{Optimal solutions $(\mathbf{x}^*, \mathbf{r}^*)$.}
\BlankLine
Initialize $LB = -\infty$, $UB = +\infty$; initialize iteration counter $k = 0$; initialize tolerance $\epsilon$.\\
\While{$(UB - LB)/LB > \epsilon$}{
    Solve the BD-MP, obtain $(\mathbf{x}_{k+1}^*, \mathbf{r}_{k+1}^*, \eta)$;\\
    Update $LB = \theta_1 (\mathbf{c}^{\mathbf{T}} \mathbf{x}_{k+1}^* + \mathbf{m}^{\mathbf{T}} \mathbf{r}_{k+1}^*) + \eta$;\\
    Call GUROBI to solve the $\text{BD-SP}(\mathbf{x}_{k+1}^*)$
    , get the optimal solution $(\delta_{q}^t(k+1))^*$ and $(\phi_{q}^t(k+1))^*$; update $UB = \min \{UB, \theta_1 (\mathbf{c}^{\mathbf{T}} \mathbf{x}_{k+1}^* + \mathbf{m}^{\mathbf{T}} \mathbf{r}_{k+1}^*) + \text{BD-SP}(\mathbf{x}_{k+1}^*)\}$;\\
    Add cuts to the BD-MP:\\
    \Indp
    $\eta \geq \sum\limits_{t\in\mathcal{T}}\sum\limits_{q\in\mathcal{Q}}(\overline{\alpha}_{q}^t - C\sum\limits_{a\in\mathcal{A}_{q,t}}\sum\limits_{p\in\mathcal{P}}l_{p}^a x_p) (\delta_{q}^t(k+1))^*+ \Tilde{\alpha}_{q}^t (\phi_{q}^t(k+1))^*;$\\
    \Indm
    Update $k = k + 1$;\\
}
\Return{Optimal solutions $(\mathbf{x}_{k+1}^*, \mathbf{r}_{k+1}^*)$ and terminate.}
\end{algorithm}

\newpage

\section{Column generation for the model proposed in Section \ref{sec:determinsticModel}}
\label{sec:cg_appendix}

In this section, we introduce the CG algorithm for solving the model proposed in Section \ref{sec:determinsticModel}.

\subsection{Column generation}
The core idea of the column generation algorithm lies in formulating the Restricted Master Problem (RMP) by incorporating only a subset of the feasible paths (denoted as $\hat{\mathcal{P}}$) into the master problem at the beginning. During the iteration process, the pricing problem is solved by solving the RMP, using the resulting dual information, and then the paths with the minimum reduced cost are added to $\mathcal{P}'$. The process is repeated until there are no paths left to add.

\subsubsection{Restricted Master Problem}
The restricted master problem can be formulated as follows:
\begin{eqnarray}\label{rmp_model}
{\rm{\big[RMP\big]}}\quad
\left\{
\begin{array}{ll}
\min\limits_{\mathbf{x}, \mathbf{r}, \mathbf{y},\mathbf{b}} \quad & \theta_{1}(\sum_{p\in\mathcal{P}}c_p x_p + \sum_{u\in\mathcal{U}_{fast}}m_u r_u) + \theta_{2}\sum_{q\in\mathcal{Q}}\sum_{t\in\mathcal{T}} y^t_{q} + \theta_{3}\sum_{q\in\mathcal{Q}}y^{\left|\mathcal{T}\right|}_{q} \\
&\\
&\begin{aligned}
&\sum\limits_{p\in\hat{\mathcal{P}}}\beta_{d,p}^{-} x_p = \sum\limits_{p\in\hat{\mathcal{P}}}\beta_{d,p}^{+} x_p, 
    &\forall d\in\mathcal{D}, \label{cons_balance_apendix}\\[6pt]
    &\sum_{a\in\mathcal{A}_{q,t}}\sum_{p\in \hat{\mathcal{P}}}l^a_p x_p\leq N^{com}, 
    &\forall q\in\mathcal{Q},t\in\mathcal{T},\label{eq:max_coposition_rmp_apendix}\\[6pt]
    &\sum\limits_{p\in\hat{\mathcal{P}}}g_{u,t,p} x_p\leq\begin{cases}
            r_u, & \text{if } u \in \mathcal{U}_{fast}, \\
            N_d^{charge} , & \text{if } u \in \mathcal{U}  \setminus \mathcal{U}_{fast}
        \end{cases},
        &\forall t\in\mathcal{T},\label{eq:capacity_charge_station_rmp_apendix}\\[6pt]
    &y_{q}^t = y_{q}^{t-1} + \alpha_{q}^t - b_{q}^t, 
    &\forall q\in\mathcal{Q}, t\in\mathcal{T},\label{eq:passenger_flow_rmp_apendix} \\[6pt]
    &b_{q}^t  \leq  C \sum_{a\in\mathcal{A}_{q,t}}\sum_{p\in\mathcal{P}}l_{p}^{a} x_p , 
    &\forall q\in\mathcal{Q}, t\in\mathcal{T},\label{eq:capacity_mv_rmp_apendix}\\[6pt]
    &0\leq x_p \leq 1, 
    &\forall p\in\hat{\mathcal{P}}, \label{eq:domain_x_rmp_apendix}\\[6pt]
    &r_u \in [0,N^{charge}_u], 
    &\forall u\in\mathcal{U}_{fast}, \label{eq:domain_r_rmp_apendix}\\[6pt]
    &y_{q}^t\geq 0, b_{q}^t\geq 0, 
    &\forall q\in\mathcal{Q}, t\in\mathcal{T}.\label{eq:domain_y_b_rmp_apendix} 
\end{aligned}
\end{array}
\right. 
\end{eqnarray}

\subsubsection{Pricing problem}
 The reduced cost of path $p$ (denoted as $R_p$) can be calculated as follows:

\begin{align}
    R_p = \theta_{oper} c_p &- \sum_{d\in\mathcal{D}}(\beta_{dp}^{-}-\beta_{dp}^{+})\iota_d - \sum_{q\in\mathcal{Q}}\sum_{t\in\mathcal{T}}\sum_{a\in\mathcal{A}_{q,t}}l_p^a\pi_{qt}- \sum_{u\in\mathcal{U}}\sum_{t\in\mathcal{T}}g_{u,t,p}\rho_{ut} +C\sum_{q\in\mathcal{Q}}\sum_{t\in\mathcal{T}} \varsigma_{qt}\sum_{a\in\mathcal{A}_{q,t}}l^a_p. \tag{35}
\end{align}

To find the path $p\in\mathcal{P}$ with the lowest reduced cost, we solve a time-dependent shortest path problem in the network. We use the label-setting algorithm to solve this problem. Denote the label of each node as $\Psi_n$. Given two nodes, $n$ and $v$, if there exists an arc $a = (n, v) \in \mathcal{A}$, the label update mechanism can be described as follows:
\begin{align}
\Psi_v= \Psi_n + \theta_{oper} * c_{(n,v)} - \beta_{(n,v)}\iota_{d(n,v)} - \pi_{l(n)l(v)t(n)} - \rho_{l(v)t(v)} + C\cdot\varsigma_{l(n)l(v)t(n)}. \tag{36}
\end{align}

\begin{algorithm}[H]
\caption{Customized label-correcting algorithm for the model proposed in Section \ref{sec:determinsticModel}}\label{alg:spfa-specific-nominal}
\SetKwInOut{Input}{Input}
\SetKwInOut{Output}{Output}
\SetKw{KwGoTo}{go to}
\SetKw{KwSet}{set}
\SetKw{KwBreak}{break}
\SetKw{Return}{return}
\SetAlgoLined

\Input{The time-space-SoC network $\mathcal{G}=\{\mathcal{N}, \mathcal{A}\}$, cost functions and dual variables.}
\Output{Paths with minimum reduced costs for each depot pair.}
\BlankLine

Initialize $Q$ to empty queue;\\
Initialize $inQueue[n] = \text{False}$ for all $n \in \mathcal{N}$;\\
\ForEach{node $n \in \mathcal{N}$}{
    \lIf{$n \in \{n^{source}_d | d \in \mathcal{D}\}$}{
        $\Psi(n) = 0$
    }
    \lElse{
        $\Psi(n) = \infty$
    }
    \If{$n \in \{n^{source}_d | d \in \mathcal{D}\}$}{
        Enqueue $n$ into $Q$;\\
        $inQueue[n] = \text{True}$;
    }
}
\While{$Q$ is not empty}{
    $n \gets$ Dequeue from $Q$;\\
    $inQueue[n] = \text{False}$;\\
    \For{each outgoing arc $a=(n, v) \in \mathcal{A}$}{
        Calculate the tentative label for $v$ as follows:\\
        \Indp
        $\Psi_v= \Psi_n + \theta_{oper} * c_{(n,v)} - \beta_{(n,v)}\iota_{d(n,v)} - \pi_{l(n)l(v)t(n)} - \rho_{l(v)t(v)} + C\cdot\varsigma_{l(n)l(v)t(n)}$;\\
        \Indm
        \If{$\Psi_{v}' < \Psi_v$}{
            $\Psi_v = \Psi_{v}'$;\\
            Update the predecessor of $v$ to $n$;\\
            \If{not $inQueue[v]$}{
                Enqueue $v$ into $Q$;\\
                $inQueue[v] = \text{True}$;
            }
        }
    }
}

\For{each depot $d \in \mathcal{D}$}{
    $p_d^* \gets$ reconstruct the path from $n^{source}_d$ to $n^{sink}_d$ using the predecessor links;\\
    \Return{$p_d^*$ and $\Psi_{p_d^*}$};
}

\end{algorithm}

\end{APPENDICES}

\end{document}